\documentclass[a4paper,11pt]{article}
\usepackage{graphicx}        
\usepackage{multicol}        
\usepackage[bottom]{footmisc}
\usepackage{epsfig}
\usepackage{authblk}
\usepackage{graphicx}
\usepackage{float}
\usepackage{diagbox}
\usepackage{nicefrac}
\usepackage[applemac]{inputenc}
\usepackage[T1]{fontenc}
\usepackage{blkarray}
\usepackage{mathtools,amsthm,amsfonts}

\usepackage{newtxtext}
\usepackage[upint,varg]{newtxmath}

\usepackage{color}

\usepackage{bm}
\usepackage{nicefrac}
\usepackage[hmargin={15mm,15mm},vmargin={17mm,22mm}]{geometry}
\usepackage{array}
\usepackage{pdfpages}
\usepackage{multirow}
\usepackage{enumitem}
\usepackage{bbold}
\usepackage{xparse}								
\usepackage{stmaryrd}							
\usepackage{sidecap}								
\usepackage{upgreek}
\usepackage{accents}
\allowdisplaybreaks
\usepackage{tcolorbox}
\usepackage{algorithmic}
\usepackage[ocgcolorlinks,allcolors={blue}]{hyperref}

\usepackage{siunitx}
\usepackage{mhchem}
\usepackage{subcaption}
\usepackage{lineno}

\usepackage[bbgreekl]{mathbbol}

\DeclareSymbolFontAlphabet{\mathbb}{AMSb}
\DeclareSymbolFontAlphabet{\mathbbl}{bbold}

\newcommand{\email}[1]{\href{mailto:#1}{#1}}

\def\[{\Bigl [}
\def\]{\Bigr ]}
\def\({\Bigl (}
\def\){\Bigr )}

\def\<{\langle}
\def\>{\rangle}

\def\bbsig{\bbsigma}
\def\bbeps{\bbespilon}
\newcommand{\jump}[1]{\llbracket #1 \rrbracket}
\def\T{\boldsymbol{T}}

\def\div{{\rm \textbf{div}}}
\def\divScal{{\rm div}}

\def\tensnabla{\bm{\nabla}}

\def\boldnabla{\bm{\nabla}}

\def\bPi{\bm{\Pi}}

\def\cells{\mathcal{M}}
\def\faces{\mathcal{F}}
\def\nodes{\mathcal{V}}
\def\edges{\mathcal{E}}

\def\R{\mathbb R}
\def\ov{\overline}

\def\aa{{\mathfrak a}}

\def\bu{\boldsymbol{u}}
\def\bv{\boldsymbol{v}}
\def\bw{\boldsymbol{w}}

\def\l{{\boldsymbol{\lambda}}}
\def\m{{\boldsymbol{\mu}}}

\def\te{{\tilde e}}

\def\tang{{\boldsymbol{\tau}}}
\def\n{{\boldsymbol{n}}}
\def\d{{\rm d}}

\def\U{\mathbf{U}}
\def\slip{\boldsymbol{S}}
\def\vslip{\dot{\boldsymbol{S}}}
\def\bV{\boldsymbol{V}}

\def\F{f_{\mathrm{fric}}}

\def\K{\mathbb{K}}

\def\Dbb{\mathbb{D}}

\def\be{\boldsymbol{e}}

\def\br{\boldsymbol{r}}

\def\j{\boldsymbol{j}}

\newcommand{\visco}{\mu_f}

\def\x{{\boldsymbol{x}}}

\def\dsp{\displaystyle}

\newcommand{\Ks}{{\mathcal{K}s}}
\def\Ksig{{K\!\sigma}}
\def\Ls{{\mathcal{L}s}}
\def\Lsig{{L\!\sigma}}
\newcommand{\sige}{{\sigma e}}
\newcommand{\MD}{\mathbf{M}_h}

\def\TT{\mathbb{T}}

\newcommand{\Poly}[1]{\mathbb{P}^{#1}}
\newcommand{\Polygras}[1]{\boldsymbol{\mathbb{P}}^{#1}}

\def\wt{\widetilde}

\newtheorem{theorem}{Theorem}
\numberwithin{theorem}{section}
\newtheorem{proposition}[theorem]{Proposition}
\newtheorem{lemma}[theorem]{Lemma}

\theoremstyle{remark}
\newtheorem{remark}[theorem]{Remark}
\theoremstyle{definition}
\newtheorem{assumption}{Assumption}

\begin{document}

\title{VEM-bubble discretisation and iterative coupling algorithm for a fluid induced seismicity model with rate-and-state friction}

\author[1]{{El Hadi Zerga}\footnote{\email{mohammed-El-Hadi.zerga@univ-cotedazur.fr}}}
\author[1]{{Roland Masson}\footnote{\email{roland.masson@inria.fr}}}
\author[2]{{Jean Paul Ampuero}\footnote{\email{ampuero@geoazur.unice.fr}}}
\author[2,3]{{Pierre Romanet}\footnote{\email{romanet@geoazur.unice.fr}}}

\affil[1]{Universit\'e C\^ote d'Azur, Inria, CNRS, LJAD, Team Galets, Nice, France}
\affil[2]{Université Côte d'Azur, Observatoire de la Côte d'Azur, IRD, CNRS, Géoazur, Valbonne, France}
\affil[3]{Cerema, Valbonne, France}

\date{}

\maketitle

\begin{abstract}

This work addresses the discretisation and iterative coupling algorithm for an induced seismicity model accounting for fault reactivation and slip driven by subsurface fluid injection. The model couples single-phase Darcy flow in both the porous matrix and the fault network with a quasi-static poroelastic mechanical response in the bulk and fault slip governed by a rate-and-state friction law with radiation damping stabilisation leading to a quas-dynamic formulation. The contact mechanics problem is discretised using the Virtual Element Method with bubble stabilisation of the facewise constant Lagrange multiplier. It is combined with a Hybrid Finite Volume scheme for the mixed-dimensional Darcy flow model. The time discretisation is implicit and the nonlinear system is solved at each time step using an iterative coupling algorithm which iterates between the Darcy flow, the frictional contact-mechanics, and the state differential equation. The contraction of this iterative coupling algorithm is established with a smallness assumption on the time step independent on the mesh size. This proves convergence of the fixed-point as well as well-posedness of the system. The frictional contact-mechanical model is solved using a semi-smooth Newton algorithm based on a local nonlinear inversion of the friction rate-and-state law and a block Gauss--Seidel preconditioner exploiting the structure of the VEM-bubble-$\mathbb{P}^0$ mixed discretisation. 
\end{abstract}

\section{Introduction}

Subsurface resources such as deep geothermal energy, underground hydrogen storage, and geological CO2 sequestration have become crucial pillars of the energy transition and sustainable development. However, injecting or extracting fluids into or from the subsurface modify the pressure and stress state of the surrounding rocks, often extending far beyond the injection point. These changes can potentially trigger fault reactivation, including earthquakes. Such phenomena pose seismic risks that must be better understood and mitigated to ensure operational safety, reduce risks to nearby populations, and enhance the social acceptance of these projects.
If it remains sufficiently low, this induced seismicity (then referred to as microseismicity) raises no issue and can even contribute to site monitoring. However, the ability to predict and mitigate the risks of induced seismicity is key for sustainable exploitation of the subsurface. In recent years, anthropogenic activities in storage reservoirs have triggered earthquakes with magnitudes up to ~$M_w$ 6 (\cite{keranen2013,kim2018,ellsworth2019}), high enough to cause damage, in regions where natural seismic activity was otherwise low. A recent example is a magnitude 3.9 earthquake in Vendenheim, Alsace, which led to the closure of this deep geothermal site in 2021 (\cite{schmittbuhl2021}).

In this context, numerical simulation stands out as a key tool for better understanding, predicting, and managing these phenomena. It allows to account for coupled multiphysical processes underlying induced seismicity, evaluating the potential impacts of human activities on the subsurface, and designing strategies to minimize associated risks. At sufficiently large spatial scales, the models are based on a representation of the faults as interfaces equipped with mechanical properties (normal stiffness, friction laws, etc.) and hydraulic properties (aperture and hydraulic conductivity). The physical model couples Darcian fluid flow in the porous matrix and the fault network, the poro-elastic deformation of the rock in the matrix domain, as well as the highly nonlinear frictional slip laws along the fault network. An essential feature is the dynamic nature of the friction law, commonly described by a Rate and State law \cite{dieterich1979a,ruina1983} and for which the friction coefficient depends on the slip velocity and a fault state variable.
The decrease of the normal traction  with the fault pressure is the source of the potential fault reactivation which can further accelerate slip and induce seismicity when the friction coefficient and consequently the resistance to slip decreases with the slip velocity.

For induced seismicity, numerical models have evolved from simplified fault reactivation approaches toward fully coupled hydro-mechanical formulations capable of reproducing the transition from fluid injection to aseismic slip and dynamic earthquake rupture. This evolution has been driven by the need to represent increasingly complex physical processes, including fluid diffusion, poroelastic stress transfer, nonlinear frictional behavior, and rupture dynamics, while maintaining computational efficiency for realistic field-scale applications. 

A first distinction among existing approaches concerns the representation of fluid flow. Many models consider fluid transport exclusively within the fault or fracture network, assuming that the surrounding rock acts primarily as an elastic medium. Such reduced-dimensional formulations significantly reduce computational cost while accurately capturing the dominant mechanisms controlling fault pressurisation in low-permeability environments. Boundary element and spectral boundary integral methods have proven particularly effective in this context, enabling efficient simulations of complex fault systems governed by rate-and-state friction \cite{Ampuero24,Dublanchet24,Lapusta2022}. In contrast, mixed-dimensional hydro-mechanical formulations explicitly account for Darcy flow in both the porous matrix and lower-dimensional fault network, allowing pore-pressure diffusion between the reservoir and the fault to be represented consistently. These approaches, generally based on finite element or finite volume discretisations, provide a more realistic description of geothermal reservoirs and subsurface storage applications where matrix diffusion strongly influences fault reactivation \cite{GKT16,Ucar2017,GH19,KEILEGAVLEN2021,Novikov2022,droniou2023bubble,Paluszny2025,geos-2025,ferronato2026}.

A second major distinction concerns the mechanical description of fault slip. Quasi-dynamic formulations neglect elastodynamic wave propagation in the surrounding rock and approximate its dissipative effect on the fault through a radiation damping term, which stabilises the fault dynamics while substantially reducing computational cost and retaining the essential physics of rupture nucleation and slip.  Owing to their efficiency, these models have become the preferred framework for long-term simulations of injection-induced fault slip and earthquake nucleation, particularly when coupled with rate-and-state friction laws \cite{Ampuero24,Dublanchet24,Lapusta2022,Paluszny2025}. Conversely, fully dynamic formulations solve the complete elastodynamic equations and are therefore able to reproduce spontaneous rupture propagation, seismic wave radiation, rupture arrest, and dynamic stress transfer. These methods are essential for investigating the transition from slow slip to seismic rupture and have been implemented using finite elements, extended finite elements, Discontinuous Galerkin methods, finite volumes, and boundary elements \cite{JinZoback18,Pampillon2018,SCHWARTZKOPFF2021,Younis23,AAGabriel2025,Cacace2026,Novikov25}.

The constitutive description of friction represents another important source of variability among numerical models. Slip-weakening laws remain widely used in fully dynamic simulations because of their numerical robustness and their ability to reproduce dynamic rupture propagation \cite{JinZoback18,Younis23,Novikov25}. However, rate-and-state friction has become the reference framework for modeling induced seismicity, as it naturally captures the evolution of fault strength with both slip rate and contact maturity, allowing the simulation of seismic cycles. Note also that, in the context of induced seismicity, rate-and-state allows to model induced aseismic slip, in contrast to slip-weakening. Consequently, most recent quasi-dynamic models rely on rate-and-state friction \cite{Ampuero24,Dublanchet24,Lapusta2022,Paluszny2025}, while several fully dynamic formulations have also successfully incorporated this constitutive law within finite element, XFEM or Discontinuous Galerkin frameworks \cite{Pampillon2018,SCHWARTZKOPFF2021,AAGabriel2025,Cacace2026}.\\

Nevertheless, despite the significant progress achieved over the past decade, numerical methods that simultaneously account for mixed-dimensional fluid flow, advanced friction laws, and fully coupled hydro-mechanical interactions remain relatively scarce. In particular, few scalable formulations incorporate the two-way coupling induced by shear dilation, whereby fault slip modifies hydraulic properties and, consequently, pore-pressure diffusion. This limitation is particularly relevant for underground energy storage applications, such as geological carbon dioxide sequestration and hydrogen storage, where fluid is injected into permeable formations and pressure diffusion within the surrounding porous rock plays a fundamental role in fault reactivation. Motivated by these applications, we propose a fully coupled hydro-mechanical framework that combines mixed-dimensional Darcy flow in the porous matrix and fault with a quasi-dynamic fault mechanics formulation governed by a rate-and-state friction law. A key feature of the proposed model is the strong two-way coupling between mechanics and flow through a shear-dilation contribution to the fault aperture evolution, allowing fault slip to dynamically modify hydraulic transmissivity and, in turn, influence pore-pressure evolution and fault stability.

Following \cite{droniou2023bubble}, the discretisation of the contact mechanical model with rate and state friction, normal stiffness and shear dilation is based on a Lagrange formulation. It combines piecewise constant Lagrange multipliers representing the fault tractions with a low-order Virtual Element Method enriched with face-bubble stabilisation to guarantee inf-sup stability and a diagonal coupling operator between the displacement and the Lagrange multiplier. The mixed-dimensional Darcy flow is discretized using a Hybrid Finite Volume scheme with interface pressure unknowns, allowing pressure discontinuities across matrix-fault interfaces while preserving local mass conservation \cite{BHMS2016}. Unlike most existing numerical approaches for induced seismicity, which rely on explicit integration of the friction evolution, we adopt a fully implicit time discretisation. This choice enables significantly larger time steps governed by physical time scales rather than numerical stability constraints, which is particularly advantageous for strongly coupled hydro-mechanical simulations involving slow pressure diffusion and episodic seismic slip.

Combining ideas from previous work on rupture dynamics models \cite{Pipping2016,DeHoop2020} and mixed-dimensional poromechanics \cite{Girault2016,Kumar2024,Kundan2025}, the resulting nonlinear system is solved at each time step through an iterative coupling algorithm based on a fixed-stress splitting strategy with diagonal approximation of the Schur complement. The fixed-stress type parameters are derived to guarantee convergence of the fixed-point iterations under a time-step restriction independent of the spatial discretisation. Compared with fully monolithic approaches, this partitioned strategy preserves the use of efficient dedicated solvers for each physical subproblem while maintaining robust convergence properties.

A second important contribution concerns the solution strategy for the nonlinear contact mechanics problem arising at each fixed-point iteration. The proposed approach is based on a semi-smooth Newton method, whose key ingredients are the local inversion of the rate-and-state friction law and a block Gauss--Seidel preconditioner for the Jacobian linear system exploiting the block structure of the mixed VEM-bubble--$\mathbb{P}_0$ discretisation. We further prove the nonsingularity of the Jacobian and preconditioner, again under a smallness assumption on the time step. \\

The remainder of this manuscript is organized as follows. Section \ref{sec:model} presents the coupled model, combining a linear mixed-dimensional single-phase Darcy flow formulation with a quasi-dynamic poroelastic model governed by rate-and-state friction. Section \ref{sec:discretisation} recalls the VEM-bubble-$\mathbb{P}^0$ mixed discretisation for contact mechanics in Section \ref{sec:vem} and the Hybrid Finite Volume (HFV) scheme for mixed-dimensional Darcy flow in Section \ref{subsec:HFV}. Section \ref{subsec:coupledmodel} then introduces the space--time discretisation of the fully coupled problem, combining these spatial discretisations with a Crank--Nicolson implicit time integration of the friction and state evolution laws. 
Section \ref{sec:iterativecoupling} presents the iterative coupling algorithm and provides its numerical analysis. In particular, we derive optimal fixed-stress stabilisation coefficients that guarantee fixed-point convergence under a smallness assumption on the time step. 
Section \ref{sec:newton} details the solution strategy for the nonlinear contact mechanics problem arising at each fixed-point iteration. 
Finally, Section \ref{sec:testcases} validates the proposed methodology against induced-seismicity benchmark problems in both two and three spatial dimensions. We then consider a simplified carbon dioxide injection scenario, inspired by \cite{cappa2011}, in both two and three dimensional geometries to assess the performance of the proposed two-way coupled formulation and to verify the theoretical convergence properties of the iterative coupling algorithm.

\section{Quasi-dynamic poromechanical model} \label{sec:model}
\begin{figure}[!htb]
  \begin{center}
	\resizebox{0.28\textwidth}{!}{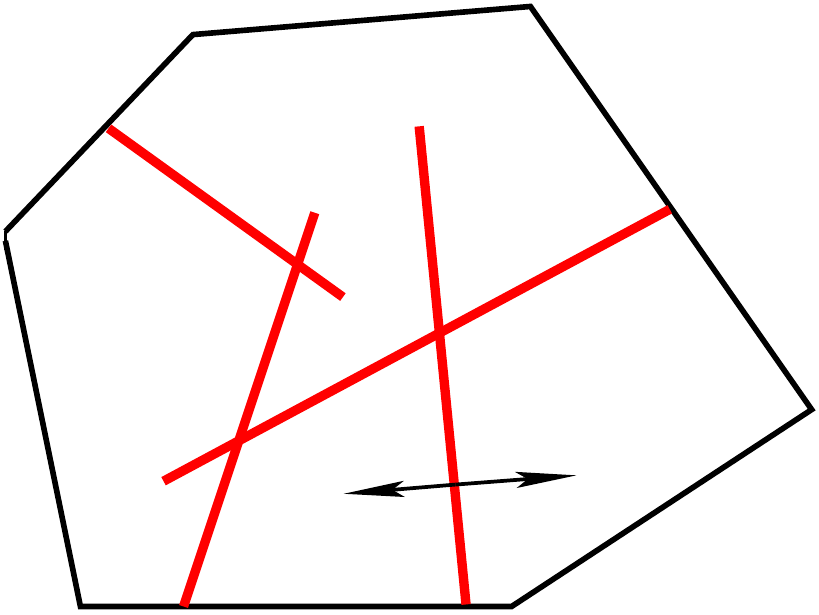}     
        \caption{Model representing the fault network as interfaces with the main variables of the model, namely the displacement field $\boldsymbol{u}$, the pressures in the matrix $p_m$, and in the fault network $p_f$.}
         \label{figdfm}
\end{center}
\end{figure}

We let $\Omega\subset\R^d$, $d\in\{2,3\}$, denote a bounded polytopal domain, partitioned
into a fault domain $\Gamma$ and a matrix domain $\Omega\backslash\overline\Gamma$.
The fault network is defined by 
$
\overline \Gamma = \bigcup_{i\in I} \overline \Gamma_i,
$  
where each fault $\Gamma_i\subset \Omega$, $i\in I$ is a planar polygonal simply connected open domain (see Figure \ref{figdfm}). 
The two sides of a given fault of $\Gamma$ are denoted by $\pm$ in the matrix domain, with unit normal vectors $\n^\pm$ oriented outward from the sides $\pm$. We denote by $\gamma_\aa$ the trace operators on the side $\aa \in \{+,-\}$ of $\Gamma$ for functions in $H^1(\Omega{\setminus}\overline\Gamma)^d$ and by $\gamma^\aa_n$ the normal trace operator for functions in $H_{\div}(\Omega{\setminus}\overline\Gamma)^d$. The jump operator on $\Gamma$ for functions $\bu$ in  $(H^1(\Omega\backslash\overline\Gamma))^d$ is defined by
$$
\jump{\bu} = {\gamma_+\bu - \gamma_-\bu}, 
$$
and we denote by $\jump{\bu}_n = \jump{\bu}\cdot\n^+$ and $\jump{\bu}_\tang = \jump{\bu} -\jump{\bu}_n \n^+$ its normal and tangential components.

The symmetric gradient operator $\bbeps$ is defined such that $\bbeps(\bv) = {1\over 2} (\tensnabla~ \bv + (\tensnabla ~\bv)^t)$ for a given vector field $\bv\in H^1(\Omega\backslash\overline\Gamma)^d$.

Let $\slip = \jump{\bu}_\tang$ denote the slip and $\vslip = \partial_t \slip$ the corresponding slip velocity along the fault network $\Gamma$. The quasi-dynamic model \cite{Uphoff22,Ampuero24,Dublanchet24} combines quasi-static mechanical equilibrium in the bulk with a normal compliance contact model accounting for shear dilation in the normal direction and a Rate and State dynamic Friction (RSF) contact model in the tangential direction:
\begin{equation}
\label{eq_edp_meca} 
\left\{
\begin{array}{lllll}
    &-\div ~ \bbsig^\top(\bu,p_m) & = \boldsymbol{f} & \mbox{ on }  (0,T)\times \Omega{\setminus}\overline\Gamma,\\[1ex]
   &  \bbsig^\top(\bu,p_m) & = \mathbb{A} \bbeps(\bu) - \alpha p_m ~\mathbb{I}   & \mbox{ on }  (0,T)\times \Omega{\setminus}\overline\Gamma,\\[1ex]
  & \T^{+} + \T^{-} &= 0  & \mbox{ on } (0,T)\times \Gamma,\\[1ex]
  &T_\n  & = - K_\n (\jump{\bu}_n + \tan \psi |\jump{\bu - \bu^0}_\tang|)  & \mbox{ on }  (0,T)\times \Gamma, \\[1ex]
  &   \T_\tang & = - \F\( (-T_\n)^+, |\vslip|,\theta \) {\vslip  \over |\vslip|} -\eta \vslip   & \mbox{ on } (0,T)\times \Gamma, 
  \end{array}
\right.
\end{equation}
where $(-T_\n)^+ :=\max(-T_\n,0)$. The coefficient $\eta$ is the radiation damping coefficient, which provides an approximation of the elastodynamic stress changes associated with seismic-wave radiation and plays a key role in controlling the slip rate in the seismic regime. In \eqref{eq_edp_meca}, $\alpha$ denotes the Biot coefficient, $\bbsig^\top$ the total stress, $\mathbb{A}$ the fourth-order elasticity tensor, $K_\n$ the normal stiffness, $\psi$ the shear dilation angle, and $\bu^0$ the initial displacement field.

The contact tractions appearing in \eqref{eq_edp_meca} are defined by
\begin{equation*}
\left\{\!\!\!\!
\begin{array}{lll}
  & \T^{\aa} = \gamma^{\aa}_n \bbsig^{\top} (\bu,p_m)  + p_{f} \n^{ \aa} & \mbox{ on } (0,T)\times \Gamma,  \aa \in \{+,-\},\\[1ex]
  & T_\n = \T^{+}\cdot \n^+, & \mbox{ on } (0,T)\times \Gamma, \\[1ex]
  & \T_\tang = \T^{+} - (\T^{+}\cdot \n^+)\n^+ & \mbox{ on } (0,T)\times \Gamma,
\end{array}
\right.
\end{equation*}
where $p_m$ and $p_f$ denote, respectively, the matrix and fault pressures (see Figure \ref{figdfm}).

The state variable $\theta$ represents the maturity of the fault asperities. It is nondimensionalized by $d_c/v_0$ and expressed on a logarithmic scale, where $d_c$ denotes the characteristic slip distance and $v_0$ the reference slip rate. Its evolution along the fault network $\Gamma$ is governed by the following ordinary differential equation:
\begin{equation}
\label{eq_edo_state} 
 \partial_t \theta + G( \theta)  = H(|\vslip|)  \quad \mbox{ on } (0,T)\times \Gamma. 
\end{equation}

The quasi-dynamic model \eqref{eq_edp_meca}-\eqref{eq_edo_state} is coupled with a mixed-dimensional Darcy flow model describing fluid transport in both the porous matrix and the fault:
\begin{equation}
	\label{eq_edp_darcy} 
	\left\{\!\!\!\!
	\begin{array}{lll}
		&  \partial_t \phi_{m} + c_f \phi_0 \partial_t p_m + \divScal \bV_{m}=h_{m} & \mbox{ on } (0,T)\times \Omega{\setminus}\ov\Gamma,\\[1ex]
		& \bV_{m}=-\frac{ \mathbb{K}_{m}}{\visco} \boldnabla p_{m}  & \mbox{ on } (0,T)\times \Omega{\setminus}\ov\Gamma,\\[1ex]
		& \partial_t \d_{f} + c_f \d_0 \partial_t p_f + \divScal_{\tang} \bV_{f}-\jump{\bV_{m} }_{n}=h_{f}  & \mbox{ on } (0,T)\times \Gamma,\\[1ex]
		& \bV_{f}= - \frac{C_{f}}{\visco} \boldnabla_{\tang} p_{f}, & \mbox{ on } (0,T)\times \Gamma,\\[1ex]
		& \gamma_n^\pm \bV_{m}^{ \pm} = {\Lambda_{f} \over \mu_f} ( \gamma^\pm p_{m} -  p_{f} ) & \mbox{ on } (0,T)\times \Gamma, 
	\end{array}
	\right.
\end{equation}
where the first two equations describe fluid storage and flow in the porous matrix, while the remaining equations account for fluid storage and tangential flow within the fault, together with matrix--fault fluid exchange.

The coupling between the hydraulic and mechanical processes is described by
\begin{equation}\label{coupling_law}
	\left\{\!\!\!\!
	\begin{array}{lll}
	& \partial_t \phi_{m}=\alpha~ \divScal \left(\partial_t \bu \right) + \frac{1}{M} \partial_{t} p_{m}  & \mbox{ on } (0,T)\times \Omega{\setminus}\ov\Gamma,\\[1ex]
	& \partial_t  \d_f = - \partial_t  \jump{\bu }_{\n} & \mbox{ on } (0,T)\times \Gamma.
	\end{array}
\right.
\end{equation}

In \eqref{eq_edp_darcy} and \eqref{coupling_law}, $\visco$ represents the dynamic viscosity of the fluid, $c_f$ its compressibility, $\phi_m$ the porosity of the matrix, $\mathbb{K}_m$ its absolute permeability tensor, and $M$ the Biot modulus. The fault hydraulic aperture is denoted by $\d_f$, its conductivity by $C_f$, and its normal transmissivity by $\Lambda_f$. The quantities $\phi_0$ and $\d_0$ denote, respectively, the reference porosity and reference fault hydraulic aperture. The terms $h_m$ and $h_f$ on the right-hand sides represent fluid injection or production sources in the matrix and fault, respectively.

The model is completed by appropriate pressure boundary and interface conditions. In particular, $p_f$ is typically assumed to be continuous at fault intersections. Boundary conditions are imposed on $p_m$ along $\partial\Omega$ and on $p_f$ along $\partial \Gamma \cap \partial \Omega$. Homogeneous Neumann boundary conditions are classically applied to $p_f$ at the fault tips on $\partial \Gamma \backslash \partial \Omega$.

For later use, we introduce the following notation for the matrix and fault compressibility terms:
\begin{equation}\label{def_compr_matrix_fracture}
\beta_m = c_f \phi_0 + \frac{1}{M} \quad \mbox{ and } \quad \beta_f = c_f \d_0. 
\end{equation}

\begin{assumption}\label{hypothesis}
  We make the following assumption on the model \eqref{eq_edp_meca}-\eqref{eq_edo_state}-\eqref{eq_edp_darcy}-\eqref{coupling_law}.
\begin{itemize}
\item  $\F(\sigma_\n,w,\theta)$ is locally Lipschitz continuous and non negative on $\R^+ \times \R^+ \times \R$.
\item $\F$ is  uniformly Lipschitz continuous w.r.t. $\sigma_\n$ and $\theta$ with Lipschitz constants denoted by $L_\sigma$ and $L_\theta$.
\item $\F$ is non decreasing w.r.t. $w$ on $\R^+$.   
\item  $\F(\sigma_\n,0,\theta) = 0$ for all $\sigma_\n\in \R^+$, $\theta \in \R$.  
\item $G$ is locally Lipschitz continuous and non decreasing on $\R$.
\item $H$ is  uniformly Lipschitz continuous on $\R^+$ with Lipschitz constant denoted by $L_v$.
\item $\beta_m > 0$, $\beta_f > 0$, $K_\n > 0$, $\alpha > 0$ are assumed constant for simplicity.
\item The elasticity tensor $\mathbb{A}$ is fourth order symmetric having the usual uniform ellipticity and boundedness properties. 
\item $C_f$ and $\Lambda_f$ are assumed uniformly bounded  from below by strictly positive constants.  The permeability tensor $\mathbb{K}_{m}$ is second order symmetric having the usual uniform ellipticity and boundedness properties.
\item We assume that the intersection of $\partial\Omega$ with the boundary of each connected component of $\Omega\backslash\overline{\Gamma}$ is of non-zero $(d-1)$-dimensional measure.  
\end{itemize}   
\end{assumption}
  
A typical example of RSF law is provided by the regularized RSF model 
\begin{equation}\label{def_FRSF}
\F\( \sigma_\n, w,\theta \)  = \sigma_\n F_{RSF}(w,\theta) \quad \mbox{ with } \quad F_{RSF}(w,\theta) = a \operatorname{arcsinh}\left[ {w \over 2 v_0} \exp\left( {1 \over a}(f_0 + b \theta) \right) \right], 
\end{equation}
combined with the aging law given by
\begin{equation}\label{def_aginglaw}
  H(w) = - { w\over d_c}, \quad G(\theta) = - {v_0 \over d_c} e^{-\theta}.
\end{equation}  
Note that the aging law satisfies the above assumptions (see \cite{piping2019}). On the other hand, to satisfy the above assumptions, both $\sigma_\n$ and  $F_{RSF}$ must be smoothly capped at prescribed upper bounds to prevent unbounded growth.

\section{Space and time discretisations}\label{sec:discretisation}

We consider a polyhedral mesh of the domain $\Omega$, conforming to the fault network $\Gamma$. For each cell $K$, we denote by $h_K$ its diameter and by $|K|$ its measure; we also denote by $|\sigma|$ the $(d-1)$-dimensional measure of a face $\sigma$. The set of cells $K$, the set of faces $\sigma$,  the set of edges $e$ and the set of nodes $s$ are denoted respectively by  $\cells$, $\faces$, $\edges$ and $\nodes$. 
Let us denote by $\faces_K$ and $\nodes_K$ the set of faces and nodes of a cell $K$, by $\edges_\sigma$ the set of edges of a face $\sigma$, and by $\cells_\sigma$ the set of cells that contain the face $\sigma$. We assume the existence of a subset of faces $\faces_{\Gamma} \subset \faces$ such that
$$
\ov{\Gamma} = \bigcup_{\sigma \in \faces_{\Gamma}} \ov{\sigma}.
$$
The mesh is assumed conforming in the sense that the set $\cells_{\sigma}$ of neighboring cells of $\sigma \in \faces$  is either $\cells_{\sigma} = \{ K,L\}$ for an interior face $\sigma \in \faces^{\text{int}}$ (in which case we write $\sigma=K|L$), or  $\cells_{\sigma} = \{ K\}$ for a boundary face $\sigma \in \faces^{\text{ext}}$.
It is assumed that $\faces_{\Gamma} \subset  \faces^{\text{int}}$ and, if $\sigma=K|L$, that $K$ and $L$ are ordered such that $\n_{\Ksig} = \n^+$ and $\n_{\Lsig} = \n^-$, where $\n_{\Ksig} $ (resp.~$\n_{\Lsig}$) is the unit  normal vector to $\sigma$ oriented outward of $K$ (resp.~$L$). 

We denote by $\nodes^{\text{ext}}$ and $\edges^{\text{ext}}$ the boundary nodes and edges. For each $\sigma \in \faces$,  $\n_{\sige}$ is the unit normal vector to $e \in \edges_{\sigma}$ in the plane spanned by $\sigma$ oriented outward to $\sigma$. \\

The space discretisation of the quasi-dynamic model is based on the VEM-bubble nodal polytopal method introduced in \cite{droniou2023bubble} for quasi-static contact-mechanical models with Coulomb friction.  This scheme is based on a mixed formulation using facewise constant Lagrange multipliers along the fault and the first order nodal polytopal Virtual Element Method (VEM) for the displacement field stabilised by discrete vectorial bubble unknowns on one side of each fault face. This method has the advantage of leading to a diagonal coupling operator for contact, which facilitates its coupling with dynamic friction models on general fault networks. It is briefly recalled in Section \ref{sec:vem} for the reader's convenience. The discretisation of the mixed-dimensional Darcy flow model is based on the Hybrid Finite Volume (HFV) scheme introduced in \cite{BHMS2016} and briefly recalled in Section \ref{subsec:HFV}. Then, Section \ref{sec:quasi-dynamic-scheme} introduces our space time discretisation of the quasi-dynamic poromechanical model \eqref{eq_edp_meca}-\eqref{eq_edo_state}-\eqref{eq_edp_darcy}-\eqref{coupling_law}.

\subsection{VEM-bubble  -- $\mathbb{P}^0$ mixed discretisation of contact-mechanics} \label{sec:vem}
\begin{figure}[!htb]
    \centering
        \includegraphics[width=0.42\textwidth]{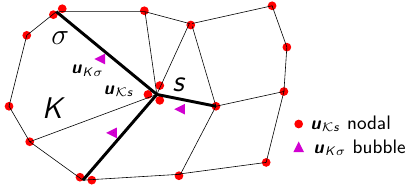}
        \includegraphics[width=0.45\textwidth]{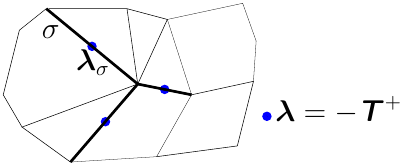}
        \caption{Example of discrete nodal $\boldsymbol{u}_\Ks$ and bubble $\boldsymbol{u}_\Ksig$ displacement  unknowns, and of Lagrange Multiplier unknown $\l_\sigma$.}
        \label{fig-vembulle}
\end{figure}

The degrees of freedom (DOFs) for the displacement are nodal (attached to the vertices of the mesh), but could be discontinuous across the fault network. To be more specific, let us first define a partition $\ov{\cells}_s$ of the set of cells $\cells_s$ around a given node $s\in \nodes$. For a given cell $K\in \cells_s$ we denote by $\Ks\in \ov{\cells}_s$ the subset of $\cells_s$ such that $\bigcup_{L\in \Ks} \ov{L}$ is the closure of the connected component of $(\bigcup_{L\in \cells_s} \ov{L})\setminus \Gamma$ containing the cell $K$. In other words, $\Ks$ is the set of cells in $\cells_s$ that are on the same side of $\Gamma$ as $K$. A nodal displacement unknown $\bv_{\Ks}$ is defined for each $\Ks \in \ov{\cells}_s$. Let us note that there is a unique nodal displacement unknown $\bv_{\Ks}$ at a node $s$ not belonging to $\Gamma$, since $\ov{\cells}_s = \cells_s$ in that case.  On the other hand, for a fault node $s$, the nodal displacement unknown $\bv_{\Ks}$ is the one on the side $K$ of the set of faults connected to $s$. 
Additionally, on the positive side of the fault, we attach to each face a vector DOF corresponding to a ``correction'' of the face values generated by the nodal DOFs; this correction plays the role, at the discrete level represented by the space of DOFs, of a bubble function and is introduced to ensure a suitable inf--sup condition between the space of Lagrange multipliers and the jump of the displacements. For each $K\in\cells$, we denote by
\begin{equation*}
	\faces_{\Gamma,K}^+ = 
	\Big\{ \sigma \in \faces_{\Gamma} \cap \faces_{K}\,:\, \n_{\Ksig} \cdot \n^+  > 0  \Big\}, 
\end{equation*}
the set of faces of $K$ that lie on the positive side of $\Gamma$ (that set is empty if $K$ does not touch $\Gamma$ or only touches it from its negative side). According to the discussion above, the discrete space of displacement is
\begin{equation*}
	\begin{aligned}\U_h = \Big\{\bv_h={}   ((\bv_{\Ks})_{\Ks\in\ov{\cells}_s,\,s\in\nodes},(\bv_{\Ksig})_{K\in\cells,\,\sigma\in\mathcal{F}_{\Gamma,K}^+})\,:
		\bv_\Ks\in\R^d\,,\;\bv_\Ksig\in\R^d \Big\},
	\end{aligned}
\end{equation*}
and its subspace with homogeneous Dirichlet boundary conditions on $\partial \Omega$ is 
\begin{equation}\label{UD_dof}
  \begin{aligned}
  \U_{h,0} = \Big\{\bv_h \in \U_h: \bv_\Ks= \boldsymbol{0} \mbox{ if $s\in\nodes^{\text{ext}}$}\Big\}.
  \end{aligned}
\end{equation}
The Lagrange multiplier plays the role of an approximation of $-\T^+$ on $\Gamma$. Its space is made of piecewise constant vectors:
$$ 
  \MD = \big\{  \l_{h} \in L^2\left(\Gamma \right)^d\,:\, \l_\sigma:=(\l_{h})_{|\sigma}\mbox{ is constant for all $\sigma\in\faces_\Gamma$} \big\}.
  $$
  For $\l_{h} \in \MD$, we define its normal and tangential components by
  $$
  \lambda_{h,n} = \l_{h} \cdot \n^+, \quad \l_{h,\tang} = \l_h - \lambda_{h,n}  \n^+.  
  $$
  Let us refer to Figure \ref{fig-vembulle} for an illustration of the displacement and Lagrange multiplier discrete unknowns. 
We first define, for each $K\in\cells$ and $\sigma\in\faces_K$, a tangential face gradient $\tensnabla^\Ksig:\U_{h,0}\to \Poly{0}(\sigma)^{d\times d}$ and tangential displacement reconstruction $\bPi^{\Ksig}:\U_{h,0}\to \Polygras{1}(\sigma)$. First, we choose nonnegative weights $(\omega_s^\sigma)_{s\in\nodes_\sigma}$ to express the center of mass $\overline{\boldsymbol{x}}_\sigma$ of $\sigma$ in terms of its vertices 
$
\overline{\boldsymbol{x}}_{\sigma} = \sum_{s \in \nodes_{\sigma} } \omega_{s}^{\sigma} \x_s\,,\quad
\sum_{s \in \nodes_{\sigma} } \omega_{s}^{\sigma}=1.
$
Then, for $\bv_h\in\U_{h,0}$, we set
\begin{equation}\label{def:def.operators.sigma}
\begin{aligned}
  \tensnabla^{\Ksig} \bv_{h} ={}& \frac{1}{|\sigma|} \sum_{e= s_1 s_2 \in \edges_{\sigma}} |e| {\boldsymbol{v}_{\Ks_1} + \boldsymbol{v}_{\Ks_2} \over 2} \otimes \n_\sige,\\
  \bPi^{\Ksig} \bv_{h}(\boldsymbol{x})  ={}& \tensnabla^{\Ksig} \bv_{h} (\boldsymbol{x} - \overline{\boldsymbol{x}}_{\sigma}) + \overline{\bv}_{\Ksig}\quad\forall\boldsymbol{x}\in\sigma,\quad\mbox{ where } \quad \overline{\bv}_{\Ksig}=\sum_{s \in \nodes_{\sigma} } \omega_{s}^{\sigma}\bv_{\Ks}.
\end{aligned}
\end{equation}
Above,  we have noted $e=s_1s_2$ to indicate that the edge $e$ has vertices $s_1$, $s_2$. If $\sigma \in \faces_{\Gamma}$ is a fault face, and $K$ (resp.~$L$) is the cell on the positive (resp.~negative) side of $\sigma$, we define the displacement jump operator on $\sigma$ as $\jump{{\cdot}}_\sigma:\U_{h,0}\to\Poly{0}(\sigma)^d$ such that, for all $\bv_h\in\U_{h,0}$,
\begin{equation}\label{eq:def.jump.sigma}
  \jump{\bv_{h}}_{\sigma}  =  \frac{1}{|\sigma|} \int_{\sigma}(\bPi^{\Ksig}  \bv_{h} - \bPi^{\Lsig} \bv_{h}) +  \bv_{\Ksig}
  = \overline{\bv}_\Ksig-\overline{\bv}_\Lsig +  \bv_{\Ksig}.
\end{equation}
The normal component of that jump is denoted by $\jump{{\cdot}}_{\sigma,n} =  \jump{{\cdot}}_{\sigma} \cdot \boldsymbol{n}_{\Ksig}$. For each cell $K\in\cells$, we select nonnegative weights $(\omega_s^K)_{s\in\nodes_K}$ of a linear decomposition of the center of mass $\overline{\boldsymbol{x}}_K$ of $K$ in terms of its vertices
$
  \overline{\boldsymbol{x}}_K=\sum_{s\in\nodes_K}\omega_s^K\boldsymbol{x}_s\,,\quad 
  \sum_{s \in \nodes_K } \omega_{s}^K=1,
$
and we design a gradient reconstruction $\tensnabla^K:\U_{h,0}\to\Poly{0}(K)^{d\times d}$ and a displacement reconstruction $\bPi^K:\U_{h,0}\to\Poly{1}(K)^d$ by setting, for $\bv_h\in\U_{h,0}$,
\begin{align}
  \label{eq:def.nablaK}
  \tensnabla^{K} \bv_{h}  ={}& \frac{1}{|K|} \sum_{\sigma \in \faces_{K}}  |\sigma| \overline{\bv}_{\Ksig} \otimes \n_\Ksig + \sum_{\sigma \in \faces_{\Gamma,K}^+}  \frac{|\sigma|}{|K|} \bv_{\Ksig} \otimes \n_\Ksig,\\
  \label{eq:def.PiK}
  \bPi^{K} \bv_{h}(\boldsymbol{x})  ={}& \tensnabla^{K} \bv_h (\boldsymbol{x} - \overline{\boldsymbol{x}}_{K}) + \overline{\bv}_{K}\quad\forall\boldsymbol{x}\in K\,,\quad
  \mbox{ where } \quad \overline{\bv}_{K} = \sum_{s \in \nodes_K } \omega_{s}^K \bv_{\Ks}.
\end{align}
These local jump, gradient and displacement reconstructions are patched together to create their global piecewise polynomial counterparts $\jump{{\cdot}}_h:\U_{h,0}\to \Poly{0}(\faces_\Gamma)^d$, $\tensnabla^h:\U_{h,0}\to\Poly{0}(\cells)^{d\times d}$ and $\bPi^h:\U_{h,0}\to\Poly{1}(\cells)^d$: for all $\bv_h\in\U_{h,0}$,
$$
(\jump{\bv_h}_{h})_{|\sigma}={}\jump{\bv_h}_\sigma\qquad\forall\sigma\in\faces_\Gamma,\quad 
(\tensnabla^h\bv_h)_{|K}={}\tensnabla^K\bv_h\qquad\forall K\in\cells,\quad 
(\bPi^h\bv_h)_{|K}={}\bPi^K\bv_h\qquad\forall K\in\cells.
$$
We also define the cellwise constant reconstruction operator 
$\wt{\bPi}^{h}\bv_{h}:\U_{h,0}\to\Poly{0}(\cells)^d$ such that  $({\wt{\bPi}^{h}\bv_{h}})_{|K} = \overline{\bv}_K.$
Finally, the discrete symmetric gradient $\bbeps_h(.)$, divergence $\divScal_h(.)$ and stress tensor $\bbsigma_h(.,.)$ are deduced from the previous operators:
$$
\bbeps_{h} =\frac{1}{2}(\tensnabla^{h}+  (\tensnabla^{h})^t),\quad \divScal_{h} = \text{Tr}\left( \bbeps_{h}\right) \quad \text{and}\quad \bbsig_{h}(\bu_h) = \mathbb{A} \bbeps_{h}(\bu_h), 
$$
where $\text{Tr}$ denotes the trace operator. Given the local stabilisation bilinear form $S_{K}:\U_{h,0}\times\U_{h,0} \to\R$ defined by
\begin{equation}\label{eq:def.SK}
	S_{K}(\bu_h,\bv_h)=  h_{K}^{-2}  \sum_{s \in \nodes_K} \left(\bu_{\Ks} -\bPi^{K}\bu_{h} (\boldsymbol{x}_s) \right) \cdot \left( \bv_{\Ks} -\bPi^{K}\bv_{h}  (\boldsymbol{x}_s) \right) 
		 + h_{K}^{-2}\sum_{\sigma \in \faces_{\Gamma,K}^+} \bu_{\Ksig}  \cdot  \bv_{\Ksig}, 
\end{equation}
we denote by $S_{\mathbb{A},h}$ the scaled stabilisation bilinear form defined by
\begin{equation*}
	S_{\mathbb{A},h}( \bu_{h}, \bv_{h}) = \sum_{K\in\cells} |K| \text{A}_K  S_{K}(\bu_h,\bv_h)
\end{equation*}
where $\text{A}_K = {1 \over |K|} \max_{i,j,k,l} \int_K \mathbb{A}_{i,j,k,l} d\x $.

\subsection{Hybrid Finite Volume discretisation of the mixed-dimensional Darcy flow}\label{subsec:HFV}
\begin{figure}[!htb]
  \centering
  \includegraphics[width=0.25\textwidth]{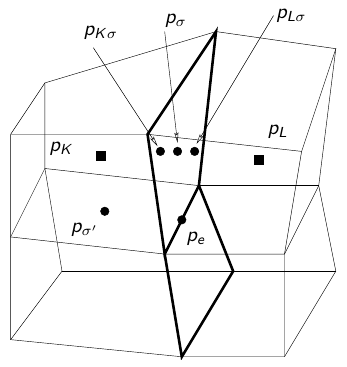}
  \caption{Illustration of the Hybrid Finite Volume (HFV) matrix ($p_K$, $p_{\sigma'}$, $p_{K\sigma}$, $p_L$, $p_{L\sigma}$) and fault ($p_\sigma$, $p_e$) discrete pressure unknowns.}
  \label{fig:HFV_dofs}
\end{figure}

The Hybrid Finite Volume (HFV) discretisation of the mixed-dimensional Darcy flow model is based on the vector space $X_h = X_{h}^m \times X_{h}^f$ of discrete pressures $p_h = (  p_{h}^m, p_{h}^f )$ (see Figure \ref{fig:HFV_dofs}) defined by 
\begin{eqnarray}
X_{h}^m = \left\{p_{h}^m= \((p_K)_{K\in\cells},(p_\sigma)_{\sigma \in \faces\setminus \faces_{\Gamma}},(p_\Ksig)_{\sigma\in\faces_\Gamma,\,K\in\cells_\sigma} \)\,:\,p_K \in \R\,,\; p_{\sigma} \in \R\,,\; p_{\Ksig} \in \R \right\}, \\
X_{h}^f = \left\{ p_{h}^f = \( (p_\sigma)_{\sigma\in \faces_\Gamma}, (p_e)_{e\in \edges_\Gamma} \)\,:\, p_\sigma \in \R\,,\; p_e \in \R \right\}. 
\end{eqnarray}
Assuming for simplicity homogeneous Dirichlet boundary conditions, we denote by $X^{m}_{h,0}$ (resp.~$X^{f}_{h,0}$)  the subspace of $X_{h}^m$ (resp.~$X_{h}^f$) with vanishing values at the boundary $\faces^{\text{ext}}$ (resp.~$\edges^{\text{ext}} \cap \edges_{\Gamma}$), and we set $X_{h,0} = X^m_{h,0}\times X^f_{h,0}$. 
The HFV scheme is obtained by replacing, in the primal variational formulation of \eqref{eq_edp_darcy}, the continuous operators by discrete reconstruction operators $\boldsymbol{\nabla}_{h}^m$, $\Pi_{h}^m$, $\TT_{h}^{m,\pm}$ in the matrix and $\boldsymbol{\nabla}_{h}^f$, $\Pi_{h}^f$  along the faults, defined as follows.

The matrix gradient reconstruction operator $\boldsymbol{\nabla}_{h}^m~:~ X_{h}^m \longrightarrow L^2 ( \Omega )^d$ is such that, for a suitable symmetric positive definite tensor $\TT_{K} = (\TT_{K}^{\nu \nu'})_{(\nu,\nu')\in I_K\times I_K}$, for all $p_{h}^m \in X_{h}^m$,
$$
\int_{K} \mathbb{K}_m  \boldsymbol{\nabla}_{h}^m p_{h}^m \cdot \boldsymbol{\nabla}_{h}^m q_{h}^m =  \sum_{\nu \in I_K} \sum_{\nu' \in I_K} \TT_{K}^{\nu \nu'}\left( p_{\nu} - p_{K}\right) \left( q_{\nu'} - q_{K}\right)\qquad\forall q_{h}^m \in X_{h}^m,
$$ 
with 
$$
I_K = \left\{    \sigma \in \faces_K\setminus \faces_{\Gamma}\right\} \cup  \left\{ \Ksig,~\sigma \in \faces_K \cap \faces_{\Gamma} \right\}, 
$$
denoting the subset of matrix degrees of freedom (DOFs) located at the boundary of the cell $K$. The fault tangential gradient operator $\boldsymbol{\nabla}_{h}^f~:~X_{h}^f \longrightarrow L^2 (\Gamma )^{d-1}$ is such that, for a suitable symmetric positive definite tensor $\TT_{\sigma} = (\TT_{\sigma}^{e e'})_{(e,e')\in \edges_\sigma \times \edges_\sigma }$, for all $p_{h}^f \in X_{h}^f$,
$$
\int_{\sigma} C_{f,h}  \boldsymbol{\nabla}_{h}^f p_{h}^f \cdot \boldsymbol{\nabla}_{h}^f q_{h}^f  =  \sum_{e \in \edges_{\sigma}} \sum_{e' \in \edges_{\sigma}} \TT_{\sigma}^{e e'}\left( p_{e'} - p_{\sigma}\right) \left( q_{e} - q_{\sigma}\right)\qquad \forall q_{h}^f \in X_{h}^f,
$$ 
with the face-wise constant approximation $C_{f,h}$ of the fault conductivity $C_f$. The detailed definitions of $\TT_K$ and $\TT_{\sigma}$ can be found in \cite{BHMS2016}.

The piecewise constant matrix and fault function reconstruction operators $\Pi_{h}^m:X_{h}^m \to   L^2(\Omega )$ and  $\Pi_{h}^f:X_{h}^f \to   L^2 (\Gamma)$ are defined by
$$
\Pi_{h}^m q_{h}^m(\x) = q_K, ~ \forall \x \in K,~K \in \cells,\quad \quad 
\Pi_{h}^f q_{h}^f(\x) = q_{\sigma}, ~ \forall \x \in \sigma,~\sigma \in \faces_{\Gamma}, 
$$
and the face-wise constant sided trace reconstruction operators $\TT^{m,\pm}_{h}:X_{h}  \to L^2 ( \Gamma )$ by  
$$
\TT_{h}^{m,+} q_{h}^m (\x)=  q_{\Ksig}, \quad 
\TT_{h}^{m,-} q_{h}^m(\x)=  q_{\Lsig}, \quad  \forall \x \in \sigma = K|L, \, \sigma\in \faces_\Gamma.    
$$
Let us also denote by  $\Lambda_{f,h}$  the face-wise constant approximation of the fault normal transmissivity.

\subsection{Space time discretisation of the quasi-dynamic poromechanical model} \label{sec:quasi-dynamic-scheme}\label{subsec:coupledmodel}
Let us denote by $\Delta t^n = t^n-t^{n-1} > 0$ the $n$-th the time step and use the superscript $n$ to denote the unknowns at time $t^n$.  As in the continuous case we will use the notations $\slip^n_h = \jump{\bu^n_h}_{h,\tang}$ for the slip and $\vslip_h^n$ for the slip rate at time $t^n$.  The state function is discretised using a facewise constant approximation on $\faces_\Gamma$ and denoted by $\theta^n_h \in M_h$ with
$$ 
  M_h = \big\{  \theta_{h} \in L^2\left(  \Gamma \right)\,:\, \theta_\sigma:=(\theta_{h})_{|\sigma}\mbox{ is constant for all $\sigma\in\faces_\Gamma$} \big\}.
  $$
  Let us define the tangential subspace of $\MD$ 
  $$
  \mathbf{M}_{h,\tang} = \{ \m_h - (\m_h\cdot\n^+) \n^+, \, \forall \, \m_h \in \MD\}, 
  $$
  and let be given the initial conditions
$$
  \bu^0_h \in \U_{h,0}, \quad p^0_h =(p_h^{m,0},p_h^{f,0})\in X_{h,0},  \quad \slip^0_h \in \mathbf{M}_{h,\tang}, \quad \vslip^0_h \in \mathbf{M}_{h,\tang}, \quad  \theta^0_h\in M_h.
 $$
  The scheme discretising the equations \eqref{eq_edp_meca}-\eqref{eq_edo_state}-\eqref{eq_edp_darcy}-\eqref{coupling_law} computes at each time step $n$ the solution
  $$
  \bu_h^n \in \U_{h,0}, \quad \l^n_{h}\in \MD, \quad \vslip_h^n \in \mathbf{M}_{h,\tang}, \quad p^n_h =(p_h^{m,n},p_h^{f,n})\in X_{h,0}, \quad \theta_h^n \in M_h, 
  $$
  of the discrete system \eqref{eq_Crank_Nicolson}-\eqref{Lagrange_meca_VB}-\eqref{GD_DarcyFlow} defined below, given the previous time step solution. The time integration of the quasi-dynamical model is based on the following implicit Crank--Nicolson time integration of the friction and state laws: 
 \begin{equation}
  \label{eq_Crank_Nicolson}
\begin{aligned}
  & {\slip_h^n - \slip_h^{n-1} \over \Delta t^n}   = {1 \over 2} \(\vslip_h^n + \vslip_h^{n-1}\),\\
  & {\theta_h^{n} -\theta_h^{n-1} \over \Delta t^n} = - G( {\theta_h^{n} +\theta_h^{n-1} \over 2}  ) + H(|{1 \over 2} (\vslip_h^n + \vslip_h^{n-1})|). 
 \end{aligned}
 \end{equation}
  It is combined with the VEM-bubble -- $\mathbb{P}^0$ mixed-variational formulation of the contact mechanical problem  \eqref{eq_edp_meca}: 
  \begin{equation}
  \label{Lagrange_meca_VB}
\begin{aligned}
  & \dsp \int_\Omega  \bbsig_h( \bu^n_h): \bbeps_h( \bv_h) +    S_{\mathbb{A},h}\left( \bu^n_{h}, \bv_{h} \right) + \int_\Gamma \l_h^n \cdot \jump{\bv_h}_h \\ 
 &  \quad \quad \quad \quad \quad \quad \quad   \dsp = \int_\Omega  \alpha~ \Pi_h^m p_h^{m,n} ~\divScal_h \bv_h -  \int_\Gamma \Pi_h^f p_h^{f,n} ~\jump{\bv_h}_{h,\n}  + \int_\Omega \boldsymbol{f}\cdot \wt{\bPi}^h\bv_h, \quad \forall \, \bv_{h} \in \U_{h,0},\\[.8em]
  & \dsp   \mbox{with } \quad \lambda^n_{h,\n} = K_\n \( \jump{\bu^n_h}_{h,\n} + \tan\psi |\jump{\bu_h^n - \bu_h^0}_{h,\tang}|\), \quad  \quad \lambda^n_{h,\tang} = \F\((\lambda^n_{h,\n})^+,|\vslip_h^n|,\theta_h^n\){ \vslip_h^n \over |\vslip_h^n|} + \eta \vslip_h^n.  
\end{aligned}
\end{equation}
The discretisation of the Darcy flow model combines the Euler implicit time integration with the HFV spatial discretisation which can be expressed as the following discrete variational formulation: for all $q_{h} = (q_h^{m},q_h^{f}) \in X_{h,0}$, it holds 
\begin{equation}\label{GD_DarcyFlow}
	\begin{array}{lll}
	  & \dsp  \int_{\Omega}  \( \beta_m \Pi_h^m ({p_h^{m,n} - p_h^{m,n-1} \over \Delta t^n}) \, \Pi_{h}^m q_{h}^m \,+\, \alpha~ \divScal_h ({\bu^n_h-\bu_h^{n-1}\over \Delta t^n})  ~\Pi_h^m q_h^{m} 
          \,+\, \frac{\mathbb{K}^{n-1}_m}{\visco} \boldnabla_{h}^m p_{h}^{m,n} \cdot \boldnabla_{h}^m q_{h}^m \)  \\
	  &	\dsp   + \int_{\Gamma}   \( \beta_f \Pi_h^f ({p_h^{f,n} - p_h^{f,n-1} \over \Delta t^n})   \, \Pi_{h}^f q_{h}^f \,-\, \jump{{\bu^n_h-\bu_h^{n-1}\over \Delta t^n} }_{h,\n} ~\Pi_h^f q_h^{f} 
          \,+\, \frac{C_{f,h}^{n-1}}{\visco} \boldnabla_{h}^f p_{h}^{f,n} \cdot \boldnabla_{h}^f q_{h}^f \) \\[2ex]
		&  \dsp + \sum_{\mathfrak{a} \in \{+,-\}} \int_{\Gamma} {\Lambda^{n-1}_{f, h} \over \mu_f} \( \TT_{h}^{m,\mathfrak{a}} p_{h}^{m,n} -  \Pi_{h}^f p^{f,n}_{h} \)   \( \TT_{h}^{m,\mathfrak{a}} q_{h}^m -  \Pi_{h}^f q_{h}^f \)  
	  \dsp   
          = \dsp \int_{\Omega} h_m \Pi_{h}^m q_{h}^m  + \int_{\Gamma} h_f \Pi_{h}^f q_{h}^f.  
	\end{array}
\end{equation}
Note that the matrix permeability and the fault tangential conductivity and normal transmissivity are discretised explicitly in time based typically on the previous time step porosity and fault hydraulic aperture. 

\section{Iterative coupling algorithm}  \label{sec:iterativecoupling}
The system \eqref{eq_Crank_Nicolson}-\eqref{Lagrange_meca_VB}-\eqref{GD_DarcyFlow} is solved at each time step using an iterative coupling algorithm based on a fixed-point scheme solving sequentially at each iteration $k$ the Darcy flow equation at given displacement $\bu_h^{n,k-1}$; the contact-mechanical model at given pressures $p_h^{n,k}$ and state $\theta_h^{n,k-1} $; and the state equation at given slip rate $\vslip^{n,k}$. Note that the current time-step superscript $n$ will be omitted in the following for clarity whenever no ambiguity arises, whereas the index $n-1$ will be indicated explicitly. 

Let us introduce the following notation for the Darcy diffusion operator $\mathcal{D}_h^{m,f}$ such that 
\begin{align*}
  \<\mathcal{D}^{m,f}_h p_h,q_h\> = & \int_{\Omega}\frac{\mathbb{K}^{n-1}_m}{\visco} \boldnabla_{h}^m p_{h}^{m} \cdot \boldnabla_{h}^m q_{h}^m \, +\, \int_{\Gamma} \frac{C_{f,h}^{n-1}}{\visco} \boldnabla_{h}^f p_{h}^{f} \cdot \boldnabla_{h}^f q_{h}^f \\
  &+\, \sum_{\mathfrak{a} \in \{+,-\}} \int_{\Gamma} { \Lambda^{n-1}_{f, h} \over \mu_f} \( \TT_{h}^{m,\mathfrak{a}} p_{h}^{m} -  \Pi_{h}^f p^{f}_{h} \)   \( \TT_{h}^{m,\mathfrak{a}} q_{h}^m -  \Pi_{h}^f q_{h}^f \) \quad \forall (p_h, q_h) \in X_{h,0}\times X_{h,0}.   
\end{align*}
Similarly, let us denote by $\K_h$ the fixed mechanics rigidity operator such that 
\begin{align*}
  \<\K_h \bu_h,\bv_h\> =  \int_\Omega  \bbsig_h( \bu_h): \bbeps_h( \bv_h) +    S_{\mathbb{A},h}\left( \bu_{h}, \bv_{h} \right), \quad \forall (\bu_h,\bv_h) \in \U_{h,0} \times \U_{h,0}, 
\end{align*}
which induces a norm on $\U_{h,0}$ from the discrete Korn inequality on $\U_{h,0}$ (see \cite{DHM25}).

Following \cite{Girault2016,Kumar2024,Kundan2025}, the Darcy flow step includes a generalized fixed-stress stabilisation of the iterative poromechanical coupling with coefficients $C_m^{\rm fs} > 0$ in the matrix and $C_f^{\rm fs}>0$ in the fault network. Let us introduce respectively the matrix cellwise constant and fault facewise constant functions: 
\begin{align*}
  e_{h}^m &= - C_m^{\rm fs}  ~\Pi_h^m p_h^{m} + \alpha~ \divScal_h \bu^n_h,\\
  e_{h}^f &= - C_f^{\rm fs}  ~\Pi_h^f p_h^{f} -  \jump{\bu^n_h}_{h,\n}.     
\end{align*}
The generalized fixed-stress iterative poromechanical coupling algorithm amounts to freeze $e_h^m$ and $e_h^f$ between iterations $k$ and $k-1$. The coefficients $C_m^{\rm fs}$ and $C_f^{\rm fs}$, playing the role of a diagonal approximation of the Schur complement, will be fixed to guarantee the convergence of the fixed-point scheme. \\

\begin{subequations} \label{eq:fixedpoint}
\noindent Darcy flow step at given $e_h^{m,k-1}$, $e_h^{f,k-1}$: compute $p_h^k = (p_h^{m,k}, p_h^{f,k})$ such that for all $q_{h} = (q_h^{m},q_h^{f}) \in X_{h,0}$, it holds 
\begin{equation}
  \begin{aligned}
   & \dsp  \int_{\Omega}  \( (C_m^{\rm fs} + \beta_m) \Pi_h^{m} p_h^{m,k} + e_h^{m,k-1}\) \, \Pi_{h}^m q_{h}^m \,+\, 
    \int_{\Gamma}   \( (C_f^{\rm fs} + \beta_f) \Pi_h^f p_h^{f,k} + e_h^{f,k-1} \) \, \Pi_{h}^f q_{h}^f \,+\,
    \Delta t^n \<\mathcal{D}^{mf}_{h} p^k_h,q_h\>  \\[2ex]
&  = \dsp \int_\Omega \( \beta_m \Pi_h^{m} p_h^{m,n-1} + \alpha ~\divScal_h \bu_h^{n-1} + \Delta t^n h^m \) \,\Pi_{h}^m q_{h}^m +  \int_{\Gamma} \( \beta_f \Pi_h^{f} p_h^{f,n-1} - \jump{\bu_h^{n-1}}_{h,\n} + \Delta t^n  h_f \) \Pi_{h}^f q_{h}^f.  
  \end{aligned}
  \label{darcy_fixedpoint}
\end{equation}
Contact-mechanical step at given $p_h^k$ and $\theta_h^{k-1}$: compute $\bu_h^k$ such that 
\begin{equation}
  \begin{aligned}
    & \dsp \<\K_h \bu^k_h,\bv_h\>  + \int_\Gamma \l_h^k \cdot \jump{\bv_h}_h 
 = \int_\Omega  \alpha~ \Pi_h^m p_h^{m,k} ~\divScal_h \bv_h -  \int_\Gamma \Pi_h^f p_h^{f,k} ~\jump{\bv_h}_{h,\n}  + \int_\Omega \boldsymbol{f}\cdot \wt{\bPi}^h\bv_h,  \quad \forall\, \bv_h\in \U_{h,0}, \\[.8em]
 & \dsp  \mbox{with } \quad  \lambda^k_{h,\n} = K_\n \( \jump{\bu^k_h}_{h,\n} + \tan\psi |\jump{\bu^k_h - \bu_h^0}_{h,\tang}|\), \\
& \dsp  \mbox{and } \quad \lambda^k_{h,\tang} = \F\( (\lambda^{k}_{h,\n})^+,|\vslip_h^k|,\theta_h^{k-1}\){ \vslip_h^k \over |\vslip_h^k|} + \eta \vslip_h^k, \quad \quad  \vslip_h^k   = 2 {\slip_h^k - \slip_h^{n-1} \over \Delta t^n} - \vslip_h^{n-1}. 
  \end{aligned}
  \label{meca_fixedpoint}
\end{equation}
State step at given $\vslip_h^{k}$: compute $\theta_h^k$ such that 
\begin{equation}
  \begin{aligned}
     & {\theta_h^{k} -\theta_h^{n-1} \over \Delta t^n} = - G( {\theta_h^{k} +\theta_h^{n-1} \over 2}  ) + H(|{1\over 2}(\vslip_h^{k} + \vslip_h^{n-1})|). 
  \end{aligned}
  \label{state_fixedpoint}
\end{equation}
\end{subequations}

The convergence analysis of the iterative coupling algorithm \eqref{eq:fixedpoint} will make use of the two constants related to the divergence and normal jump operators as well as to the physical coefficients $\alpha$, $K_\n$  and $\mathbb{A}$: 
\begin{equation}
\label{gamma_mf}
\gamma_m = \max_{\bv_h \in \U_{h,0}} { \| \alpha \divScal_h \bv_h \|^2_{L^2(\Omega)} \over \|\bv_h\|^2_E },\quad 
\gamma_f = \max_{\bv_h \in \U_{h,0}} { \|\jump{\bv_h}_{h,\n} \|^2_{L^2(\Gamma)} \over \|\bv_h \|^2_E }, 
\end{equation}  
with the following elasticity discrete energy scalar product on $\U_{h,0}$ 
$$
\<\bv_h, \bw_h\>_E = \dsp \<\K_h \bv_h,\bw_h\> + \int_\Gamma  K_\n \jump{\bv_h}_{h,\n} \jump{\bw_h}_{h,\n}, 
$$
and associated norm $\|\bv_h\|_E$  such that $\|\bv_h\|_E^2 = \<\bv_h, \bv_h\>_E$. 
Note that we have the usual estimate for homogeneous isotropic elasticity tensor with Lam\'e coefficients $\mu$ and $\lambda$ and a constant Biot coefficient $\alpha$:
$$
\gamma_m \leq {d \alpha^2 \over 2 \mu + d\lambda}.  
$$
Similarly, we clearly have the following estimate for a constant normal stiffness $K_\n$: 
$$
\gamma_f \leq {1 \over K_\n}. 
$$

The following Theorem states the contraction property and the convergence of the iterative  coupling algorithm for a sufficiently small time step $\Delta t^n$.
Note that $\beta_m$ and $\beta_f$  are assumed spatially constant in this analysis. In the numerical experiments in Section \ref{sec:testcoupled2D} and \ref{sec:testcoupled3D}, $\phi_0$ and $\d_0$ will be defined respectively by the previous time step solutions $\phi^{n-1}_h$ and $\d_{f,h}^{n-1}$. It results that $\beta_m$ and $\beta_f$  exhibit in practice spatial variations. A simple but not optimal way to adapt the algorithms to such spatial variations is readily obtained by taking the minimal values of $\beta_m$ and $\beta_f$ in space in the estimates. Let us refer to \cite{Both2017} for a more optimal analysis, based on weighted norms, of the fixed-stress algorithm in the heterogeneous case without fractures. 

\begin{theorem}[Convergence of the iterative coupling algorithm \eqref{eq:fixedpoint}]\label{th_conv_itcoupling}
  For any parameter $ 0 < \rho < 1$, setting
  $$
  C_m^{\rm fs} = {\gamma_m \over \rho}, \quad C_f^{\rm fs} = {\gamma_f \over \rho},
$$
  the fixed-point algorithm \eqref{eq:fixedpoint} is a contraction given by the estimate 
 \begin{equation}\label{meca_theta_darcy_estfinal}
 \begin{aligned}
   & \sum_{i\in \{m,f\}} {1 \over C_i^{\rm fs}} \|\delta^k e_h^i\|^2_{L^2(\Omega_i)} +  2(1-\rho) \|\delta^k \bu_h\|_E^2 + {\eta \over \Delta t^n} \|\delta^k \slip_h\|^2_{L^2(\Gamma)} \\
   & +  {\eta \over  L_v^2 \Delta t^n} \|\delta^k \theta_h\|^2_{L^2(\Gamma)} + 2 C_r \Delta t^n \< \mathcal{D}_h^{mf} \delta^k p_h, \delta^k p_h\> \\
  & \leq 
   2 {\Delta t^n  \over \eta} \( K_\n (L^2_\sigma + \tan^2\psi) \|\delta^{k} \bu_h\|^2_E + K_\n^2 L^2_\sigma \tan^2\psi \|\delta^k \slip_h\|^2_{L^2(\Gamma)} + L^2_\theta  \|\delta^{k-1} \theta_h\|^2_{L^2(\Gamma)} \) \\
    & +  C_r ^2 \sum_{i\in \{m,f\}} {1 \over C_i^{\rm fs}}   \|\delta^{k-1} e_h^i\|^2_{L^2(\Omega_i)}. 
 \end{aligned}  
\end{equation}
  with $C_r =  \dsp \max_{i\in \{m,f\}}{ 1 \over 1 + {\beta_i \over C_i^{\rm fs}}  } < 1$ and any time step $\Delta t^n$ such that 
  \begin{equation}\label{cond_dt}
    \Delta t^n < \min\( { (1-\rho)\eta \over (\tan^2\psi + L_\sigma^2) K_\n }, { \eta \over \sqrt{2} L_v L_\theta}, {\eta \over \sqrt{2} L_\sigma K_\n \tan\psi} \).  
  \end{equation}
  It results that $(\bu_h^k, p_h^k, \theta_h^k)$ converges to the unique solution $(\bu_h^n, p_h^n, \theta_h^n)$ of \eqref{eq_Crank_Nicolson}-\eqref{Lagrange_meca_VB}-\eqref{GD_DarcyFlow}.  
\end{theorem}  

The proof will use the following Lemmae.
\begin{lemma}\label{lemma_friction}
  Let $F_i$, $i\in \{1,2\}$ two continuous functions from $\R^+$ to $\R^+$ such that $F_i(0)=0$. Then for any $\j_i \in \R^d$, we have the estimate
  \begin{equation}
    \(F_1(|\j_1|) {\j_1\over |\j_1|} - F_2(|\j_2|){\j_2\over |\j_2|}\) \cdot (\j_2-\j_1) \leq  \(F_1(|\j_1|) - F_2(|\j_2|)\)  (|\j_2|-|\j_1|). 
  \end{equation} 
\end{lemma}
\begin{proof}[Proof of Lemma \eqref{lemma_friction}]
  We first note that the functions $F_i(\j) {\j \over |\j|}$ are extended continuously by $\boldsymbol{0}$ at $\j = \boldsymbol{0}$. Using $\j_1\cdot \j_2\leq |\j_1||\j_2|$ and $F_i \geq 0$, we have
  \begin{align*}
    \(F_1(|\j_1|) {\j_1\over |\j_1|} - F_2(|\j_2|){\j_2\over |\j_2|}\) \cdot (\j_2-\j_1) & =
    F_1(|\j_1|) \(-|\j_1| + {\j_1\cdot \j_2\over |\j_1|}\) + F_2(|\j_2|) \(-|\j_2| + {\j_1\cdot \j_2\over |\j_2|}\)\\
    & \leq  \(F_1(|\j_1|) - F_2(|\j_2|)\)  (|\j_2|-|\j_1|). 
  \end{align*}
  Note that the cases $\j_1= 0$ and/or $\j_2 = 0$ are easily dealt using $F_i(0)=0$. 
\end{proof}

\begin{lemma}\label{lemma_monotonecontact}
  Let $\boldsymbol{J}_h \in \mathbf{M}_{h,\tang}$, $\theta_h \in M_h$, $\lambda_{h,\n} \in M_h$  be given and let us define the operator $\boldsymbol{B}_h: \U_{h,0} \rightarrow \U_{h,0}$ from the Riesz representation theorem by
  $$
  \<\boldsymbol{B}_h\bv_h, \bw_h\>_E =  \<\bv_h, \bw_h\>_E + { 2 \over \Delta t^n} \int_\Gamma \eta \jump{\bv_h}_{h,\tang}  \cdot \jump{\bw_h}_{h,\tang}  + \int_\Gamma  \F\( (\lambda_{h,\n})^+,|\bV_h|,\theta_h\){ \bV_h \over |\bV_h|}   \cdot \jump{\bw_h}_{h,\tang},  
  $$
  for all $\bw_h \in \U_{h,0}$  with $\bV_h = { 2 \over \Delta t^n} \jump{\bv_h}_{h,\tang} + \boldsymbol{J}_h$.
  Then, the operator $\boldsymbol{B}_h$ is a bijection. 
\end{lemma}
\begin{proof}[Proof of Lemma \eqref{lemma_monotonecontact}]
  Following \cite{chouly-hild-12}, we will show that $\boldsymbol{B}_h$ is a monotone and hemicontinuous operator. Then, it suffices to apply Corollary 15 (p 126) of \cite{Brezis68} to conclude that $\boldsymbol{B}_h$ is a one to one operator. The monotonicity of $\boldsymbol{B}_h$ derives from the  monotonicity property of $\F$  combined with Lemma \ref{lemma_friction} yielding the estimate
  $$
  \<\boldsymbol{B}_h\bv^{(2)}_h - \boldsymbol{B}_h\bv^{(1)}_h, \bv_h^{(2)} -  \bv_h^{(1)} \>_E \geq \|\bv_h^{(2)} -  \bv_h^{(1)}\|_E^2, 
  $$
  for all $(\bv^{(1)}_h, \bv^{(2)}_h) \in \U_{h,0} \times \U_{h,0}$. The hemicontinuity of $\boldsymbol{B}_h$ states that for all $\bv_h$ and $\bw_h$ in  $\U_{h,0}$, the function $\varphi: [0,1] \rightarrow \R$ defined by
  $$
  \varphi(t) = \<\boldsymbol{B}_h( \bv_h - t\bw_h), \bw_h\>_E 
  $$
  is continuous.
  Since $\F\( (\lambda_{h,\n})^+,0,\theta_h\) = 0$ and $\F\( (\lambda_{h,\n})^+,.,\theta_h\)$ is locally Lipschitz w.r.t. its second argument, it results that the function
  $$
  \bV_h \rightarrow \F\( (\lambda_{h,\n})^+,|\bV_h|,\theta_h\) {\bV_h \over |\bV_h|},
  $$
  is locally Lipschitz on $\mathbf{M}_{h,\tang}$ and globally Lipschitz on any bounded ball of $\mathbf{M}_{h,\tang}$, which proves the hemicontinuity of $\boldsymbol{B}_h$. 
\end{proof}

\begin{proof}[Proof of Theorem \ref{th_conv_itcoupling}]
Let us first derive the estimate \eqref{meca_theta_darcy_estfinal}. We denote by $\delta^k r = r^k - r^{k-1}$ the incremental error between iterations $k$ and $k-1$, where $r$ is either a function or a vector of discrete unknowns.

  For $i\in \{m,f\}$ and any $0 < \rho < 1$, let us set $\tilde \gamma_i = {\gamma_i \over \rho}$ and choose $C^{\rm fs}_i = \tilde \gamma_i$.
For parameters $\xi_i > 0$, let us define the functions 
  $$
  \te_h^i = {1 \over \xi_i} e_h^i, \quad i\in \{m,f\}, 
  $$
  and set $\Omega_m = \Omega$ and $\Omega_f = \Gamma$ for convenience. Subtracting \eqref{darcy_fixedpoint} at iteration $k-1$ from \eqref{darcy_fixedpoint} at iteration $k$ with the test vector $\delta^k p_h$, we obtain
  \begin{equation*}
   \sum_{i\in \{m,f\}} (C^{\rm fs}_i +\beta_i) \| \Pi_h^i \delta^k p_h^i\|^2_{L^2(\Omega_i)} +
       \Delta t^n \< \mathcal{D}_h^{mf} \delta^k p_h, \delta^k p_h\> \\
       = - \sum_{i\in \{m,f\}} \int_{\Omega_i} (\delta^{k-1} e_h^i)  \Pi_h^i \delta^k p_h^i. 
  \end{equation*}
Using Young's inequalities with parameters $d_i > 0$ and the definition of $\te_h^i$, we obtain the estimate: 
 \begin{equation}\label{darcy_est}
   \sum_{i\in \{m,f\}} (C^{\rm fs}_i +\beta_i) (2-{1 \over d_i^2}) \| \Pi_h^i \delta^k p_h^i\|^2_{L^2(\Omega_i)} +
       2 \Delta t^n \< \mathcal{D}_h^{mf} \delta^k p_h, \delta^k p_h\> \\
       \leq \sum_{i\in \{m,f\}}  {d_i^2 \xi_i^2 \over C^{\rm fs}_i + \beta_i} \|\delta^{k-1} \te_h^i\|^2_{L^2(\Omega_i)}.   
  \end{equation}
 Using the monotonicity and Lipschitz properties of $\F$, the property $\F(.,0,.) = 0$ combined with Lemma \ref{lemma_friction}, and noticing that $\delta^k \vslip_h = {2 \over \Delta t^n} \delta^k \slip_h$,  we have the following estimates on the tangential contact term
\begin{align*}
   & \dsp - \int_\Gamma \( \F\( (\lambda^{k}_{h,\n})^+,|\vslip_h^k|,\theta_h^{k-1}\){ \vslip_h^k \over |\vslip_h^k|}  -  \F\( (\lambda^{k-1}_{h,\n})^+,|\vslip_h^{k-1}|,\theta_h^{k-2}\){ \vslip_h^{k-1} \over |\vslip_h^{k-1}|} \)
   \cdot \delta^k \slip_h \\
   &\leq - {\Delta t^n \over 2} \int_\Gamma \( \F\( (\lambda^{k}_{h,\n})^+,|\vslip_h^k|,\theta_h^{k-1}\)  -  \F\( (\lambda^{k-1}_{h,\n})^+,|\vslip_h^{k-1}|,\theta_h^{k-2}\) \)
   \cdot \delta^k |\vslip_h| \\
   & \leq   \( L_\sigma K_\n ( \|\delta^{k} \jump{\bu_h}_{h,\n}\|_{L^2(\Gamma)} + \tan\psi \|\delta^k \slip_h\|_{L^2(\Gamma)}) + L_\theta  \|\delta^{k-1} \theta_h\|_{L^2(\Gamma)} \) \|\delta^k \slip_h\|_{L^2(\Gamma)}\\
   & \leq   \( L_\sigma \sqrt{K_\n} \|\delta^{k} \bu_h\|_E + L_\sigma K_\n \tan\psi \|\delta^k \slip_h\|_{L^2(\Gamma)} + L_\theta  \|\delta^{k-1} \theta_h\|_{L^2(\Gamma)} \) \|\delta^k \slip_h\|_{L^2(\Gamma)}. 
\end{align*}
The normal contact term is estimated as follows
\begin{align*}
   & \dsp \int_\Gamma \delta^k \lambda_{h,\n} \jump{\delta^k \bu_h}_{h,\n} = \int_\Gamma  K_\n \( \jump{\delta^k \bu_h}_{h,\n} + \tan\psi \delta^k |\jump{\bu_h - \bu_h^0}_{h,\tang}| \) \jump{\delta^k \bu_h}_{h,\n} \\
   &\geq K_\n \|\delta^{k} \jump{\bu_h}_{h,\n}\|^2_{L^2(\Gamma)} - \tan\psi  \sqrt{K_\n} \|\delta^{k} \bu_h\|_E \|\delta^k \slip_h\|_{L^2(\Gamma)}. 
\end{align*}
 Let us use the short notations $g_h^{m} = \alpha \divScal_h \bu_h$ and $g_h^f = - \jump{\bu_h}_{h,\n}$.  
 Subtracting \eqref{meca_fixedpoint} at iteration $k-1$ from \eqref{meca_fixedpoint} at iteration $k$ with the test vector $\delta^k \bu_h$, and using the previous contact term estimates, we obtain 
 \begin{eqnarray*}
   & \|\delta^k \bu_h\|_E^2 + {2\eta \over \Delta t^n} \|\delta^k \slip_h\|^2_{L^2(\Gamma)}  \\
   & \leq 
   \( (L_\sigma+\tan\psi) \sqrt{K_\n} \|\delta^{k} \bu_h\|_E + L_\sigma K_\n \tan\psi \|\delta^k \slip_h\|_{L^2(\Gamma)} + L_\theta  \|\delta^{k-1} \theta_h\|_{L^2(\Gamma)} \) \|\delta^k \slip_h\|_{L^2(\Gamma)}\\
   & + \dsp \sum_{i\in \{m,f\}} \int_{\Omega_i} (\delta^k g_h^i) \Pi_h^i \delta^k p_h^i.   
 \end{eqnarray*}
 Using Young's inequality, we obtain the estimate
\begin{equation}\label{meca_est}
 \begin{aligned}
   & \|\delta^k \bu_h\|_E^2 + {\eta \over \Delta t^n} \|\delta^k \slip_h\|^2_{L^2(\Gamma)} \\
   & \leq {\Delta t^n \over \eta} \( (L^2_\sigma + \tan^2\psi) K_\n    \|\delta^{k} \bu_h\|^2_E + L^2_\sigma K^2_\n \tan^2\psi \|\delta^k \slip_h\|^2_{L^2(\Gamma)} + L^2_\theta  \|\delta^{k-1} \theta_h\|^2_{L^2(\Gamma)} \) \\
   & + \dsp \sum_{i\in \{m,f\}} \int_{\Omega_i} (\delta^k g_h^i) \Pi_h^i \delta^k p_h^i.   
 \end{aligned}
\end{equation}
Subtracting \eqref{state_fixedpoint} at iteration $k-1$ from \eqref{state_fixedpoint} at iteration $k$, multiplying by the test function $\delta^k \theta_h$, integrating over $\Gamma$, and using the monotonicity property of $G$ and the Lipschitz property of $H$, we obtain the estimate
\begin{equation*}
     {1 \over \Delta t^n} \|\delta^k \theta_h\|^2_{L^2(\Gamma)}  \leq {L_v \over \Delta t^n} \|\delta^k \slip_h\|_{L^2(\Gamma)}   \|\delta^k \theta_h\|_{L^2(\Gamma)}, 
\end{equation*}
from which we deduce that 
\begin{equation}\label{theta_est}
     {\eta \over 2 L_v^2  \Delta t^n} \|\delta^k \theta_h\|^2_{L^2(\Gamma)}  \leq {\eta \over 2\Delta t^n} \|\delta^k \slip_h\|^2_{L^2(\Gamma)}.   
\end{equation}
Combining \eqref{meca_est} and \eqref{theta_est}, we get
\begin{equation}\label{meca_theta_est}
 \begin{aligned}
   & \|\delta^k \bu_h\|_E^2 + {\eta \over 2\Delta t^n} \|\delta^k \slip_h\|^2_{L^2(\Gamma)} + {\eta \over 2 L_v^2 \Delta t^n} \|\delta^k \theta_h\|^2_{L^2(\Gamma)} \\
   & \leq {\Delta t^n \over \eta} \( (L^2_\sigma + \tan^2\psi) K_\n    \|\delta^{k} \bu_h\|^2_E + L^2_\sigma K^2_\n \tan^2\psi \|\delta^k \slip_h\|^2_{L^2(\Gamma)} + L^2_\theta  \|\delta^{k-1} \theta_h\|^2_{L^2(\Gamma)} \) 
     \\  & + \dsp \sum_{i\in \{m,f\}} \int_{\Omega_i} (\delta^k g_h^i) \Pi_h^i \delta^k p_h^i  
 \end{aligned}
\end{equation}
From the definition \eqref{gamma_mf} of $\gamma_i$, and setting $\tilde \gamma_i = C_i^{\rm fs} = \dsp {\gamma_i \over \rho}$, $i\in \{m,f\}$, we have
for a given parameter $c>0$ to be optimized later
$$
\sum_{i\in \{m,f\} } {c \over \tilde\gamma_i}  \|\delta^k g_h^i\|^2_{L^2(\Omega_i)} + 2c(1-\rho) \|\delta^k \bu_h\|_E^2 \leq 2 c \|\delta^k \bu_h\|_E^2. 
$$
Summing $2c$ times \eqref{meca_theta_est} and \eqref{darcy_est}, and using the previous inequality, we obtain the estimate:  
\begin{equation*}
 \begin{aligned}
   & \sum_{i\in \{m,f\}} \( {c \over \tilde\gamma_i} \|\delta^k g_h^i\|^2_{L^2(\Omega_i)} + (\tilde \gamma_i +\beta_i) (2-{1 \over d_i^2}) \| \Pi_h^i \delta^k p_h^i\|^2_{L^2(\Omega_i)} - 2c \int_{\Omega_i} (\delta^k g_h^i) \Pi_h^i \delta^k p_h^i\) \\
   & +  2c(1-\rho) \|\delta^k \bu_h\|_E^2 + 2c{\eta \over 2\Delta t^n} \|\delta^k \slip_h\|^2_{L^2(\Gamma)} + 2c {\eta \over 2 L_v^2 \Delta t^n} \|\delta^k \theta_h\|^2_{L^2(\Gamma)} + 2 \Delta t^n \< \mathcal{D}_h^{mf} \delta^k p_h, \delta^k p_h\> \\
  & \leq 
   2c {\Delta t^n \over \eta} \( (L^2_\sigma + \tan^2\psi) K_\n    \|\delta^{k} \bu_h\|^2_E + L^2_\sigma K^2_\n \tan^2\psi \|\delta^k \slip_h\|^2_{L^2(\Gamma)} + L^2_\theta  \|\delta^{k-1} \theta_h\|^2_{L^2(\Gamma)} \) \\
&     + \sum_{i\in \{m,f\}}  {d_i^2 \xi_i^2 \over \tilde \gamma_i + \beta_i} \|\delta^{k-1} \te_h^i\|^2_{L^2(\Omega_i)}. 
 \end{aligned}
\end{equation*}

Choosing $d_i$ and $\xi_i$ such that 
$$
(\tilde \gamma_i +\beta_i) (2-{1 \over d_i^2}) = c \tilde\gamma_i, \quad \xi_i^2 = {\tilde \gamma_i \over c},  
$$
we have
$$
 {c \over \tilde\gamma_i} \|\delta^k g_h^i\|^2_{L^2(\Omega_i)} + (\tilde \gamma_i +\beta_i) (2-{1 \over d_i^2}) \| \Pi_h^i \delta^k p_h^i\|^2_{L^2(\Omega_i)} - 2c \int_{\Omega_i} (\delta^k g_h^i) \Pi_h^i \delta^k p_h^i = \|\delta^k \te_h^i\|^2_{L^2(\Omega_i)},
 $$
 which gives the estimate

 \begin{equation}\label{meca_theta_darcy_est}
 \begin{aligned}
   & \sum_{i\in \{m,f\}} \|\delta^k \te_h^i\|^2_{L^2(\Omega_i)} +  2c(1-\rho) \|\delta^k \bu_h\|_E^2 + 2c{\eta \over 2\Delta t^n} \|\delta^k \slip_h\|^2_{L^2(\Gamma)} \\
   & + 2c {\eta \over 2 L_v^2 \Delta t^n} \|\delta^k \theta_h\|^2_{L^2(\Gamma)} + 2 \Delta t^n \< \mathcal{D}_h^{mf} \delta^k p_h, \delta^k p_h\> \\
  & \leq 
   2c {\Delta t^n \over \eta} \( (L^2_\sigma + \tan^2\psi) K_\n    \|\delta^{k} \bu_h\|^2_E + L^2_\sigma K^2_\n \tan^2\psi \|\delta^k \slip_h\|^2_{L^2(\Gamma)} + L^2_\theta  \|\delta^{k-1} \theta_h\|^2_{L^2(\Gamma)} \)\\
&     + \sum_{i\in \{m,f\}}  {d_i^2 \xi_i^2 \over \tilde\gamma_i + \beta_i} \|\delta^{k-1} \te_h^i\|^2_{L^2(\Omega_i)}. 
 \end{aligned}
\end{equation}

 It remains to optimize $c>0$ in order to minimize $\dsp \max_{i\in \{m,f\}} {d_i^2 \xi_i^2 \over \tilde\gamma_i + \beta_i}$. We obtain
 $$
 c = \min_{i\in \{m,f\}} {\tilde\gamma_i + \beta_i \over  \tilde\gamma_i }, \quad \max_{i\in \{m,f\}} {d_i^2 \xi_i^2 \over \tilde\gamma_i + \beta_i} = \max_{i\in \{m,f\}} \({\tilde\gamma_i \over \tilde\gamma_i + \beta_i  }\)^2 = {1 \over c^2} < 1.  
 $$
It is easily checked that $0 < d_i \leq 1$ for both $i\in \{m,f\}$.  This proves the contraction property of the fixed-point algorithm \eqref{eq:fixedpoint} under the time step smallness condition  \eqref{cond_dt}.   \\

Let us now prove that, given the previous iterates $(\bu_h^{k-1}, p_h^{k-1}, \theta_h^{k-1})$, each sub-problem  \eqref{darcy_fixedpoint}, \eqref{meca_fixedpoint}, \eqref{state_fixedpoint} has a unique solution under the smallness condition \eqref{cond_dt} on the time step.
For \eqref{darcy_fixedpoint}, existence and uniqueness holds whatever the time step by linearity and from the coercivity of $\mathcal{D}_h^{mf}$.  For \eqref{state_fixedpoint}, existence and uniqueness is a consequence of the monotonicity and Lipschitz properties of $G$ and $H$ (see \cite{piping2019}).  
For \eqref{meca_fixedpoint}, existence and uniqueness is a direct consequence of Lemma \ref{lemma_monotonecontact} provided that a lagged normal traction in the friction law and a lagged shear dilation term in the normal contact term are used. Then, a fixed-point scheme can be defined based on these lagged terms for the contact-mechanical subproblem alone which can be proved to be a global contraction under the smallness condition \eqref{cond_dt} on the time step. The proof follows the same steps as above restricted to the contact-mechanical equations at given pressures and $\theta_h$. This proves existence and uniqueness of the contact-mechanical subproblem  \eqref{meca_fixedpoint} under the smallness time step condition \eqref{cond_dt}.
To conclude the proof of Theorem \ref{th_conv_itcoupling}, the convergence of the fixed-point to a solution of \eqref{eq_Crank_Nicolson}-\eqref{Lagrange_meca_VB}-\eqref{GD_DarcyFlow} and the uniqueness of the solution is classical using global contraction, completeness and continuity arguments.
\end{proof}

\begin{remark} The convergence proof extends when considering an exact integration of the state ODE over the time step at given slip rate rather than the Crank--Nicolson time integration scheme  \eqref{state_fixedpoint}. Let us refer to \cite{piping2019} where an estimate on the state as a function of the slip rate is provided based on the monotonicity and Lipschitz properties of $G$ and $H$. This is in particular well suited for the aging and slip laws which can be integrated analytically at given slip rate. This will be used in the numerical section.   
\end{remark}


\section{Semi-smooth Newton solver for the contact-mechanical sub-problem} \label{sec:newton}

We first remark that thanks to the monotonicity property of $\F$ and to the radiation damping term, the friction law can be inverted such that
  $$
  \l_{\tang} = \F \( \sigma_{\n}, |\vslip|,\theta\) ~{\vslip \over |\vslip| } +\eta \vslip \quad \Leftrightarrow \quad \vslip = \overline V( \sigma_\n, |\l_\tang|,\theta) { \l_{\tang} \over |\l_{\tang}| },
  $$
  where the function $\overline V$ is uniformly Lipschitz from $\R^+ \times \R^+ \times \R$ to $\R^+$.
  The function $\overline V$ is implicitly defined as the solution $w$ of the scalar nonlinear system 
  \begin{equation}\label{eq:solvew}
  |\l_{\tang}| = \F \( \sigma_{\n}, w,\theta\)  +\eta w.  
  \end{equation}
  This nonlinear equation can be efficiently solved on each fault face and for any given $|\l_\tang|\geq 0$, $\sigma_\n \geq 0$ and $\theta \in \R$. If $|\l_\tang| = 0$, we just set $w=0$; else the equation \eqref{eq:solvew} is solved using $\log(w)$  as variable and a dichotomy accelerated by a Newton algorithm.

  A second ingredient is to consider the change of variable for the vector of displacement unknowns defined by 
$$
\left(\begin{array}{c}
  \bv_h^v = (\bv_{\Ks})_{\Ks\in\ov{\cells}_s,\,s\in\nodes}\\
  \jump{\bv_h}_h = (\jump{\bv_h}_\sigma)_{\sigma \in\faces_{\Gamma}}
\end{array}\right)
$$
recalling that $\jump{\bv_{h}}_{\sigma}  = \sum_{s \in \nodes_{\sigma} } \omega_{s}^{\sigma}( \bv_{\Ks} - \bv_{\Ls}) +  \bv_{\Ksig}$ for $\sigma = K|L\in \faces_\Gamma$.
In other words we substitute the fault face bubble unknowns by the fault face jumps. This change of variable decouples the nodal equilibrium equations from the Lagrange multiplier $\l_h$.  
Collecting all source terms in the vector $\boldsymbol{F}_h$, the contact-mechanical step  \eqref{meca_fixedpoint} can be formulated as follows (the fixed-point superscript $k$ is omitted):  
\begin{equation}\label{meca_fixedpoint_V}
  \begin{aligned}
    & \left( \begin{array}{ccc} \K_{vv} &  \K_{vj} &  0
         \\ \K_{jv} & \K_{jj} & \Dbb_{j\l} \end{array}\right)
    \left( \begin{array}{c} \bu^v_h \\ \jump{\bu_h}_h \\ \l_h   \end{array}\right) - \boldsymbol{F}_h = 0, \\
 &    K_\n (\jump{\bu_h}_{\sigma,\n} + \tan\psi |\jump{\bu_h - \bu_h^0}_{h,\tang}| )- \lambda_{\sigma,\n} = 0, \quad \forall \sigma \in \faces_\Gamma, \\
 & \gamma \[\slip_{\sigma} - \slip_{\sigma}^{n-1} - {\Delta t^n \over 2}\( \overline V ( (\lambda_{\sigma,\n})^+, |\l_{\sigma,\tang}|,\theta^{k-1}_\sigma) { \l_{\sigma,\tang} \over |\l_{\sigma,\tang}| } + \vslip_\sigma^{n-1}\)\] = 0, \quad \forall \sigma \in \faces_\Gamma,
\end{aligned}
\end{equation}
with  $(\Dbb_{j\l})_{\sigma,\sigma'} = |\sigma| \delta_{\sigma}^{\sigma'}$, $(\sigma,\sigma') \in \faces_\Gamma \times \faces_\Gamma$, and the block representation of the rigidity matrix $$\K_h =  \left( \begin{array}{cc} \K_{vv} &  \K_{vj} 
  \\ \K_{jv} & \K_{jj}  \end{array}\right).
$$ 
The factor $\gamma >0$ in \eqref{meca_fixedpoint_V} is a rescaling parameter chosen of the order of the shear modulus. 
The nonlinear system \eqref{meca_fixedpoint_V} is solved by a semi-smooth Newton algorithm \cite{Qi1993,QiSun1993} which requires the solution of the following Jacobian system at each iteration 
\begin{equation}\label{jac_syst}
  \begin{aligned}
    & \left( \begin{array}{ccc}
      \K_{vv} &  \K_{vj} &  0 \\
      \K_{jv} & \K_{jj} & \Dbb_{j\l} \\
      0 & \Dbb_{cj} & \Dbb_{c\lambda} \\ 
         \end{array}\right)
    \left( \begin{array}{c} \delta \bu^v_h \\ \jump{\delta \bu_h}_h \\ \delta \l_h   \end{array}\right) = - 
    \left( \begin{array}{c} \boldsymbol{R}^v_h(\bu_h) \\ \boldsymbol{R}^j_h(\bu_h,\l_h)  \\ \boldsymbol{R}^c_h( \jump{\bu_h}_h,\l_h)  \end{array}\right),
\end{aligned}
\end{equation}
where $\boldsymbol{R}^v_h$, $\boldsymbol{R}^j_h$, $\boldsymbol{R}^c_h$ are respectively the nodal, jump and contact components of the residual vector, and $\Dbb_{cj}$ and $\Dbb_{c\l}$ are $d\times d$ block diagonal matrices i.e. local to each fault face. To solve the Jacobian linear system \eqref{jac_syst}, we first eliminate the Lagrange multiplier $\delta \l_h$ from the second equation without any fill-in.
Let us set
$$
\wt \K_{cv} = - \Dbb_{c\lambda} \Dbb_{j\l}^{-1} \K_{jv}, \quad \wt \K_{cj} =  \Dbb_{cj} - \Dbb_{c\lambda} \Dbb_{j\l}^{-1} \K_{jj},
\quad  \wt{\boldsymbol{R}}^{c}_h = \boldsymbol{R}^c_h - \Dbb_{c\lambda} \Dbb_{j\l}^{-1} \boldsymbol{R}^j_h,
$$
the reduced linear system writes 
\begin{equation}\label{jac_red}
  \left( \begin{array}{cc}
      \K_{vv} &  \K_{vj}  \\
      \wt \K_{cv} & \wt \K_{cj}
    \end{array}\right)
    \left( \begin{array}{c} \delta \bu^v_h \\ \jump{\delta \bu_h}_h  \end{array}\right) =
    - \left( \begin{array}{c} \boldsymbol{R}^v_h(\bu_h) \\ \wt{\boldsymbol{R}}^{c}_h  \end{array}\right).
\end{equation}
It is solved using a GMRES iterative solver combined with the block Gauss--Seidel preconditioner 
 \begin{equation}\label{jac_red_prec}
    \left( \begin{array}{cc}
      \mbox{\rm PREC}(\K_{vv}) &  0 \\
      \wt \K_{cv} & \mbox{\rm LU}^{-1}(\wt \K_{cj})
    \end{array}\right),
\end{equation}
 where an LU factorisation is performed on $\wt \K_{cj}$  and  where ${\rm PREC}$ can typically either be obtained from a Cholesky factorisation in dimension $d=2$ or a multigrid preconditioner in dimension $d=2,3$ applied on the fixed SPD matrix $\K_{vv}$.
 Let us remark that the matrix $\wt \K_{cj}$ is block diagonal except at fault intersections for which all fault faces belonging to a given cell are connected.   

\begin{proposition}
  The Jacobian matrix and $\wt \K_{cj}$  are both non-singular provided that the time step satisfies the following smallness condition
  \begin{equation}\label{cond_dt_jac}
\Delta t^n \leq {2 \eta \over K_\n ( (L_\sigma + \tan\psi)^2 + 2 L_\sigma \tan\psi) }.  
\end{equation}
  It results that the semi-smooth Newton convergence theorem applies to \eqref{meca_fixedpoint_V} leading to a local superlinear convergence property assuming a $C^1$ regularity of $\F$ on $[0,+\infty)\times [0,+\infty)\times \R$. 
\end{proposition}

\begin{proof}
In order to simplify the analysis let us write the Jacobian obtained with the friction law expressing the tangential traction as a function of the slip rate. Using the notations
$$
\be_{\sigma,\tang} = {\vslip_\sigma \over |\vslip_\sigma|}, \quad \br_{\sigma,\tang} =  { \jump{\bu_h}_{\sigma,\tang} - \jump{\bu^0_h}_{\sigma,\tang}  \over |\jump{\bu_h}_{\sigma,\tang} - \jump{\bu^0_h}_{\sigma,\tang}|}, 
$$
the differentiation,  at each fault face $\sigma \in \faces_\Gamma$, of the following normal and friction contact equations
\begin{align*}
\lambda_{\sigma,\n} - K_\n \( \jump{\bu_h}_{\sigma,\n} + \tan\psi |\jump{\bu_h - \bu_h^0}_{\sigma,\tang}| \) & = 0,\\
\l_{\sigma,\tang} - \F\( (\lambda_{\sigma,\n})^+,|\vslip_\sigma|,\theta_\sigma^{k-1}\){ \vslip_\sigma \over |\vslip_\sigma|} - \eta \vslip_\sigma &= 0,
\end{align*}  
w.r.t. the normal and tangential jumps and normal and tangential Lagrange multipliers in this order yields the matrix:
$$
\left(\begin{array}{cccc}
  - K_\n & - K_\n \tan\psi ~\br_{\sigma,\tang}^T & 1 & 0 \\
  0 & - A_\sigma & -(\partial_{\sigma_\n} \F) \mathcal{H}(\lambda_{\sigma,\n})\be_{\sigma,\tang} & \mathbb{I}_{d-1}
 \end{array} 
\right),  
$$
where $\mathcal{H}$ is the Heaviside function, and with the SPD matrix 
$$
A_\sigma = {2 \over \Delta t^n} \( (\partial_w \F) \be_{\sigma,\tang}\otimes \be_{\sigma,\tang} + {\F \over |\vslip_\sigma| } (\mathbb{I}_{d-1} - \be_{\sigma,\tang}\otimes \be_{\sigma,\tang} ) + \eta \mathbb{I}_{d-1} \) 
$$
such that
\begin{equation}\label{coercivity:Asigma}
\< A_\sigma \xi,\xi\> \geq \eta {2 \over \Delta t^n} |\xi|^2 \quad \forall \xi \in \R^{d-1}. 
\end{equation}
Note that $\Dbb_{cj}$ and $\Dbb_{c\lambda}$ are the block diagonal matrices defined by 
$$
(\Dbb_{cj})_{\sigma,\sigma} = \left(\begin{array}{cc}
  K_\n & K_\n \tan\psi ~\br_{\sigma,\tang}^T  \\
  0 & \gamma \mathbb{I}_{d-1} 
 \end{array} 
\right), \quad (\Dbb_{c\lambda})_{\sigma,\sigma} = - \left(\begin{array}{cc}
  1 & 0 \\
  - \gamma (\partial_{\sigma_\n} \F) \mathcal{H}(\lambda_{\sigma,\n}) A_\sigma^{-1} \be_{\sigma,\tang} &  \gamma A_\sigma^{-1}
 \end{array} 
\right).   
$$
The elimination of the Lagrange multipliers w.r.t. the jumps from the contact equations and their substitution in the contact term yields the bilinear form
\begin{align*}
  \delta\l_\sigma \cdot \jump{\delta\bv_h}_{\sigma}  = & \< A_\sigma \jump{\delta\bu_h}_{\sigma,\tang}, \jump{\delta\bv_h}_{\sigma,\tang}\> + K_\n \jump{\delta\bu_h}_{\sigma,\n} \jump{\delta\bv_h}_{\sigma,\n} \\
  &  + K_\n \tan\psi  \<\br_{\sigma,\tang}, \jump{\delta\bu_h}_{\sigma,\tang}\> \jump{\delta\bv_h}_{\sigma,\n} \\
  & +(\partial_{\sigma_\n} \F)  \mathcal{H}(\lambda_{\sigma,\n}) K_\n \(\jump{\delta\bu_h}_{\sigma,\n}  + \tan\psi  \<\br_{\sigma,\tang}, \jump{\delta\bu_h}_{\sigma,\tang}\>\) \<\be_{\sigma,\tang}, \jump{\delta\bv_h}_{\sigma,\tang}\>  
\end{align*}
From Cauchy Schwarz inequality, the Lipschitz property of $\F$ and \eqref{coercivity:Asigma}, it results that
$$
\delta\l_\sigma \cdot \jump{\delta\bu_h}_{\sigma}  \geq {2 \eta \over \Delta t^n} |\jump{\delta\bu_h}_{\sigma,\tang} |^2 + K_\n |\jump{\delta\bu_h}_{\sigma,\n} |^2 - L_\sigma \tan\psi K_\n |\jump{\delta\bu_h}_{\sigma,\tang} |^2 - K_\n (L_\sigma + \tan\psi) |\jump{\delta\bu_h}_{\sigma,\tang} ||\jump{\delta\bu_h}_{\sigma,\n} |. 
$$
Then, from Young's inequality and the smallness assumption on the time step \eqref{cond_dt_jac}, we obtain the estimate  
\begin{align*}
  \delta\l_\sigma \cdot \jump{\delta\bu_h}_{\sigma}   &\geq \( {2\eta \over \Delta t^n} -  K_\n ( {1 \over 2} (L_\sigma + \tan\psi)^2 +L_\sigma \tan\psi) \) |\jump{\delta\bu_h}_{\sigma,\tang} |^2 + {K_\n \over 2} |\jump{\delta\bu_h}_{\sigma,\n} |^2 \\
  & \geq {\eta \over \Delta t^n} |\jump{\delta\bu_h}_{\sigma,\tang} |^2 + {K_\n \over 2} |\jump{\delta\bu_h}_{\sigma,\n} |^2 
\end{align*}
This shows that the contact term obtained after eliminating the traction from the contact equations yields a coercive contribution under the smallness condition on the time step \eqref{cond_dt_jac}. Since $A_\sigma$ is non singular, this coercivity property of course still holds when starting from the inverse friction formulation. Together with the SPD property of the rigidity matrix $\K_h$, this proves that the full Jacobian matrix in \eqref{jac_syst} and the reduced Jacobian matrix in \eqref{jac_red_prec} are non singular under the condition \eqref{cond_dt_jac}. 
It also implies that
$$
\K_{jj} - \Dbb_{j\l} \Dbb_{c\lambda}^{-1}  \Dbb_{cj},
$$
is non singular. Since the block diagonal matrix $\Dbb_{c\lambda}$  is non singular, this proves that
$$
\wt \K_{cj} =  \Dbb_{cj} - \Dbb_{c\lambda} \Dbb_{j\l}^{-1} \K_{jj},
$$
is also non singular under the condition \eqref{cond_dt_jac}.

To establish the local superlinear convergence of the semi-smooth Newton algorithm applied to \eqref{meca_fixedpoint_V}, it remains to prove the semi-smoothness of the nonlinear terms
$$
\left|\jump{\bu_h-\bu_h^0}_{h,\tang}\right|
\quad\text{and}\quad
\overline V\bigl((\lambda_{\sigma,\n})^+,|\l_{\sigma,\tang}|,\theta^{k-1}_\sigma\bigr)
\frac{\l_{\sigma,\tang}}{|\l_{\sigma,\tang}|}
$$
in system \eqref{meca_fixedpoint_V}. The semi-smoothness of the first term is a classical result.

For the second term, the $C^1$-regularity and monotonicity assumptions on $\F$ together with $\eta > 0$ and the property $\F(\sigma_\n,0,\theta^{k-1}_\sigma) =0$ imply that the mapping
$$
(\sigma_\n,\mathbf t)\longmapsto G(\sigma_\n,\mathbf{t}) = 
\overline V(\sigma_\n,|\mathbf{t}|,\theta^{k-1}_\sigma)
\frac{\mathbf{t}}{|\mathbf{t}|}
$$
is $C^1$ on $[0,+\infty)\times\mathbb R^2$, with the values at $\mathbf{t}=\mathbf{0}$ understood by continuous extensions with $G(\sigma_\n,\mathbf{0}) = \mathbf{0}$ and
  $DG(\sigma_\n,\mathbf{0}) = \left( \mathbf{0}, {1 \over \partial_w \F(\sigma_\n,0,\theta^{k-1}_\sigma) + \eta} \mathbb{I}_{d-1} \right)$. Since the positive-part mapping $\sigma_\n\mapsto \sigma_\n^+$ is semi-smooth, the composition
$$
\overline V\bigl((\lambda_{\sigma,\n})^+,|\l_{\sigma,\tang}|,\theta^{k-1}_\sigma\bigr)
\frac{\l_{\sigma,\tang}}{|\l_{\sigma,\tang}|}
$$
is also semi-smooth. This establishes the semi-smoothness of both nonlinear terms.
      
\end{proof}   

\begin{remark}\label{remark_elimination_traction}
We could alternatively eliminate the Lagrange multiplier directly from the Jacobian linear system using the contact law. The elimination strategy adopted here is motivated by its greater flexibility, as it remains applicable to contact laws for which the traction block $\Dbb_{c\lambda}$ is singular. This situation arises, for example, when enforcing a non-penetration condition corresponding to an infinite normal stiffness.
\end{remark}

\section{Numerical experiments}  \label{sec:testcases}
The following test cases all consider the regularised RSF law defined by \eqref{def_FRSF} combined with the aging law \eqref{def_aginglaw}. 
The adaptive time stepping is defined by an initial time step $\Delta t^{0}$, a maximum time step $\Delta t_{max}$ and a parameter $M_{\Delta t}$ chosen in the interval $[1,16]$ in order to assess the numerical behavior of the algorithms with respect to various time step sizes. It is defined by  
\begin{equation}
  \label{eq-timestepping}
\Delta t^n = \min\(1.2 \Delta t^{n-1}, {d_c \over M_{\Delta t} \|\vslip\|_\infty },\Delta t_{max}\). 
\end{equation}
Note that the time step sufficient condition \eqref{cond_dt} for the convergence of the iterative coupling algorithm is obtained based on global Lipschitz constants and is consequently not expected to be sharp especially during slow or aseismic slip. 
For fully coupled hydro-mechanical test cases, the matrix stabilisation coefficient is defined by $C_m^{\rm fs} = \gamma_m$ (corresponding to $\rho=1$) with $\gamma_m$ estimated by ${d \alpha^2 \over 2 \mu + d\lambda}$ where $d$ is the space dimension. The fault stabilisation coefficient is defined by $C_f^{\rm fs} = \gamma_f$ with $\gamma_f$ estimated by ${1 \over K_\n}$. 
To improve the convergence of the iterative coupling algorithm in the case with shear dilation, we also test an adaptive choice of the stabilisation coefficient $C_f^{\rm fs}$ based on probing. In that case $C_f^{\rm fs}$ depends on the fault face and is computed in order to be a diagonal approximation of the linearised Schur complement operator expressing the normal jump as a function of the fault pressure. This approximation is computed in order to be exact on constant fault pressures and require the solution of a linearised contact-mechanical system (see \cite{Jonval2026} for details in the case of static Coulomb friction). In the following experiments it is updated only at each time step. Note that, in the case without shear dilation, the probing stabilisation is roughly constant along the fault and equal to ${1 \over K_\n}$. Consequently, it does not provide any significant improvement and is not considered in the following tests. \\

The convergence of the fixed-point iterations is reached if the weighted sum of the state, pressure and displacement increment norms is below $10^{-4}$. For the contact-mechanics nonlinear solver, the stopping criterion is achieved if the relative $l^2$ norm of the residual or the weighted sum of the displacement and traction increment maximum norms is below $10^{-4}$.

The linear solvers are implemented using the library PETSc \cite{balay2018petsc} both for the Darcy flow and the contact-mechanics.
The computations were performed on a machine equipped with a single NVIDIA RTX 3090 GPU (24 GB), a 6-core CPU, and 30 GB of RAM. The linear systems were solved on the GPU using PETSc preconditioning. For the nodal elasticity block, the GMRES solver was preconditioned with either a Cholesky factorisation in the 2D tests or the PCMG geometric multigrid method in the 3D tests, while for the Darcy problem, GMRES was combined with the BoomerAMG HYPRE preconditioner.

\subsection{Validation of the quasi-dynamic discretisation}\label{sec:testvalidation}

\subsubsection{Community induced  seismicity 2D benchmark}

\begin{figure}[!htb]
  \centering
  \begin{center}
    \includegraphics[width=0.5\textwidth]{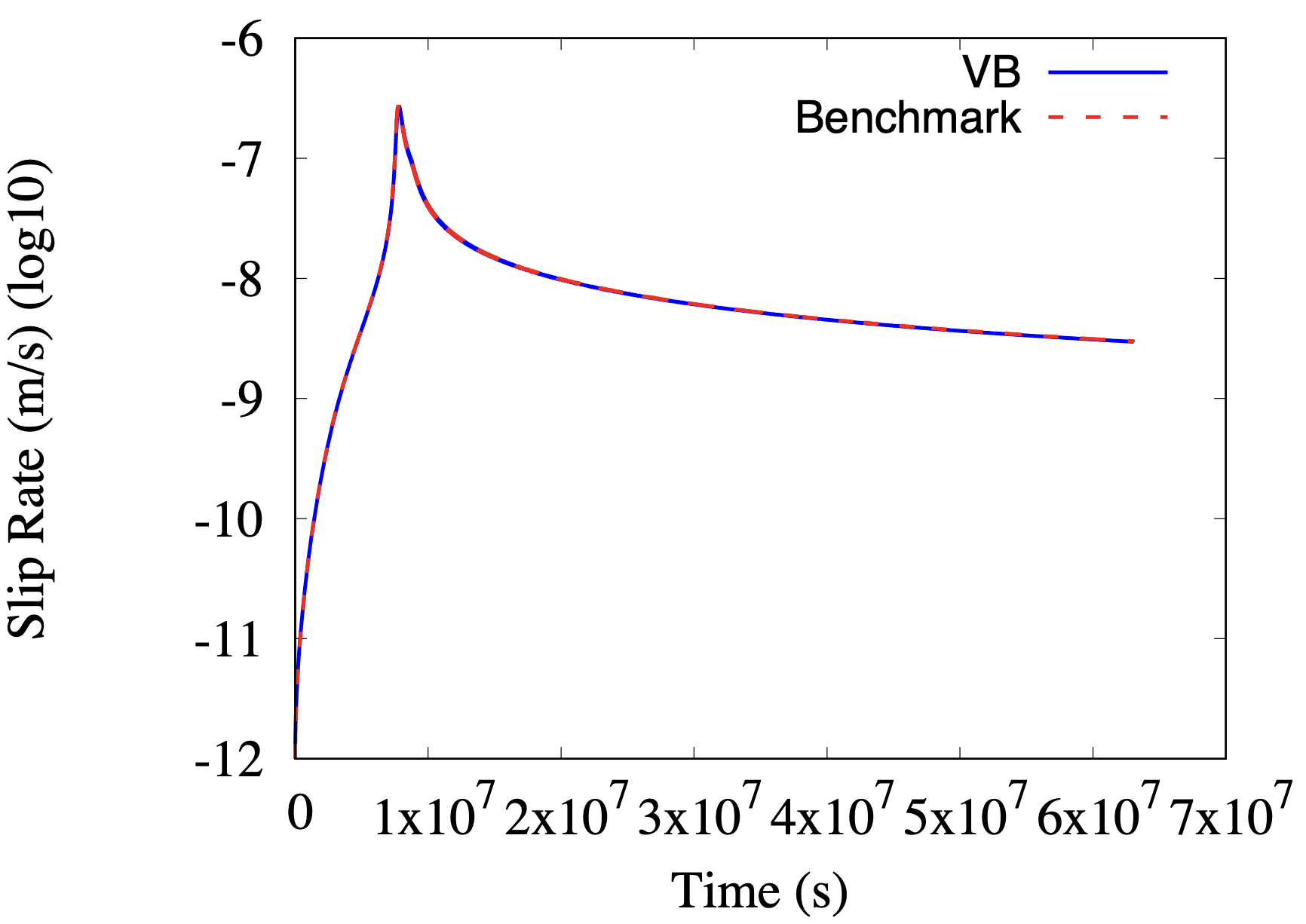}
     \includegraphics[width=0.4\textwidth]{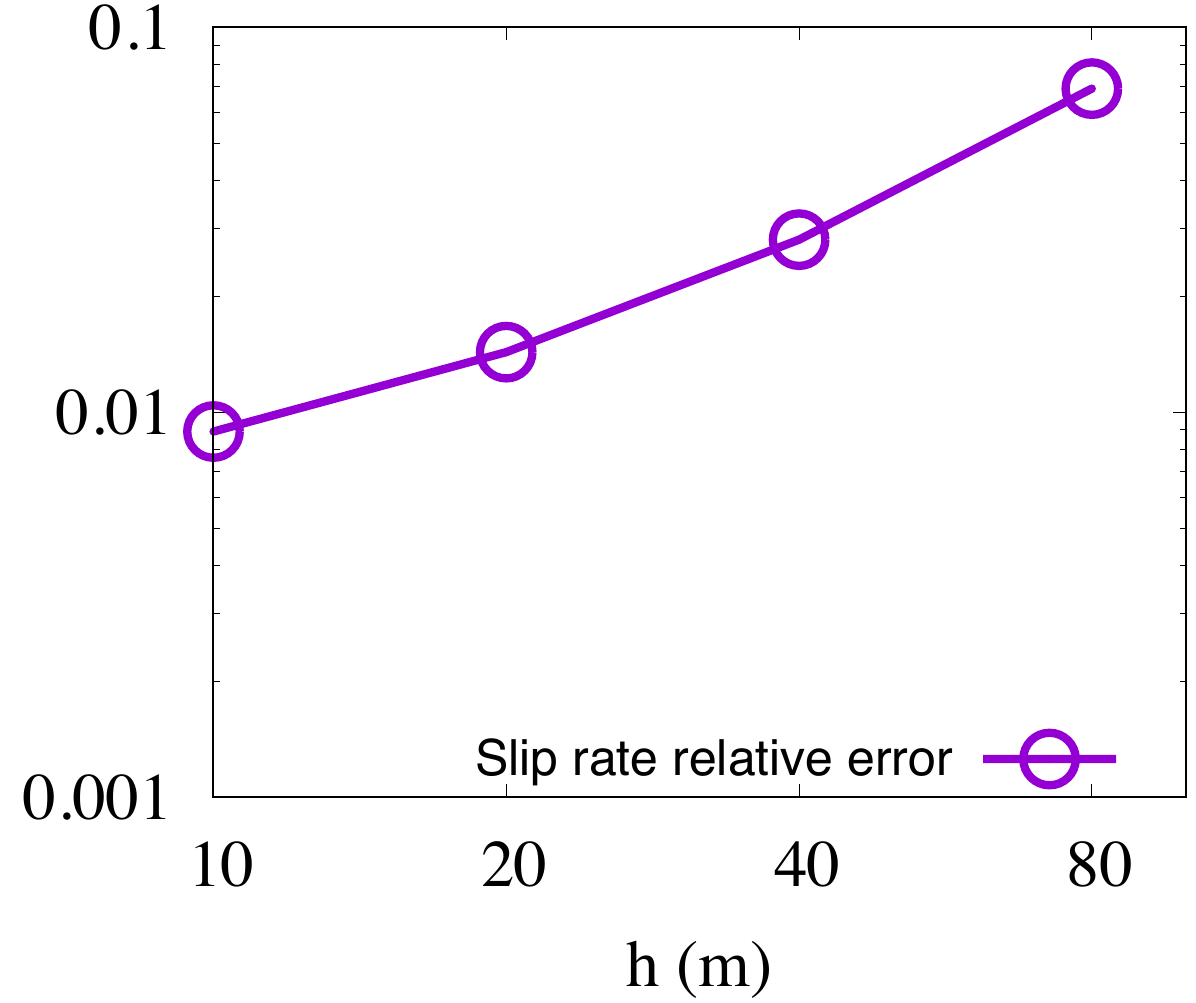}
  \caption{(left) Comparison of the maximum slip rate in log10 scale as a function of time obtained with our scheme (VB) and with the FastCycles benchmark solution on a mesh of step size $10$ m along the fault. (right) Convergence of the relative error for the maximum slip rate evaluated w.r.t. the fine reference benchmark solution on a family of meshes of size $h$ along the fault.}
  \label{fig-slip}
\end{center}
\end{figure}

The numerical model is first validated using the induced seismicity community benchmark SEAS BP6-QD-A \cite{bench2025} which considers a 2D antiplane problem for a planar fault embedded in a homogeneous unbounded space with rate and state friction and aging state laws. The aseismic fault slip is induced by the injection of fluid during $100$ days in the middle of the fault. Fluid diffusion occurs only in the fault (impervious matrix) at fixed hydraulic aperture $\d_f$ and fault conductivity $C_f$ leading to a one way coupled model. Figure \ref{fig-slip} compares the maximum slip rate as a function of time obtained with our scheme (VB) to a reference solution from the benchmark obtained with the Spectral Boundary Element code FastCycles \cite{romanet2021} on a mesh of step size 10 m along the fault. It exhibits a very good match on the finest mesh of step size $10$ m along the fault and roughly a first order convergence on a family of refined meshes as expected.

\subsubsection{Two intersecting faults 2D benchmark}
To validate the spatial and temporal discretisations of the quasi-dynamic model with intersecting faults, we consider the two-dimensional benchmark problem introduced in \cite{Ampuero24} (Appendix C). The benchmark consists of two intersecting faults, each $10$ m long, forming an angle of $30^\circ$. After an initial period of three years of regional tectonic loading at a low rate, fluid is injected at a constant rate for $2000$ s. For this one-way coupled test case, the surrounding rock matrix is assumed to be impermeable and infinite, while the hydraulic properties of the faults remain constant throughout the simulation. The complete set of model parameters, including the fault rate-and-state friction parameters and the fault strengthening rate, is provided in \cite{Ampuero24} (Table 4).

\begin{figure}[H]
    \centering
        \includegraphics[width=0.35\textwidth]{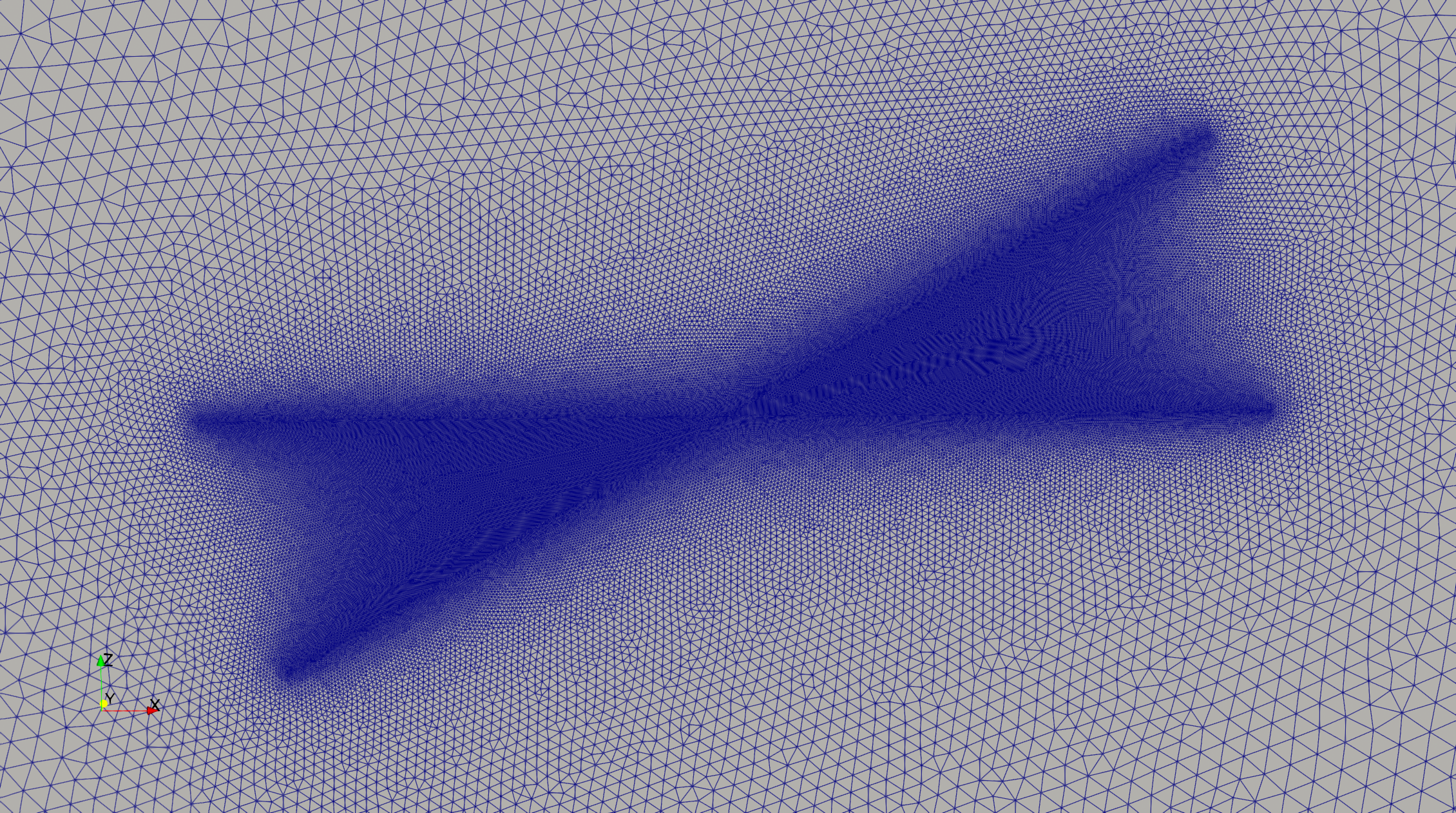}
        \caption{Zoom on the mesh with $1600$ fault edges and refined in the neighborhood of the two intersecting faults.}
        \label{fig-mesh-2failles}
\end{figure}

\begin{figure}[H]
    \centering
    \includegraphics[width=0.45\textwidth]{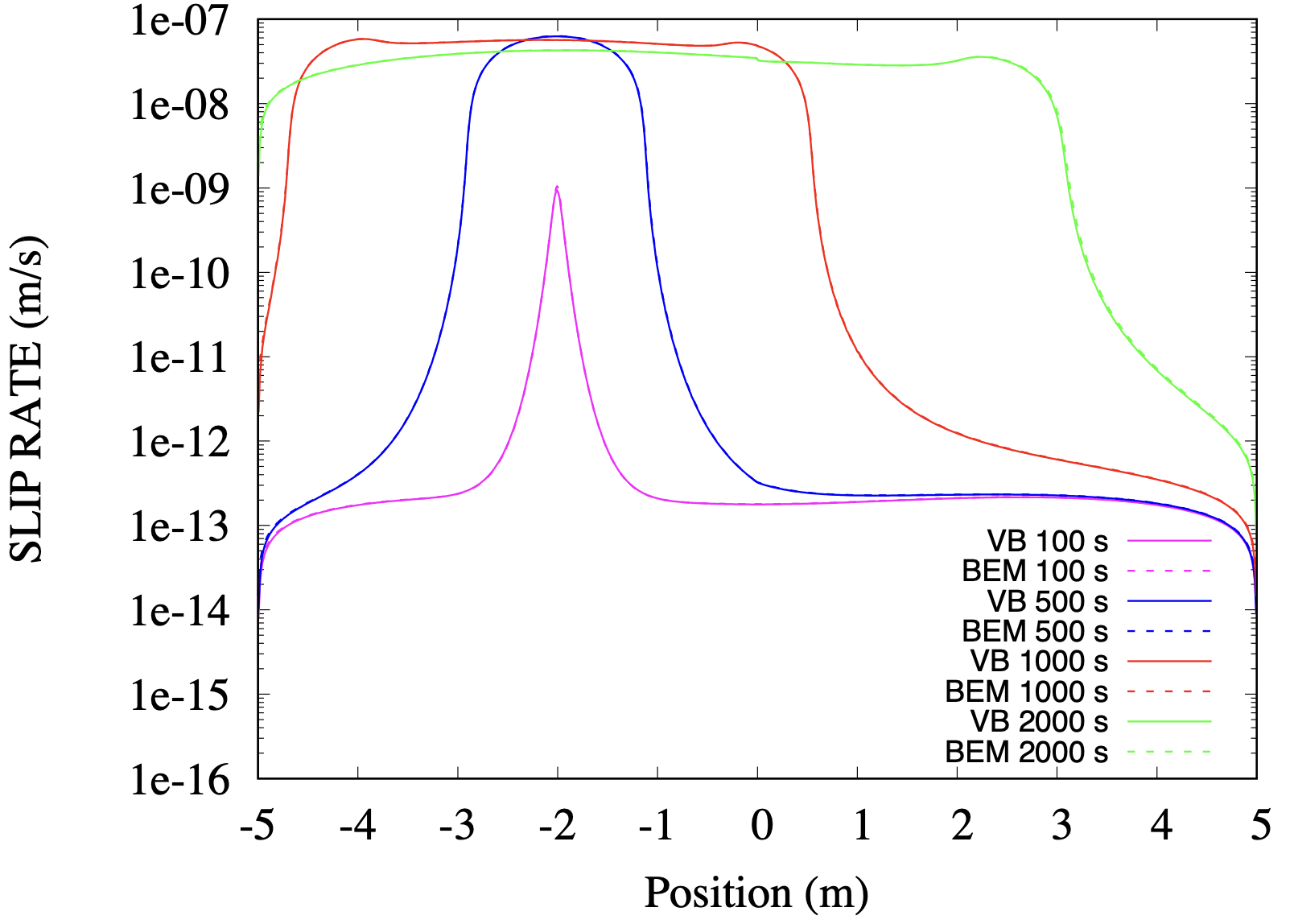}
     \includegraphics[width=0.45\textwidth]{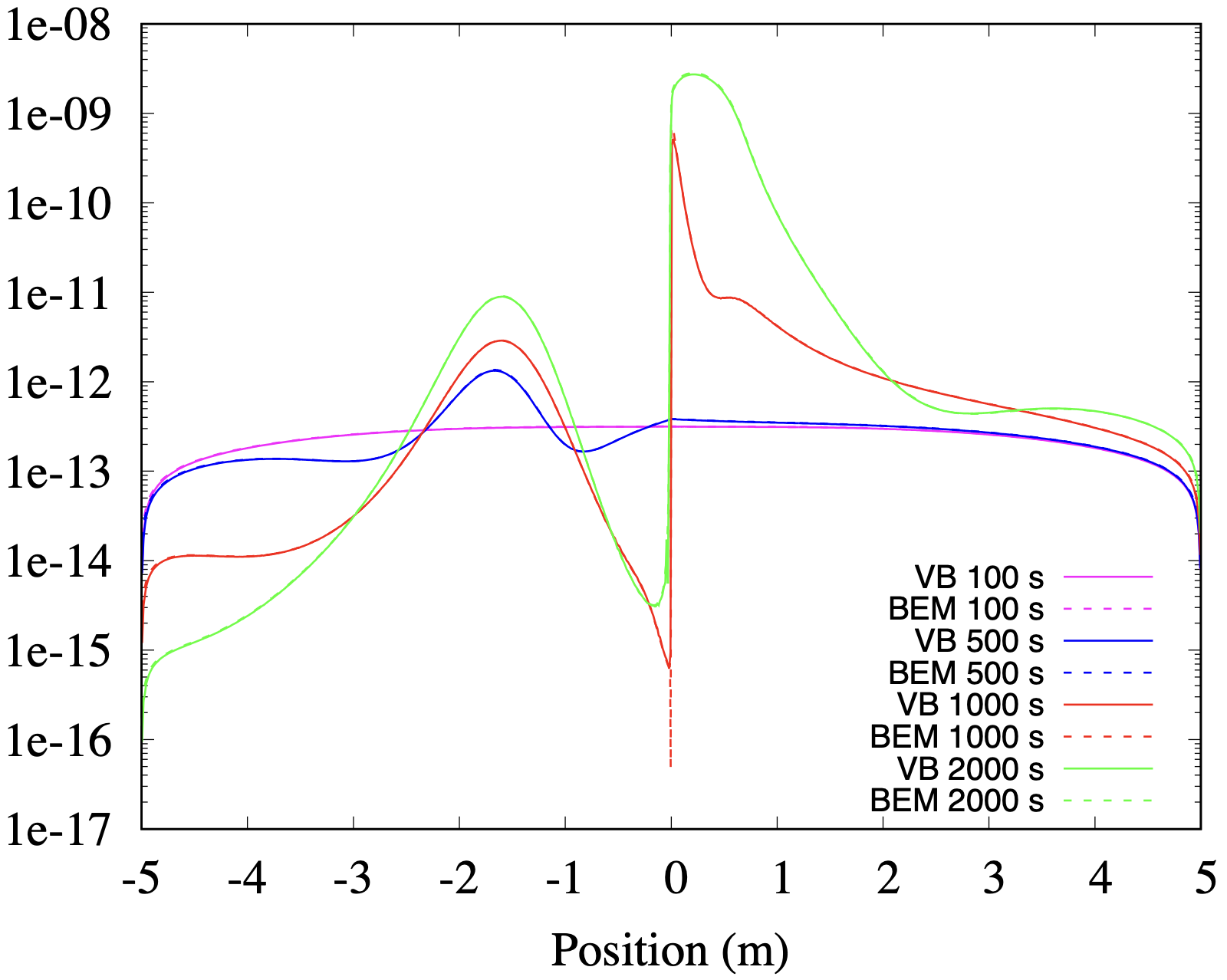}
        \caption{Comparison of the slip rate along fault $1$ (left) and $2$ (right) at times $100$ s, $500$ s, $1000$s and $2000$ s after injection, obtained with the reference BEM and our discretisation (VB).}
        \label{fig-sliprates-2failles}
\end{figure}

The computational domain is chosen as $\Omega =(-100,100)^2$ m$^2$, which is sufficiently large to ensure that the solution is only negligibly affected by the homogeneous Dirichlet boundary conditions imposed on the displacement field. The domain is discretized using a conforming triangular mesh comprising approximately 211000 cells. The mesh conforms to the fault geometry and employs a uniform discretisation along both faults, resulting in a total of $1600$ fault edges. Figure \ref{fig-mesh-2failles} illustrates the mesh and the local refinement in the vicinity of the faults.
Figure \ref{fig-sliprates-2failles} demonstrates the excellent agreement between the slip-rate distributions along both faults obtained with the proposed discretisation and those computed using the reference Boundary Element Method (BEM) code FastCycles \cite{romanet2021}. The comparison is performed at several times using the same fault discretisation in both models, highlighting the accuracy of the proposed approach.

\subsubsection{3D benchmark}
We consider a slightly modified version of the SEAS BP8 Quasi Dynamic 3D benchmark \cite{seasBP8} with a single 2D fault $\Gamma = (-400 \,{\rm m} ,400 \,{\rm m} )^2 \times \{0\}$ embedded into the full space $\R^3$.
The setup corresponds to an early version of the benchmark and differs from the final configuration, most notably in its hydraulic parameters, higher injection rate, and shorter injection shut-off and simulation times. 
The matrix is considered impervious and the fluid is injected during 36000 s in the fault with fixed hydraulic properties and a constant Gaussian flow rate
$$
Q = {q_0 \over 2 \pi D^2} e^{- {r^2\over 2 D^2}} 
$$
with $r$ the distance to the origin, $q_0 = 0.05$ m$^3$.s$^{-1}$, $D = 50$ m.  
The Lam\'e coefficients are set to $\mu = \lambda = 32.04$ GPa and the rate-strengthening parameters are set to $a = 0.016$, $b=0.01$, $d_c = 5~10^{-4}\, {\rm m}$, $v_0 = 10^{-6}\, {\rm m}.{\rm s}^{-1}$, $f_0 = 0.6$, and the radiation damping coefficient is defined by $\eta = {\mu\over 2 c_s}$ with $c_s = \sqrt{ \mu \over \rho_s}$, $\rho_s = 2670$ kg.m$^{-3}$.
The fault is pre-stressed with a compressive normal traction set to $25~10^6$ Pa and a tangential traction set to $14.6 ~10^6$ Pa in the y direction. The initial fault slip and pressure are set to zero and the initial slip rate is fixed to $10^{-12}$ in the direction $y$ of the initial tangential traction. The initial  state $\theta$ is computed in such a way that the RSF friction law is satisfied at initial state. 

The fault conductivity is set to $C_f = d_f K_f$ with fault permeability $K_f = 10^{-14}$ m$^2$ and hydraulic aperture fixed to $\d_f = 5 \,{\rm m}$ with no penetration ($K_\n = \infty$) and no shear dilation ($\psi = 0$). The fluid compressibility is $c_f = 10^{-8}$ Pa$^{-1}$ and its dynamic viscosity is set to $\mu_f = 0.89 ~10^{-3}$ Pa.s. 

\begin{figure}[H]
    \centering
        \includegraphics[width=0.25\textwidth]{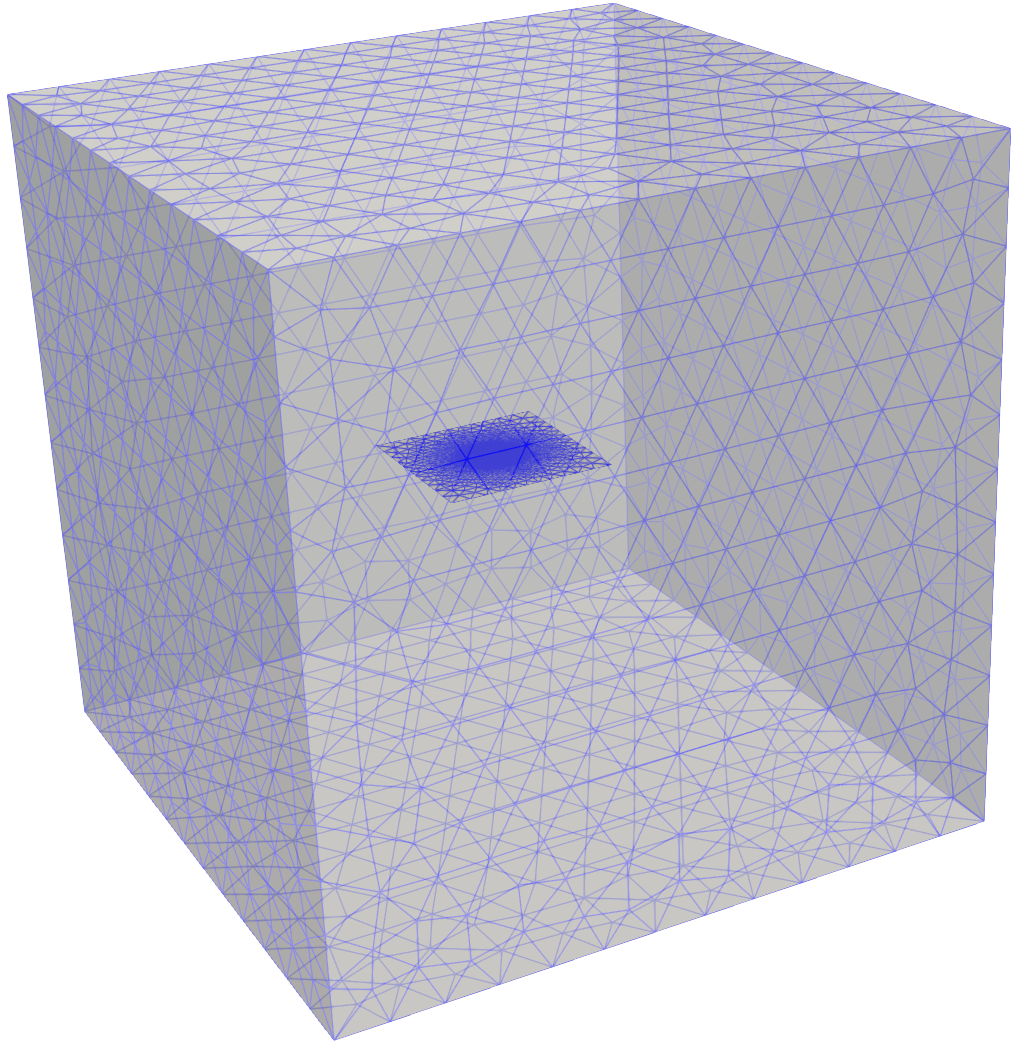}
        \caption{Tetrahedral mesh  with 6000k cells, 622k nodes, 31k fault faces, locally refined with a size down to 2.5 m at the fault in a ball of radius 150 m centered at the origin. }
        \label{fig-mesh-3Dbench}
\end{figure}

\begin{figure}[H]
    \centering
    \includegraphics[width=0.45\textwidth]{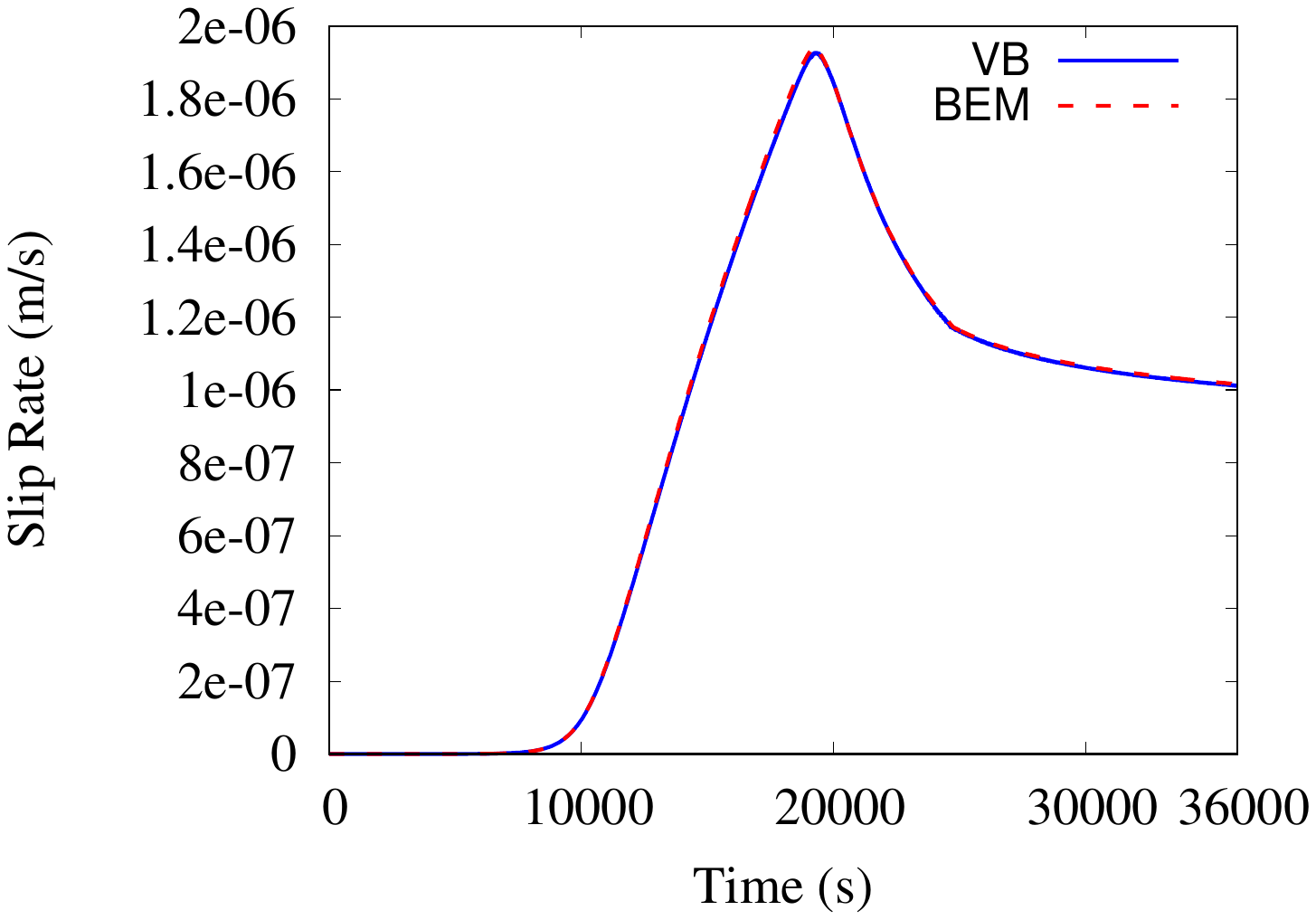}
    \includegraphics[width=0.45\textwidth]{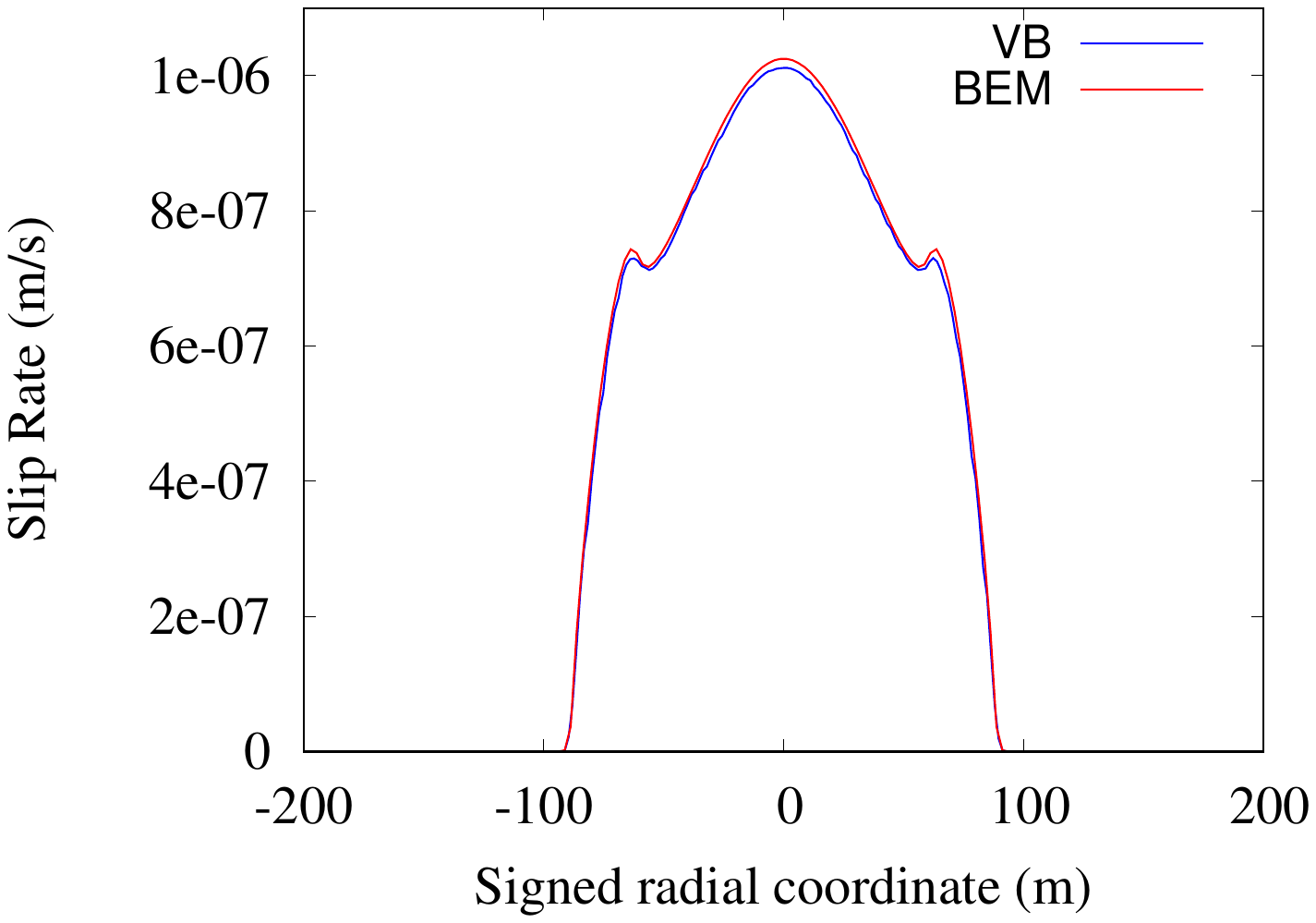}    
        \caption{(Left): comparison of the maximum slip rate as a function of time obtained with the reference BEM and our discretisation (VB). (Right): comparison of the diagonal cut of the $x$ component of the slip rate $\vslip$ at final time obtained with the reference BEM and our discretisation (VB).}
        \label{fig-3Dbench-VB-BEM}
\end{figure}

Our solution computed on the bounded domain $\Omega = (-1600 \,{\rm m} ,1600 \,{\rm m} )^3$ using the gmesh tetrahedral mesh  exhibited in Figure \ref{fig-mesh-3Dbench} is compared with the solution obtained using FastCycles3D, the 3D extension of the BEM code FastCycles \cite{romanet2021}, on a uniform cartesian mesh of size 2.5 m. Figure  \ref{fig-3Dbench-VB-BEM} exhibits a very good match of the maximum slip rate as a function of time and of the slip rate radial cut at final time. 


\subsection{Fluid injection 2D test case}\label{sec:testcoupled2D}

\begin{figure}[H]
    \centering
    \begin{subfigure}{0.38\textwidth}
        \centering
        \includegraphics[scale=0.15]{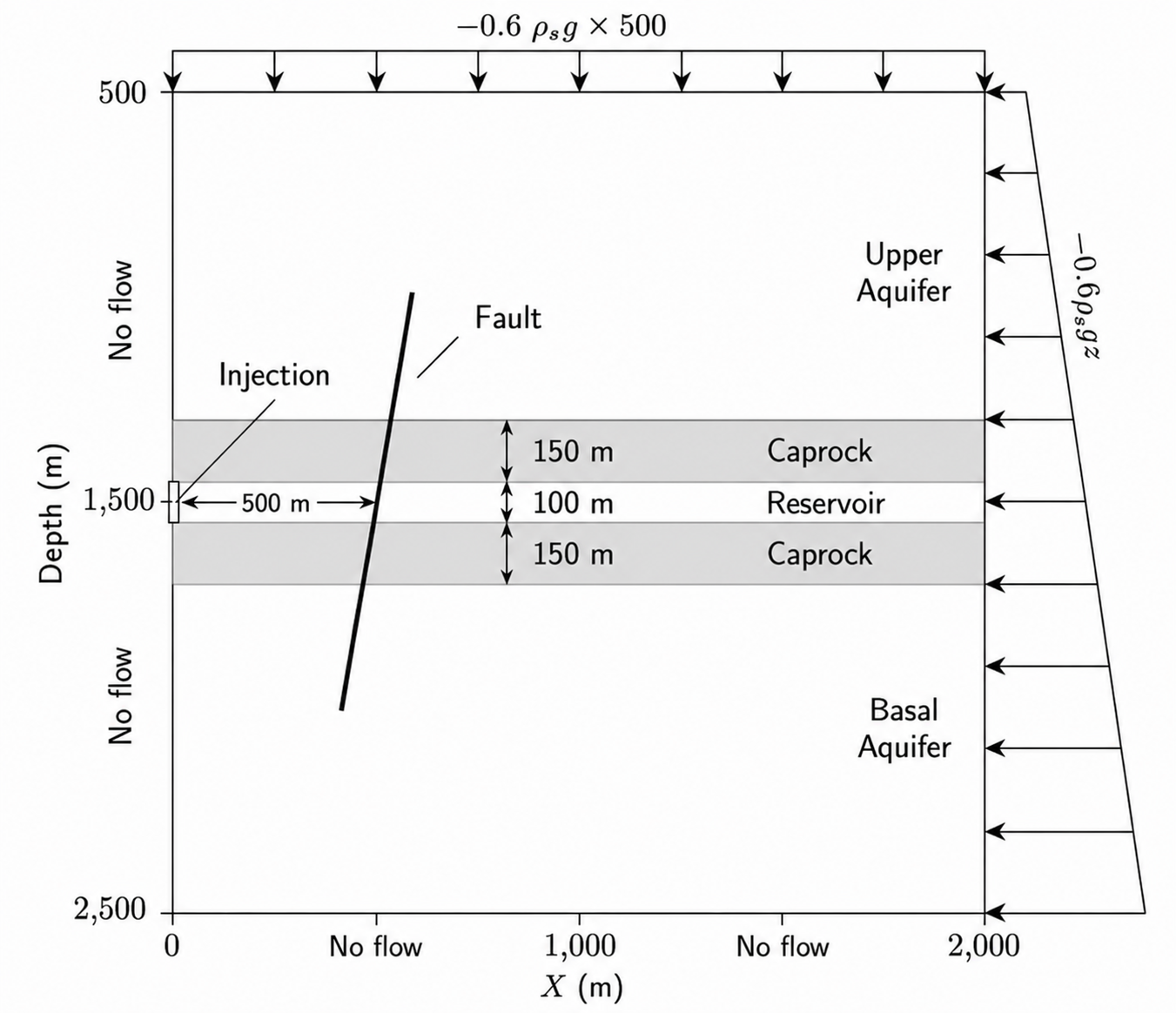}
        \caption{}
    \end{subfigure}
    \begin{subfigure}{0.38\textwidth}
        \centering
        \includegraphics[scale=0.2]{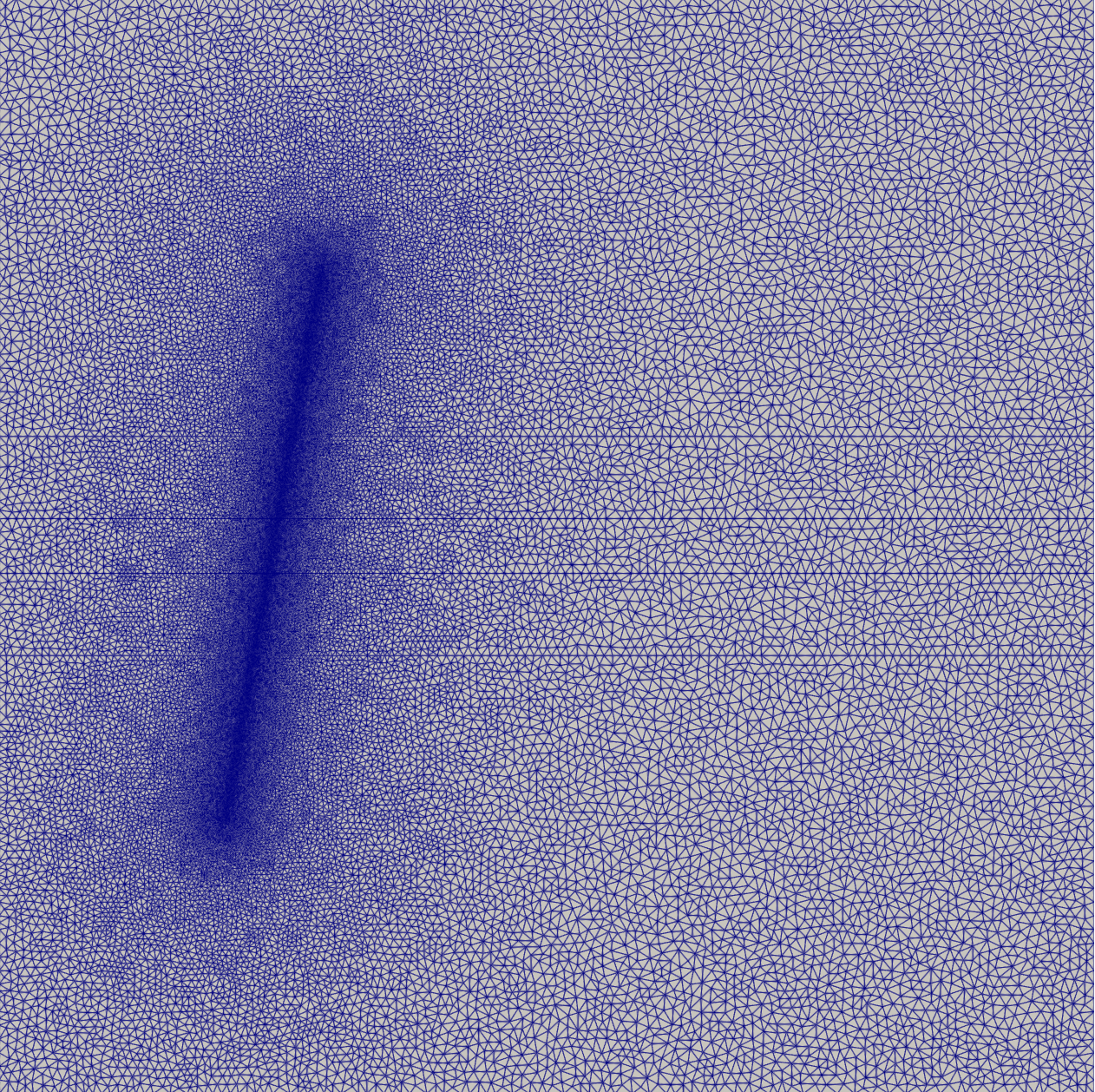}
        \caption{}
    \end{subfigure}
    \caption{(a) Model set up for the fluid injection test case from \cite{cappa2011}. (b) coarsest triangular 2D mesh refined along the fault, featuring 677 fault faces and 77k cells. }
       \label{fig:setup_mesh}
\end{figure}

We present a fluid-induced 2D fault reactivation model that simplifies the \ce{CO2} injection benchmark from \cite{cappa2011} by assuming single-phase isothermal flow for a slightly compressible liquid and taking into account the gravity term in the Darcy law. The geometry (see Figure~\ref{fig:setup_mesh}(a))
comprises three distinct layers whose properties are adopted from Table~1 of \cite{cappa2011}, featuring an aquifer with fixed isotropic permeability $k_1 = \SI{1e-15}{m^2}$ and initial porosity $\phi^0_m = 0.1$, a reservoir with fixed isotropic permeability $k_2 = \SI{1e-13}{m^2}$ and initial porosity $\phi^0_m = 0.1$, and a caprock with fixed isotropic permeability $k_3 = \SI{1e-17}{m^2}$ and initial porosity $\phi^0_m = 0.01$. Note that these matrix permeabilities are considered fixed for simplicity in this test case. 
The  shear dilation angle will be set to either $\psi = \ang{5}$ or $\psi = \ang{0}$ and the initial hydraulic aperture to $\d_f^0 = \SI{1e-5}{m}$. 
The model configuration incorporates key hydromechanical parameters, including the Biot coefficient $\alpha=0.8$, the Biot modulus $M = 10$ GPa, the Lam\'e coefficients $\mu = \lambda = 4$ GPa and the normal stiffness $K_\n = 60$ GPa.m$^{-1}$.
The fluid compressibility is $c_f = 10^{-10}$ Pa$^{-1}$, the density in the gravity term is fixed to $\rho_{ref} = 10^3$ kg.m$^{-3}$; and the fluid dynamic viscosity is set to $\mu_f = 10^{-3}$ Pa.s.
For the flow, a fixed flow rate $q_{\mathrm{inj}}=$\SI{2e-3}{kg.m^{-2}.s^{-1}} is injected at the left boundary of the reservoir layer; the right boundary is open to the flow with a fixed hydrostatic pressure, and the remaining boundaries are considered impervious. The injection is stopped at time $t_{\mathrm{inj}} =10$ days and the simulation time is fixed to $T=15$ days. 

For the mechanics, a zero normal displacement together with a zero tangential component of the traction vector are imposed on the left and bottom boundaries. At the top boundary, the traction vector is given by the gravity load with an average rock and fluid density fixed to $2500$ kg.m$^{-3}$ and a depth of $500$ m. At the right boundary, the traction vector is imposed assuming an extensional stress regime with $\sigma_h = 0.6\sigma_v$. The initial displacement and tractions are computed at given initial hydrostatic pressures $p_m^0$, $p_f^0$ from the mechanical equilibrium with contact conditions taking into account the normal stiffness term and a zero slip (stick mode). The initial slip rate is fixed to $10^{-12}$ in the direction of the initial tangential traction, and the initial  state $\theta$ is computed in such a way that the RSF friction law is satisfied at initial state. The RSF law is defined by the parameters $d_c = 5~10^{-4}\, {\rm m}$, $v_0 = 10^{-6}\, {\rm m}.{\rm s}^{-1}$, $f_0 = \tan(\ang{25})$, and the radiation damping coefficient is defined by $\eta = {\mu\over 2 c_s}$ with $c_s = 3464$ m.s$^{-1}$.

In the following test cases, we will consider two configurations with either the rate-strengthening $a = 0.016$, $b=0.01$ RSF parameters combined with shear dilation angle $\psi = \ang{5}$ or the rate-weakening $a = 0.005$, $b=0.02$ RSF parameters combined with no shear dilation. The fault conductivity is defined by Poiseuille's law with $C_f = {(d_f)^3   \over 12}$ and the fault normal transmissivity $\Lambda_f$ is defined by the matrix permeability at the interface divided by twice the initial fault aperture which roughly amounts to a pressure continuity assumption at the matrix fault interface.  Table \ref{tab:parameters_2D} in the Appendix provides the detailed set of parameters and data used for the 2D fluid injection test cases. 

In order to investigate the robustness of the nonlinear and linear solvers with respect to the mesh, the simulations are performed on a family of three 2D triangular meshes with 677, 1263, 2432 fault edges and respectively 77k, 304k and 1214k cells (see Figure~\ref{fig:setup_mesh} for the coarsest mesh). Note that all the simulations are run on a 3D code using an extruded mesh with one layer of cell of width 1 m. 

\begin{figure}[H]
  \centering
  \includegraphics[scale=0.23]{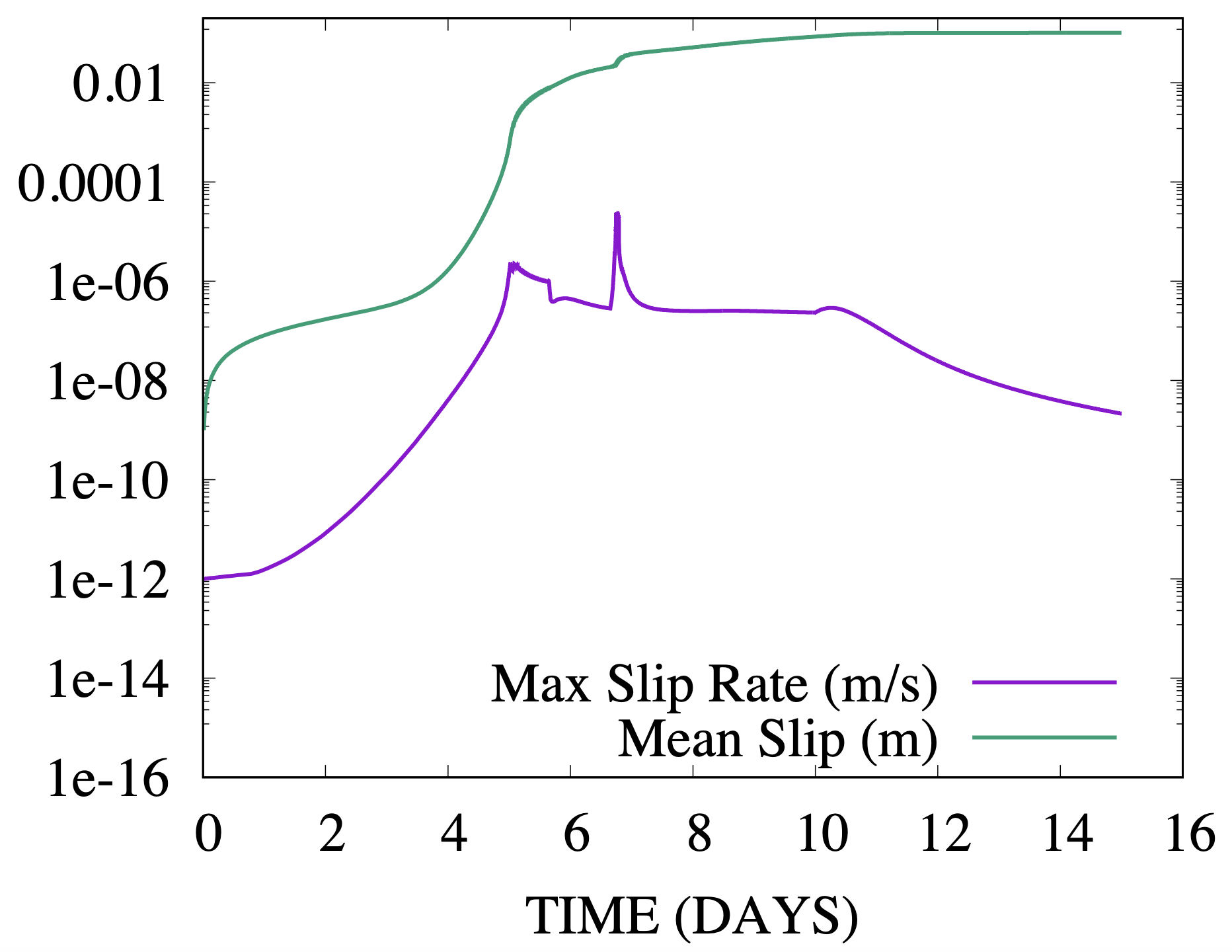}
  \includegraphics[scale=0.25]{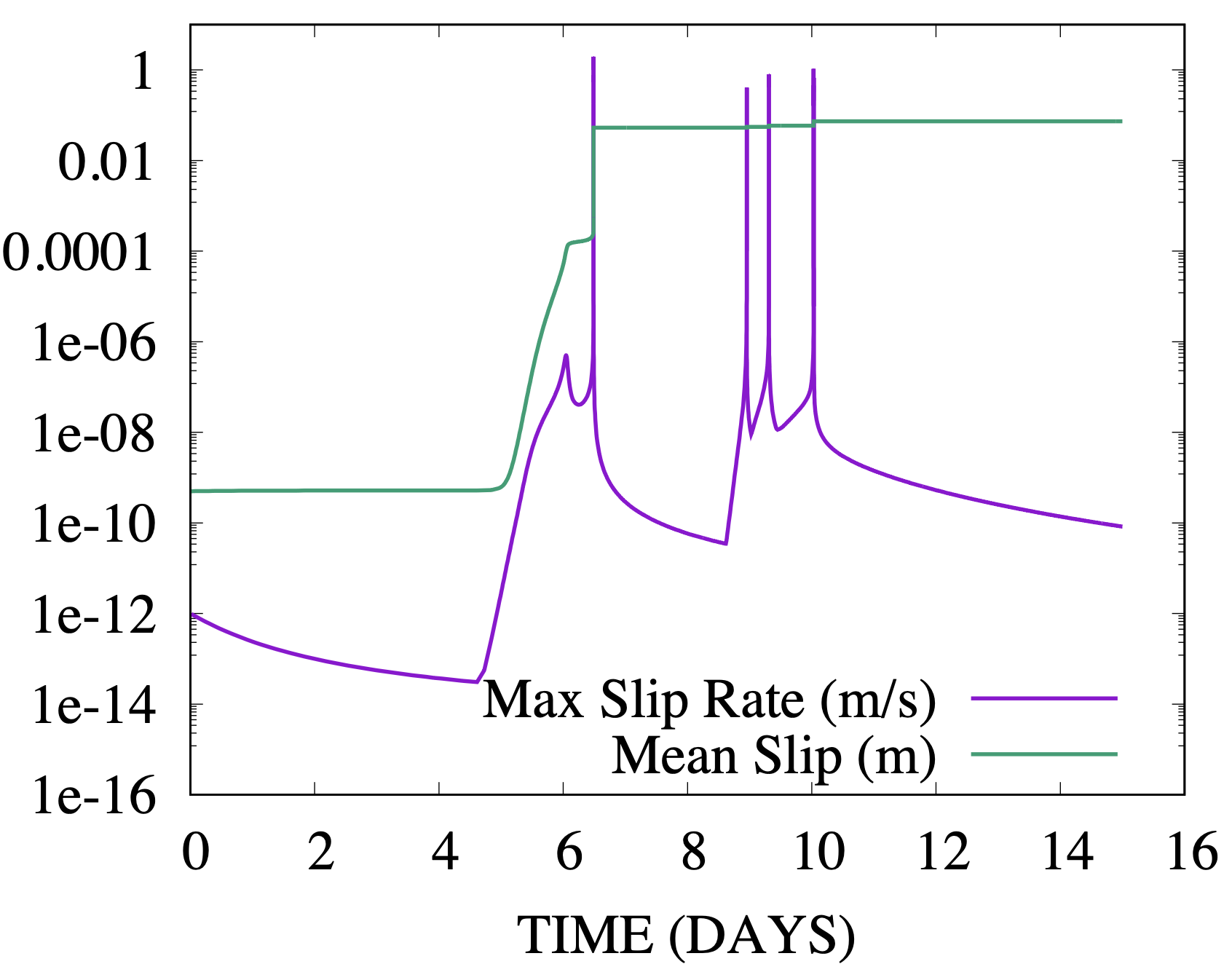}
\caption{Maximum slip rate and mean slip as a function of time for the rate-strengthening with shear dilation (left) and the rate-weakening without shear dilation (right) configurations of the 2D fluid injection test case. The simulations are performed on the finest mesh with $2432$ fault faces and with the time stepping parameter $M_{\Delta t} = 4$.}
       \label{fig:maxsliprate-meanslip}
\end{figure}

\begin{figure}[H]
  \centering
  \includegraphics[scale=0.2]{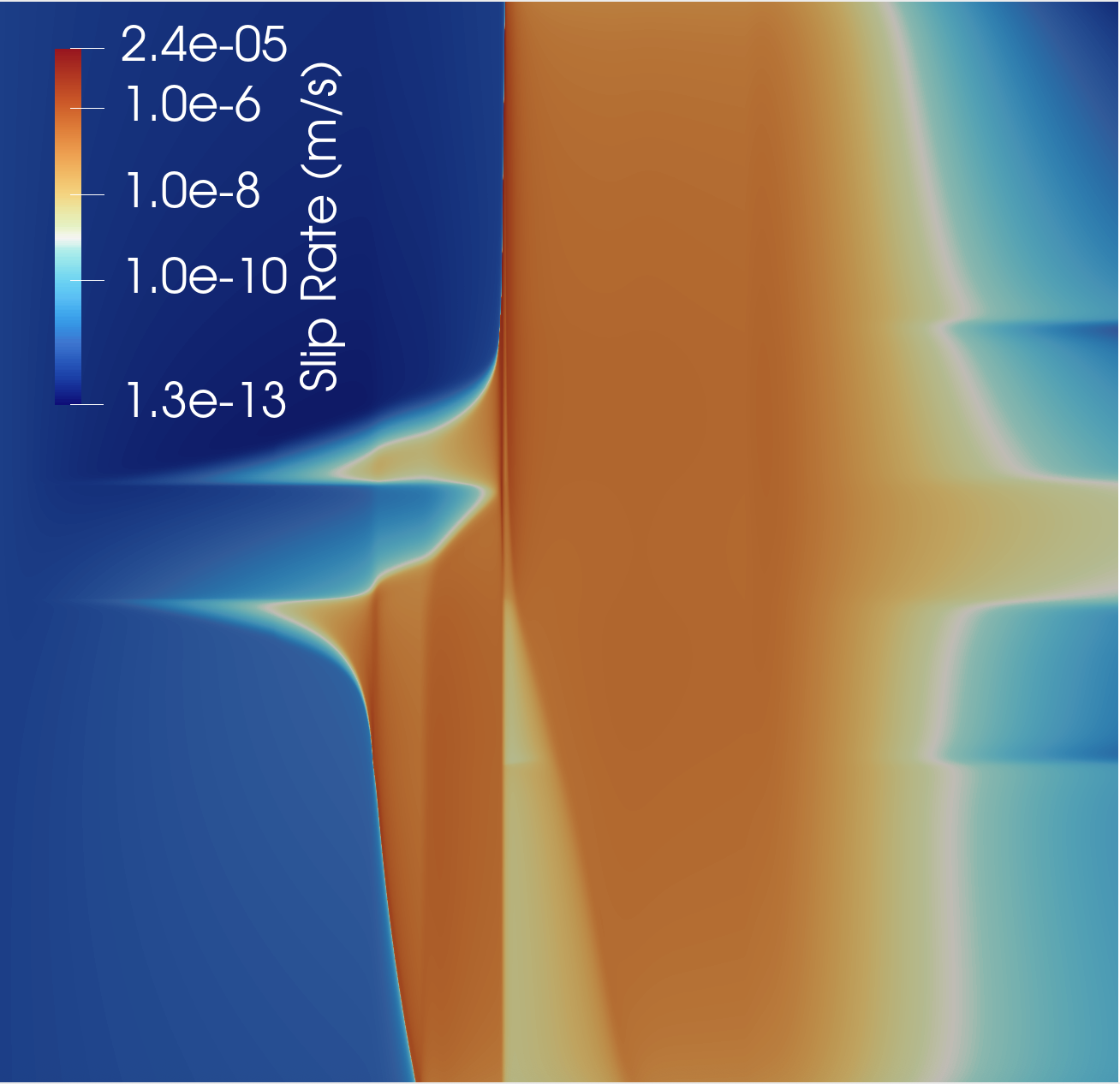}
  \includegraphics[scale=0.2]{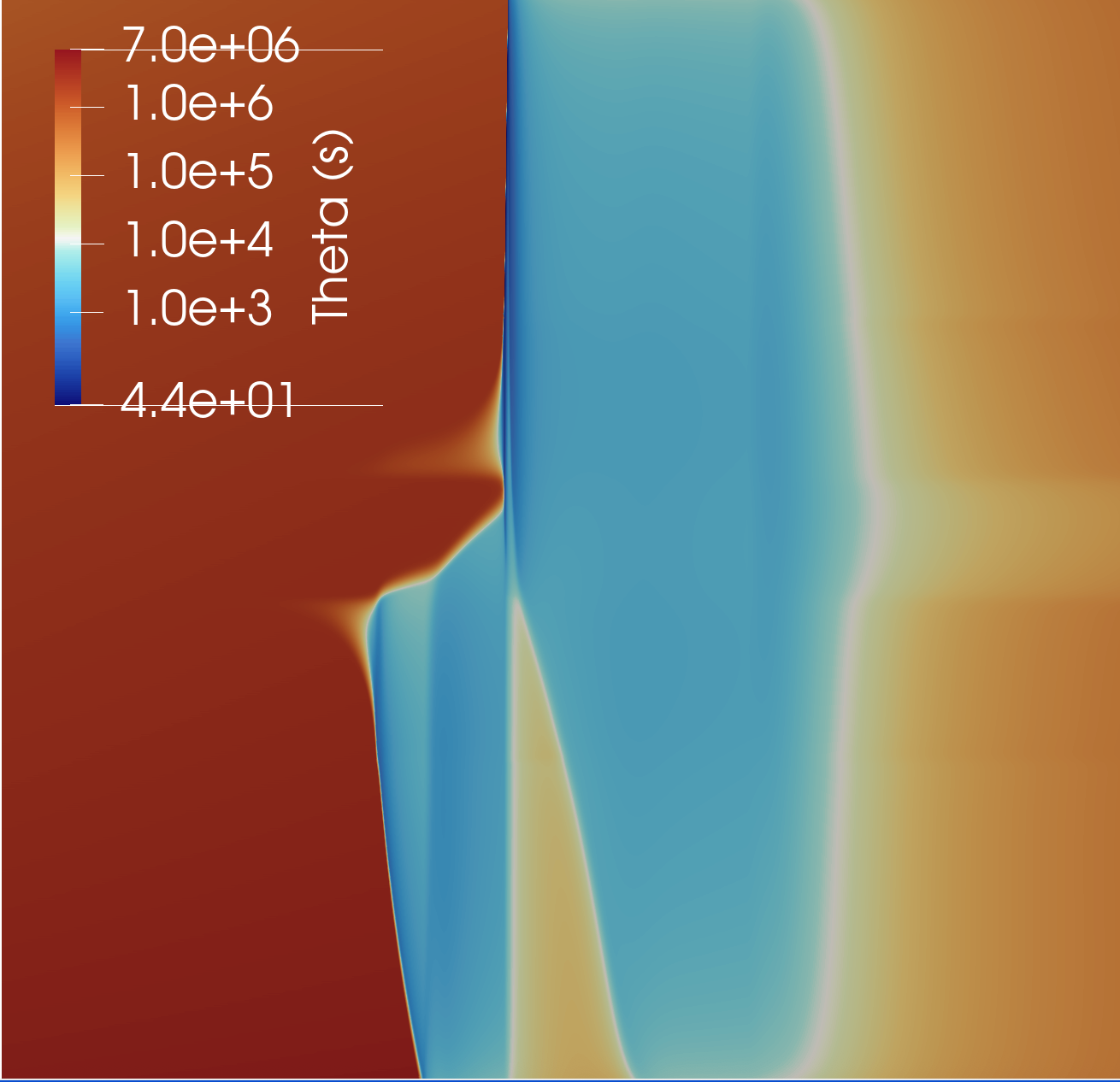}
  \includegraphics[scale=0.2]{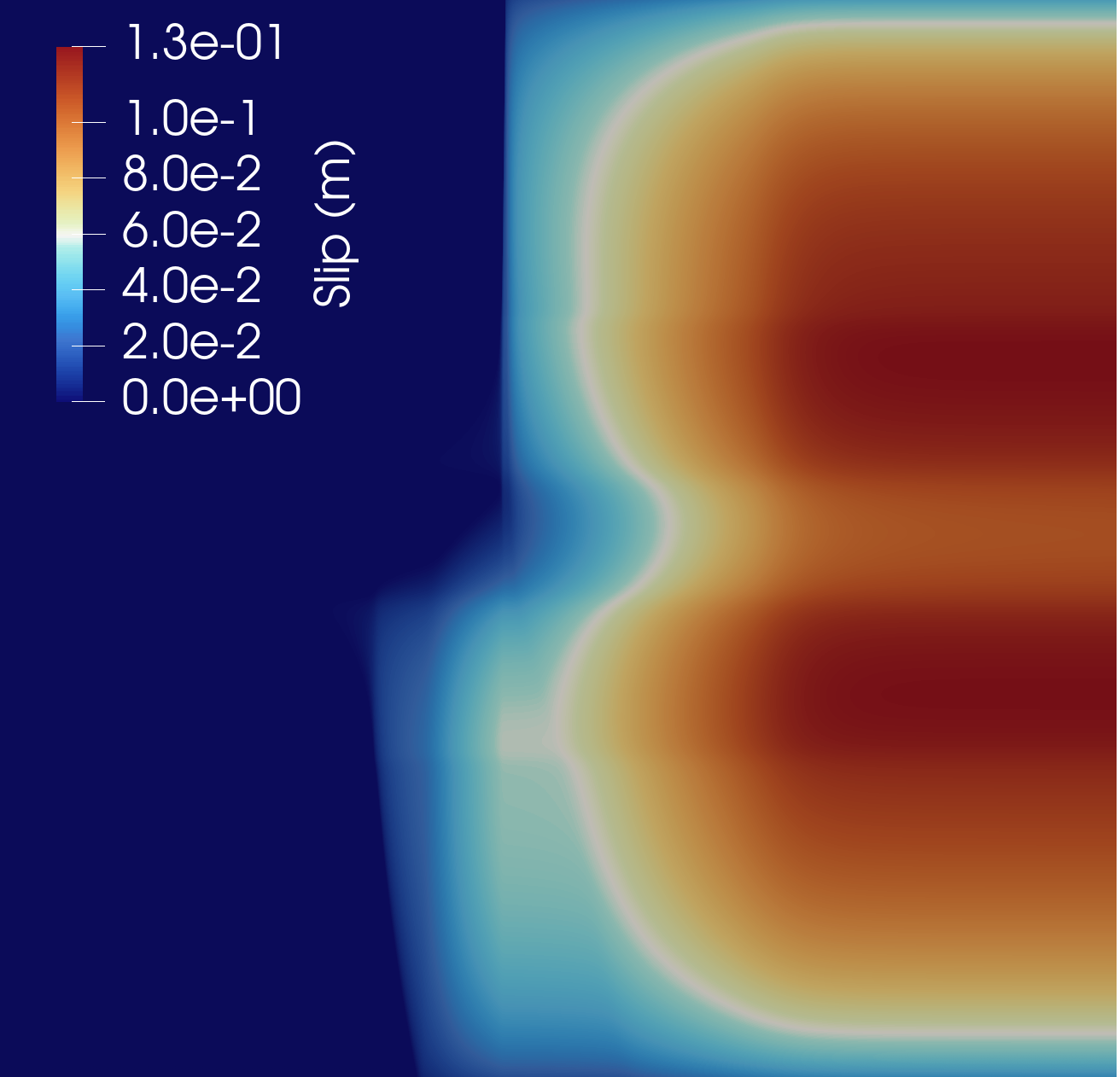}\\
  \includegraphics[scale=0.2]{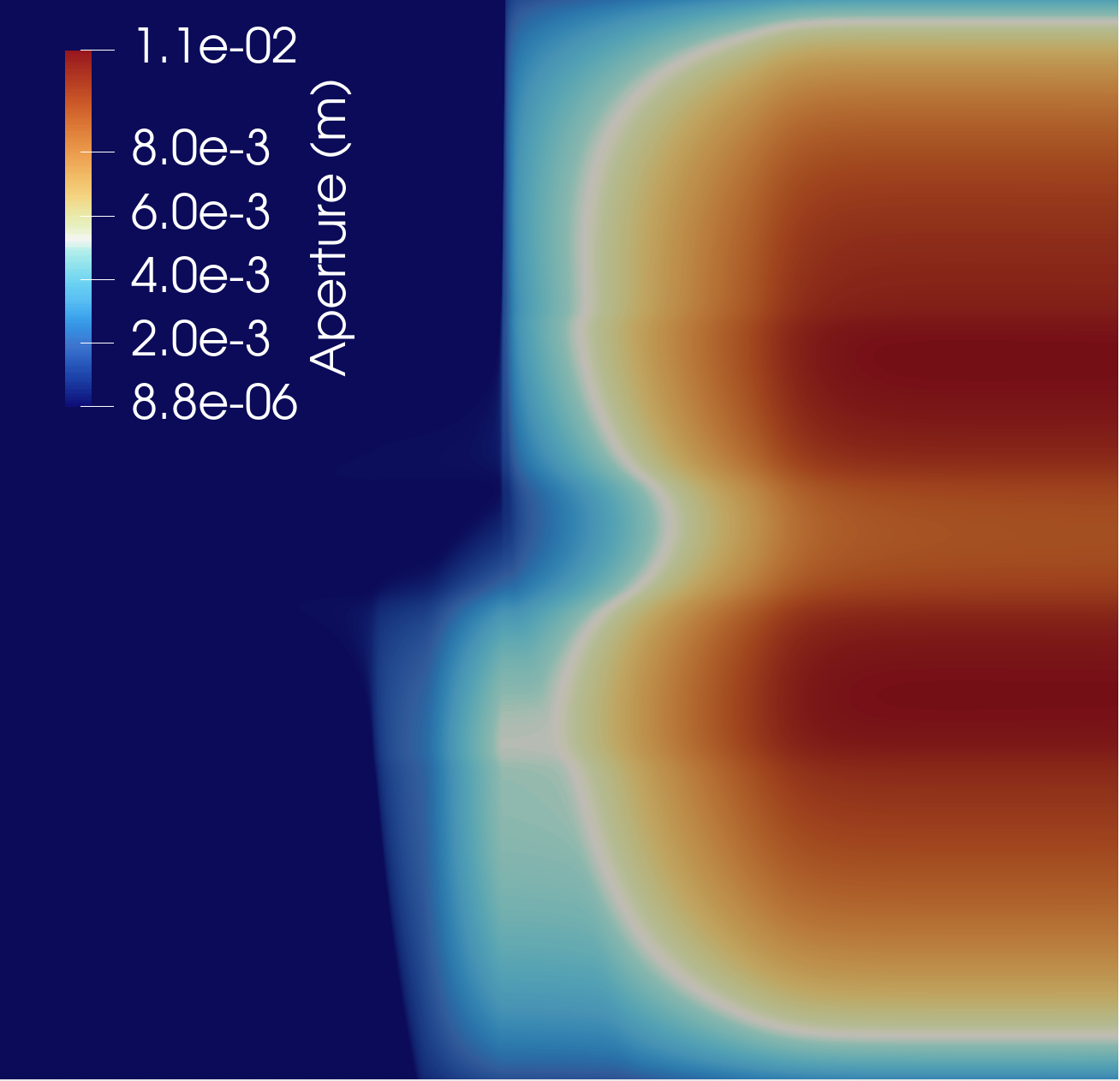}
  \includegraphics[scale=0.2]{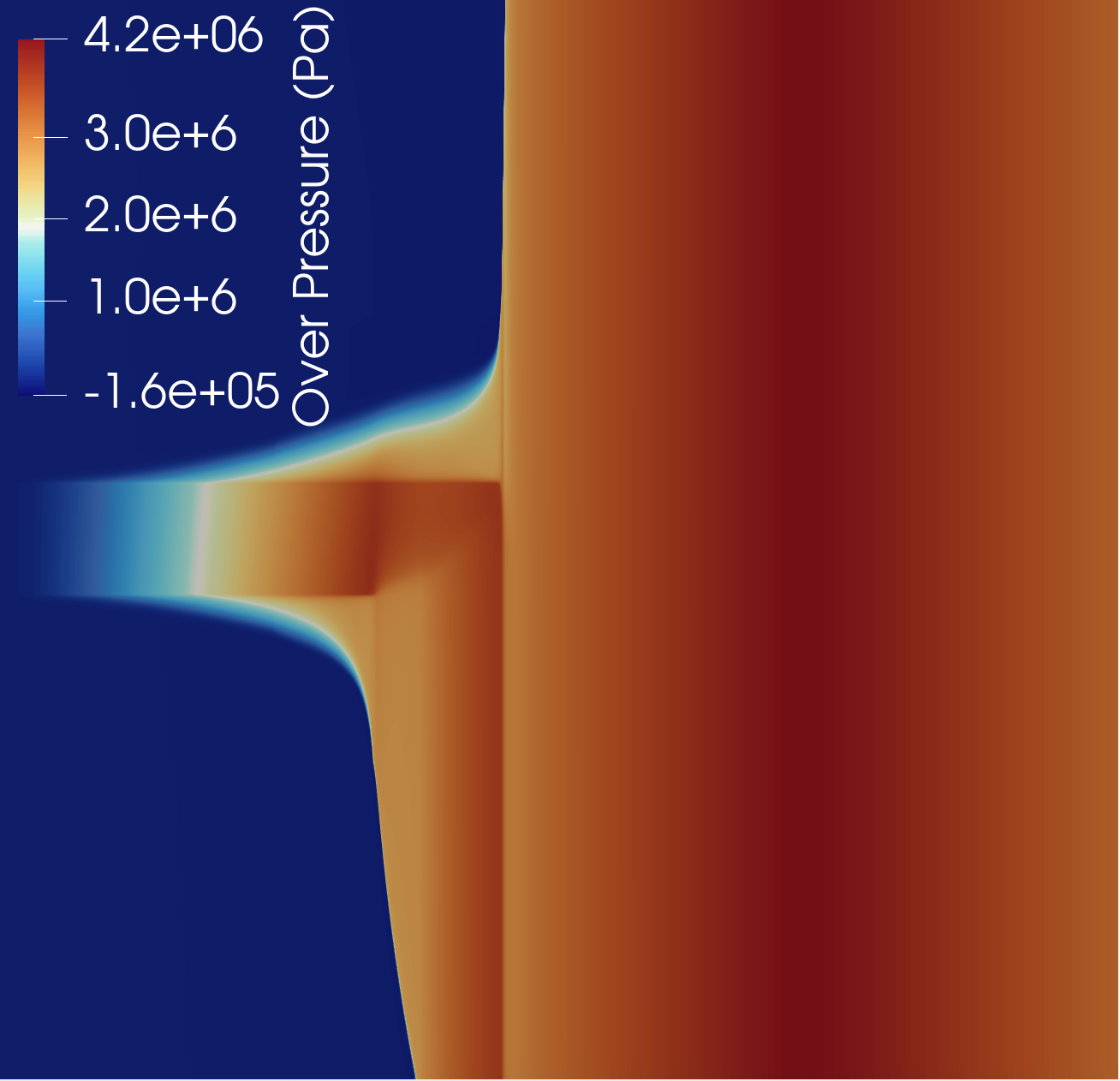}
   \includegraphics[scale=0.2]{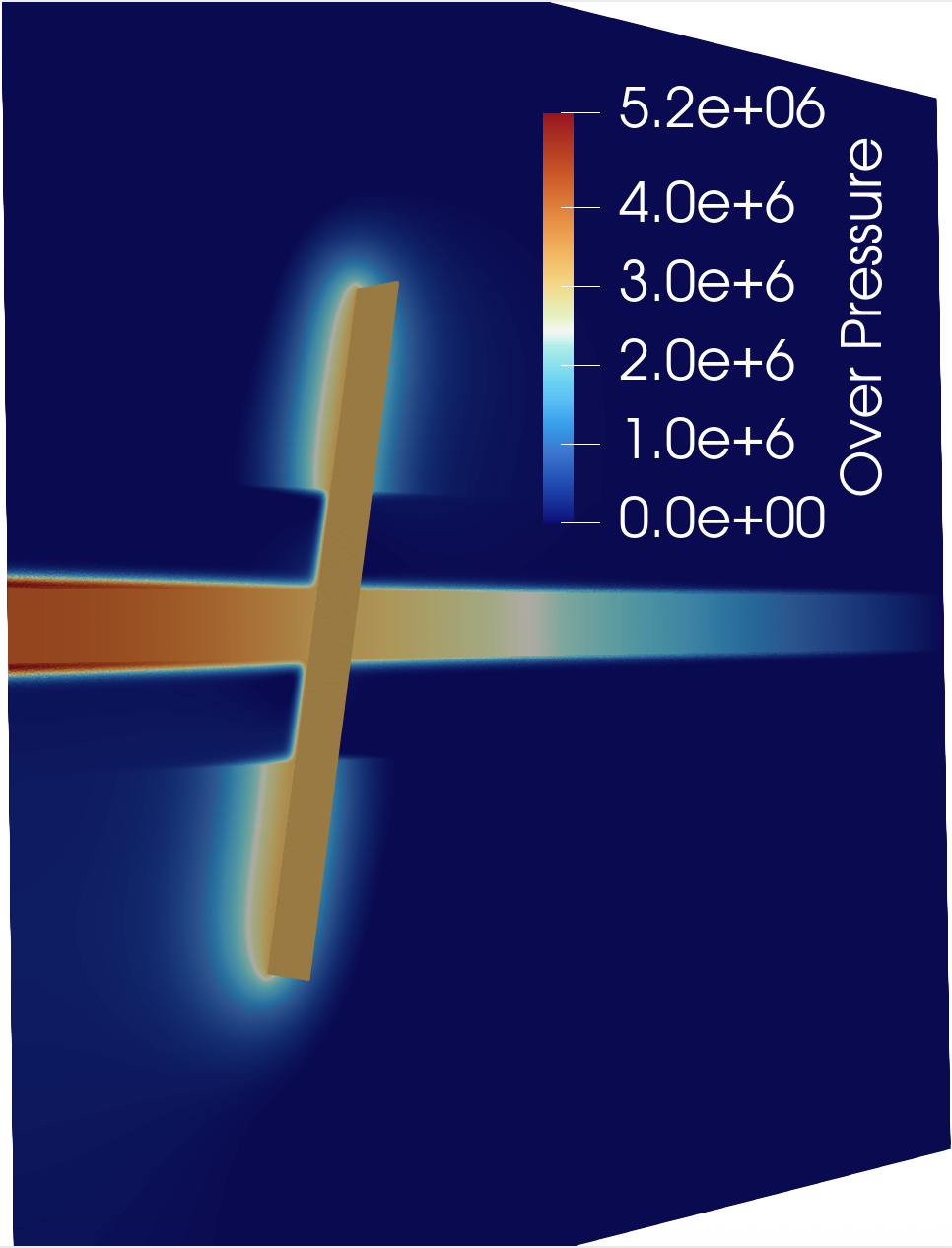}
\caption{2D fluid injection test case with rate-strengthening and shear dilation configuration simulated on the mesh with $2432$ fault faces and with adaptive time stepping defined by $M_{\Delta t}=4$, $\Delta t_{max} = \Delta t^0=10^3$ s. The five first figures from left to right and top to bottom plot in this order $\vslip$, $\theta$, $\slip$, $d_f$ and the over pressure as a function of time in the x direction and of space along the fault in the z direction. The last bottom right figure plots the over pressure at final time both in the matrix and along the fault. }
       \label{fig:strengh-sd}
\end{figure}

\begin{figure}[H]
  \centering
  \includegraphics[scale=0.2]{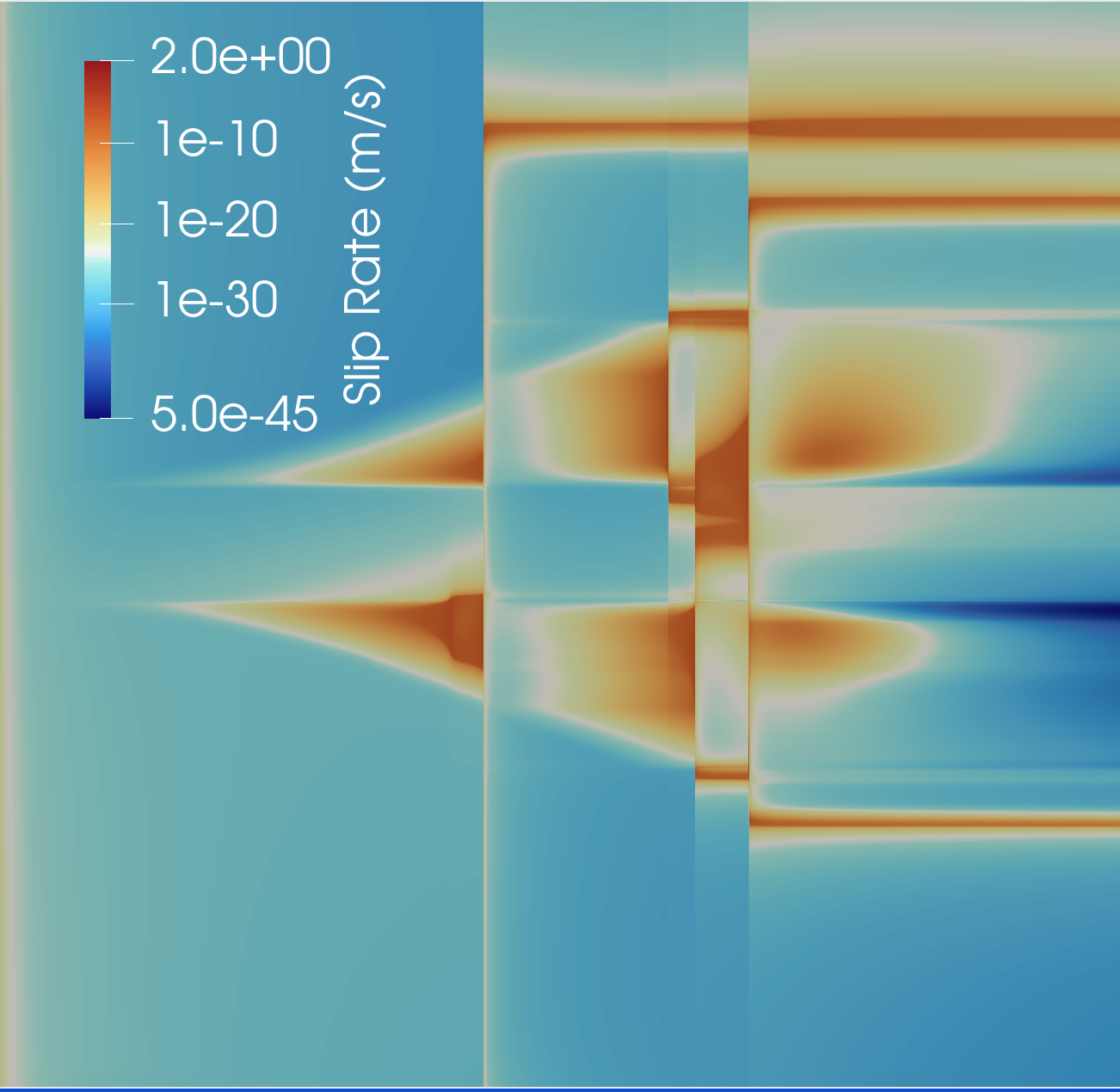}
  \includegraphics[scale=0.2]{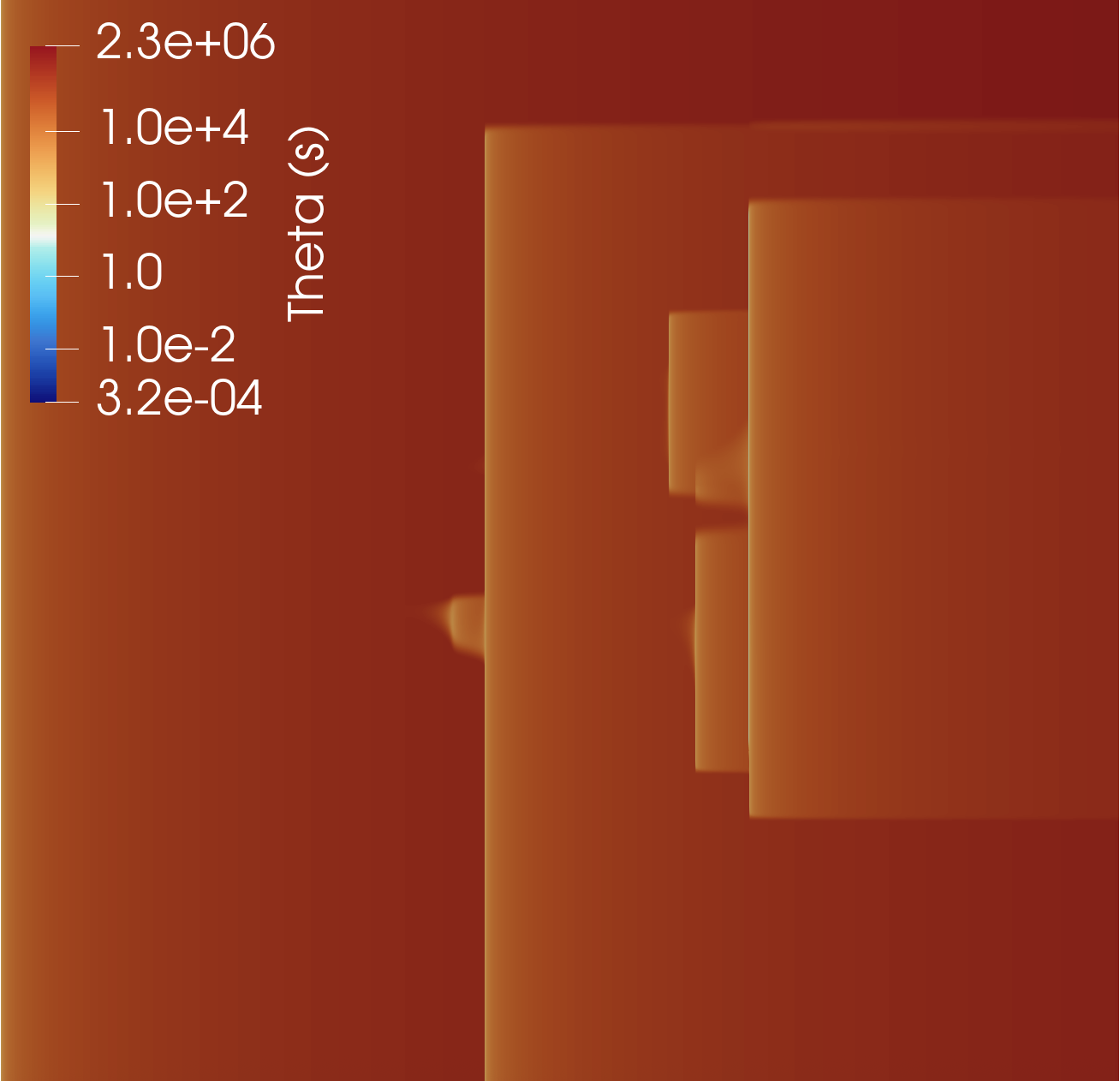}
  \includegraphics[scale=0.2]{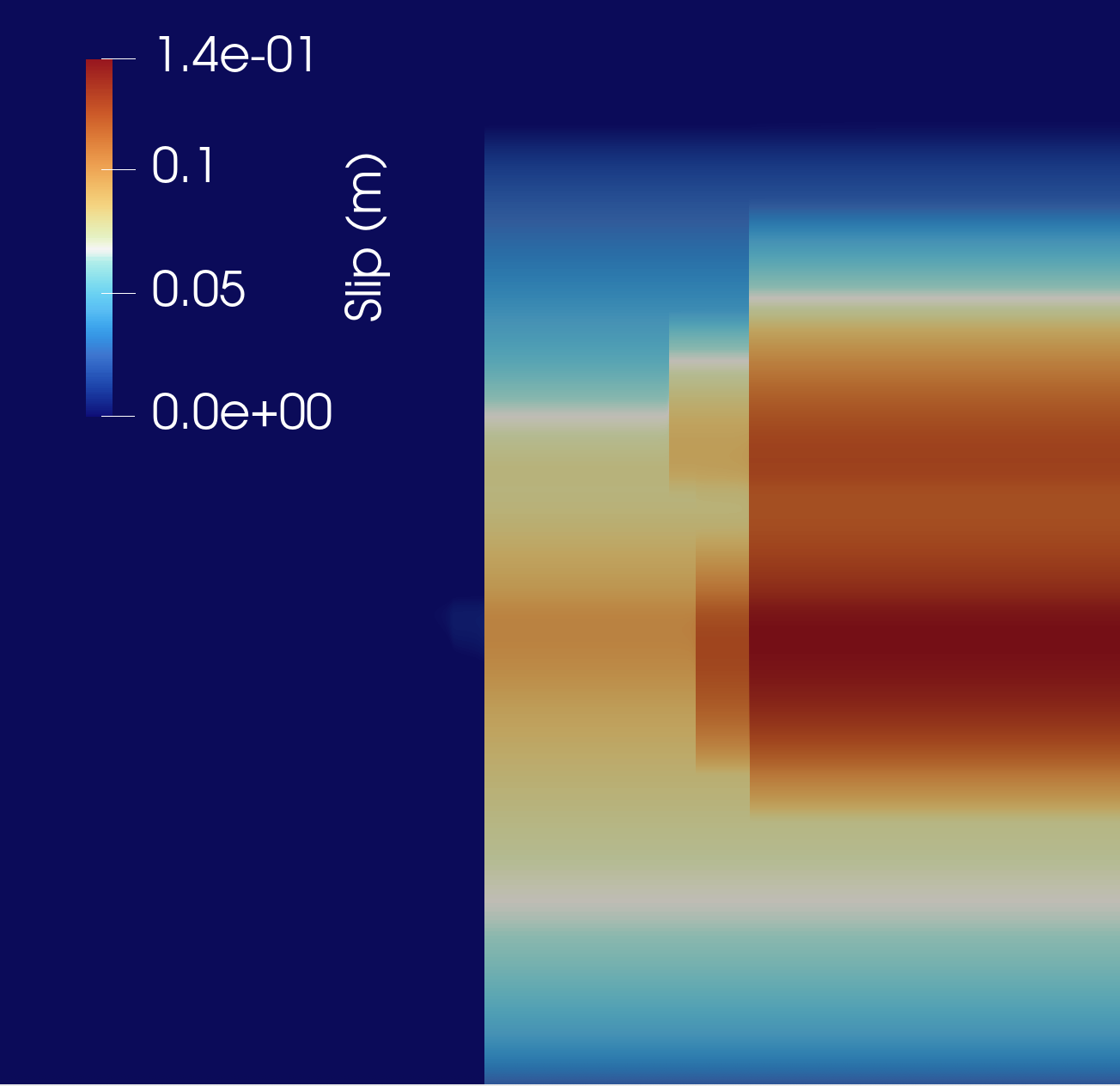}\\
  \includegraphics[scale=0.2]{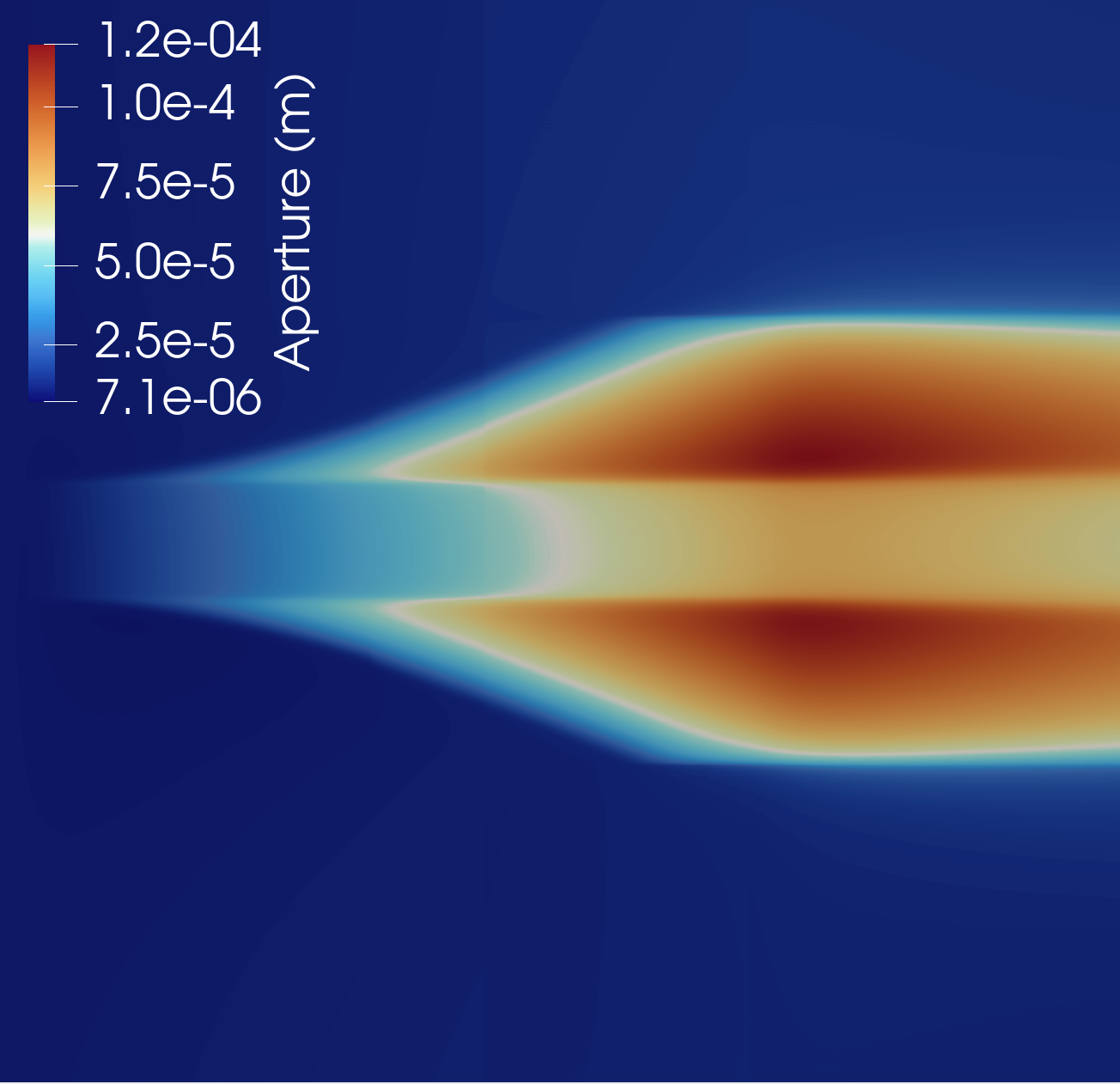}
  \includegraphics[scale=0.2]{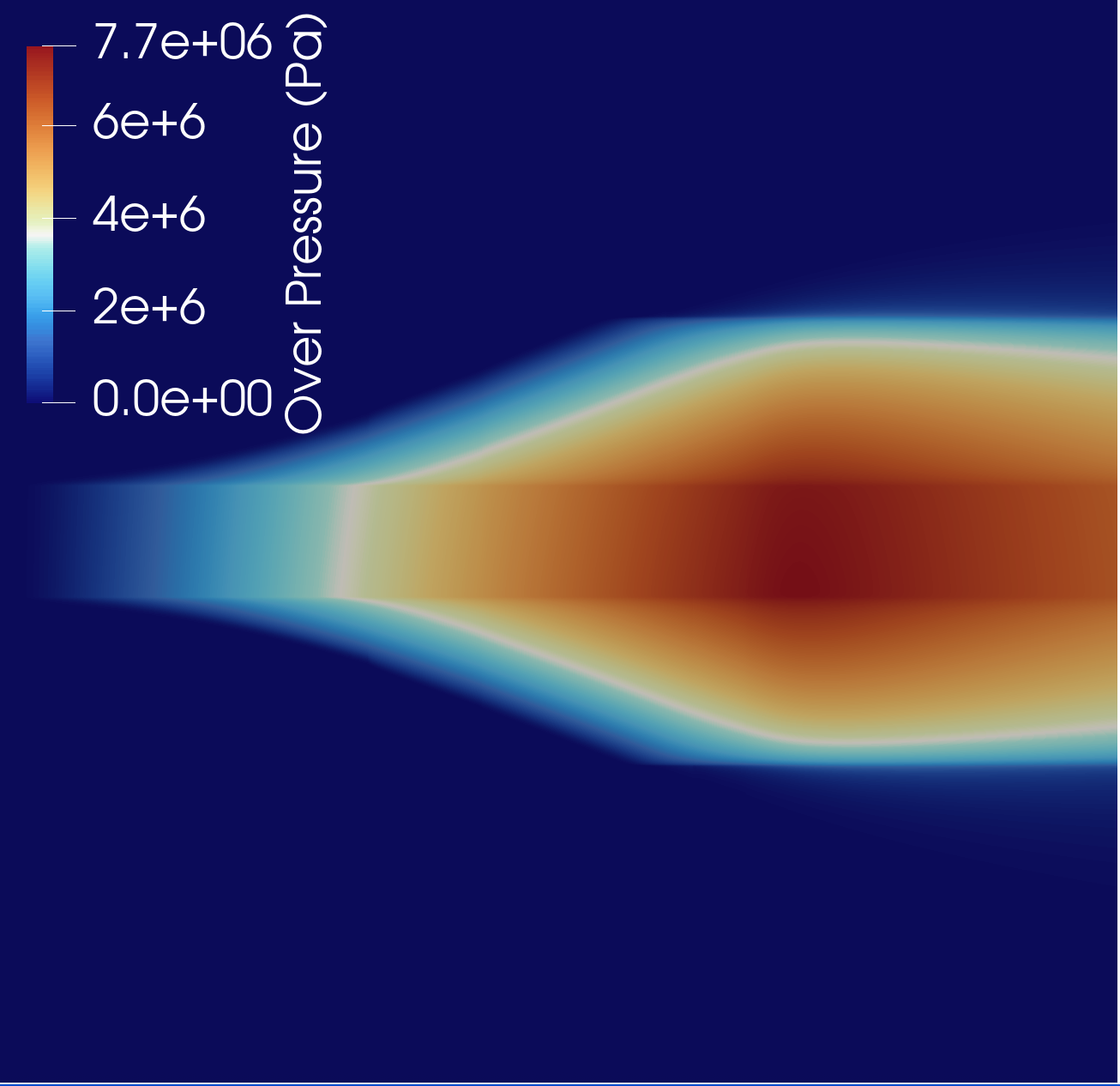}
   \includegraphics[scale=0.2]{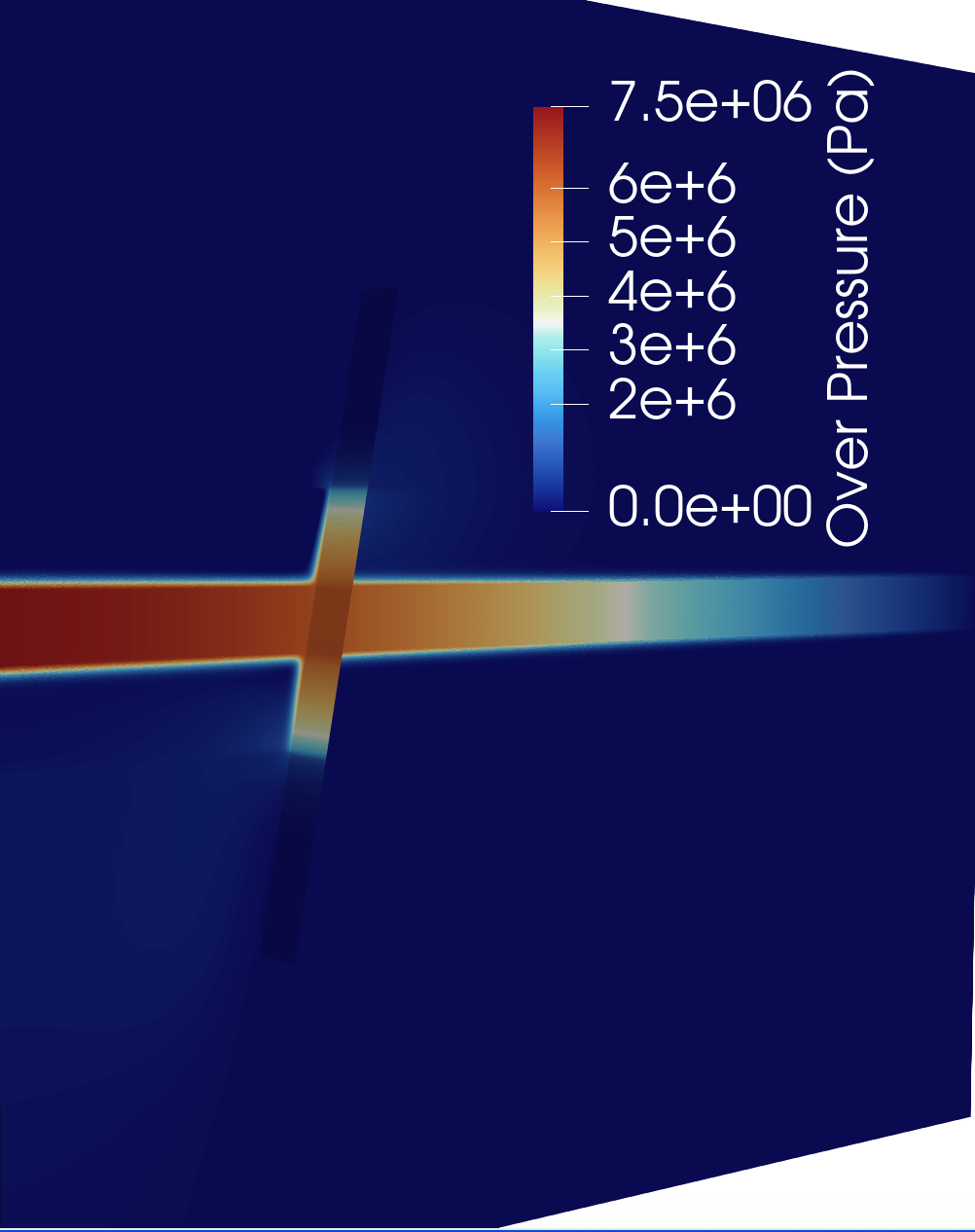}
\caption{2D fluid injection test case with rate-weakening and no shear dilation configuration simulated on the mesh with $2432$ fault faces and with adaptive time stepping  defined by $M_{\Delta t}=4$, $\Delta t_{max} = \Delta t^0=10^3$ s. The five first figures from left to right and top to bottom plot in this order $\vslip$, $\theta$, $\slip$, $d_f$ and the over pressure as a function of time in the x direction and of space along the fault in the z direction. The last bottom right figure plots the over pressure at final time both in the matrix and along the fault. }
       \label{fig:weak-nosd}
\end{figure}

Figure  \ref{fig:maxsliprate-meanslip} shows the maximum slip rate and mean slip as functions of time, while Figures \ref{fig:strengh-sd} and  \ref{fig:weak-nosd} present the solutions along the fault in the time-space domain $(t,z) \in (0,15,\rm{days}) \times (0,1000,\rm{m})$.

In both configurations, slip nucleates on the fault at the lower reservoir-barrier interface and initially propagates backward along the fault-barrier interface. This behavior results from the combination of the high fault pressure within the reservoir and the lower matrix pressure within the barrier.

The comparison between the two configurations highlights the key role of shear dilation in the coupled hydro-mechanical processes. When shear dilation is included, as shown in Figure \ref{fig:strengh-sd} and in the left panel of Figure \ref{fig:maxsliprate-meanslip}, the fault aperture propagates with slip along the fault. This promotes pressure diffusion along the fault and results in a moderate pressure increase under the imposed constant injection flow rate. As a consequence, slip remains aseismic, with a maximum slip rate of order $10^{-5}$ m/s. However, the enhanced pressure diffusion also promotes pressure propagation outside the barrier, potentially leading to leakage.

In contrast, without shear dilation, the fault aperture is coupled only to the normal traction, which essentially follows the fault pressure profile (see Figure  \ref{fig:weak-nosd}). This results in a much smaller aperture and prevents significant pressure diffusion outside the barriers. Under the fixed injection flow rate, the reduced hydraulic transmissivity consequently leads to a much larger pressure increase. This, in turn, explains the occurrence of seismic events, characterized by maximum slip rates of order $1$ m/s, as shown in Figure  \ref{fig:weak-nosd} and in the right panel of Figure \ref{fig:maxsliprate-meanslip}. It is worth noting that we have verified that combining shear dilation with the rate-weakening RSF parameters still results in aseismic slip, with a maximum slip rate of up to $10^{-4}$ m/s.

Tables \ref{table_strengh_sd_mesh}, \ref{table_strengh_sd_dt}, \ref{table_weak_nosd_mesh}, and \ref{table_weak_nosd_dt} investigate the efficiency and robustness of the nonlinear and linear solvers for both configurations, considering the family of refined meshes and a range of time-step sizes parameterized by $M_{\Delta t} \in \{1,2,4,8,16\}$, as defined in \eqref{eq-timestepping}. The initial and maximum time steps are both set to $\Delta t^0 = \Delta t_{max} = 1000$ s. 

The scalability of both the iterative coupling algorithm and the Mechanics Newton solver with respect to mesh size is very good. A significant improvement is obtained through the probing stabilisation coefficient $C^{\rm fs}_f$, particularly in terms of the maximum number of fixed-point iterations required during the simulations. The preconditioned iterative linear solver for the mechanics problem also exhibits very good scalability, with only a mild dependence on mesh size, approximately of order $h^{-1/2}$, where $h$ denotes the fault face size. Overall, the simulation time scales very favorably with mesh refinement. 
As expected from the analysis, both the linear and nonlinear solvers are sensitive to the time-step size. Nevertheless, this dependence remains relatively moderate. Once again, the probing technique provides a significant improvement in the configuration with shear dilation.



\begin{table}[h!]
		\begin{center}
                  \begin{tabular}{ccccc||cccc}
		    \hline
                    	                             &&        & Probing  &       &  Fixed & Stress  &   \\\hline
				$\# \faces_\Gamma$    &&  677   & 1263    & 2432   & 677    & 1263    & 2432 \\\hline
				$N_{\Delta t}$         &&   3245 & 3186    & 3193   &  3245     &   3186       &   3193   \\
                    $Mean_{FPt}$                      &&   6     &   5.8   & 5.8     & 8.0      & 7.2        & 6.8     \\
                     $Max_{FPt}$                      &&   8     &   8     &  7      & 66      & 38        & 30      \\
				$N_{Newton}$          &&    10   &  9.8    &  9.7   &  12.3      & 10.7        & 10.1       \\
		                 $N_{GMRES}$          &&    13.7 &  18.3   &  25.2  &  13.7      & 18.3        & 25.3      \\
                    Run time (s)                          && 2.8k & 10k & 44k & 3.7k & 12k     &    48k   \\\hline
			\end{tabular}
			\caption{2D fluid injection test cases with rate-strengthening and shear dilation configuration simulated on the meshes with $\# \faces_\Gamma$ fault faces and with adaptive time stepping defined by $M_{\Delta t}=4$, $\Delta t_{max} = 10^3 =\Delta t^0 $ s. Number of successful time steps $N_{\Delta t}$,  average number of Newton iterations per time step $N_{Newton}$, average number of fixed-point iterations per time step $Mean_{FPt}$, maximum number of fixed-point iterations per time step $Max_{FPt}$, average number of GMRES  iterations per Newton step $N_{GMRES}$, run time for the overall simulation. }
			\label{table_strengh_sd_mesh}
		\end{center}
	\end{table}

\begin{table}[h!]
		\begin{center}
                  \begin{tabular}{ccccccc||cccccc}
		    \hline
                    	                        &&         &      & Probing   &   &     &    & Fixed    & Stress &                  &    \\\hline
				$M_{\Delta t}$    &&  1     & 2   & 4           & 8      & 16      &  1        & 2 & 4            & 8       & 16 \\\hline
				$N_{\Delta t}$    &&  1524  & 1965     &  3245  & 5884   &  11223  & 1550      & 1982   & 3245    &  5884     &  11223  \\
                    $Mean_{FPt}$                 &&  6.6   & 6.4    &  6       & 5.2    &   5.0    &  9.8      & 9.6  &  8        &  6.7      & 5.9   \\
                     $Max_{FPt}$                 &&  12    &  9   &  8         &  8     &   8     &  259      & 313  &  66        &  18      &  9  \\
				$N_{Newton}$     &&  12.6   & 12.0    &  10     &  8.1  &  7.0     &   21.2   & 18.9  &   12.3   &  9.0      &  7.2   \\
		                $N_{GMRES}$     &&   16.4   &  15.7   &  13.7   &  11.4 &  9.2     &  17.1    & 16.1  &  13.7    &  11.4     & 9.2    \\\hline
			\end{tabular}
			\caption{2D fluid injection test cases with rate-strengthening and shear dilation configuration simulated on the mesh with $677$ fault faces and with adaptive time stepping defined by $\Delta t_{max} = \Delta t^0=10^3$ s and $M_{\Delta t}=1,2,4,8,16$. Number of successful time steps $N_{\Delta t}$,  average number of Newton iterations per time step $N_{Newton}$, average number of fixed-point iterations per time step $Mean_{FPt}$, maximum number of fixed-point iterations per time step $Max_{FPt}$, average number of GMRES  iterations per Newton step $N_{GMRES}$. }
			\label{table_strengh_sd_dt}
		\end{center}
	\end{table}

\begin{table}[h!]
		\begin{center}
			\begin{tabular}{l*{7}{c}r}
				\hline
				$\# \faces_\Gamma$    &&  677   & 1263    & 2432   \\\hline
				$N_{\Delta t}$        &&  6550    & 6711  &  7139       \\
                                $Mean_{FPt}$      &&   6.1    & 6.1      &   6.1       \\
                                 $Max_{FPt}$      &&   12       &  12    &   12      \\
				$N_{Newton}$          &&   7.7    & 8.1  &  8.8       \\
				$N_{GMRES}$          &&   8.4    & 10.5  &  14.2        \\
                                Run time (s)         &&   3.6k    & 15k  &  70k       \\\hline
			\end{tabular}
			\caption{2D fluid injection test cases with rate-weakening and no shear dilation configuration simulated on the meshes with $\# \faces_\Gamma$ fault faces  with adaptive time stepping defined by $M_{\Delta t}=4$, $\Delta t_{max} = \Delta t^0=10^3$ s. Number of successful time steps $N_{\Delta t}$,  average number of Newton iterations per time step $N_{Newton}$, average number of fixed-point iterations per time step $Mean_{FPt}$, maximum number of fixed-point iterations per time step $Max_{FPt}$, average number of GMRES  iterations per Newton step $N_{GMRES}$, run time for the overall simulation. }
			\label{table_weak_nosd_mesh}
		\end{center}
	\end{table}

\begin{table}[h!]
		\begin{center}
			\begin{tabular}{l*{8}{c}r}
				\hline
				$M_{\Delta t}$    &&  1     & 2       & 4       & 8  & 16 \\\hline
				$N_{\Delta t}$    &&  2466   & 3907    &  6650   & 12018    &  23048  \\
                                $Mean_{FPt}$      &&  7.4   & 6.9     &   6.1   &  5.9  &  5.0  \\
                                $Max_{FPt}$       &&  168   &  150    &   12    &  8  &    8   \\
				$N_{Newton}$       && 11.5   &  10.6   &   7.7   & 6.8   &  5.4  \\
				$N_{GMRES}$       &&  14.8  &  13.1    &   8.4   & 6.7   &  5.7 \\\hline
			\end{tabular}
			\caption{2D fluid injection test cases with rate-weakening and no shear dilation configuration simulated on the mesh with $677$ fault faces and with adaptive time stepping defined by $M_{\Delta t}=1,2,4,8,16$, $\Delta t_{max} = \Delta t^0=10^3$ s. Number of successful time steps $N_{\Delta t}$,  average number of Newton iterations per time step $N_{Newton}$, average number of fixed-point iterations per time step $Mean_{FPt}$, maximum number of fixed-point iterations per time step $Max_{FPt}$, average number of GMRES  iterations per Newton step $N_{GMRES}$. }
			\label{table_weak_nosd_dt}
		\end{center}
	\end{table}

\subsection{Fluid injection 3D test case}\label{sec:testcoupled3D}
\begin{figure}[H]
  \centering
  \includegraphics[scale=0.25]{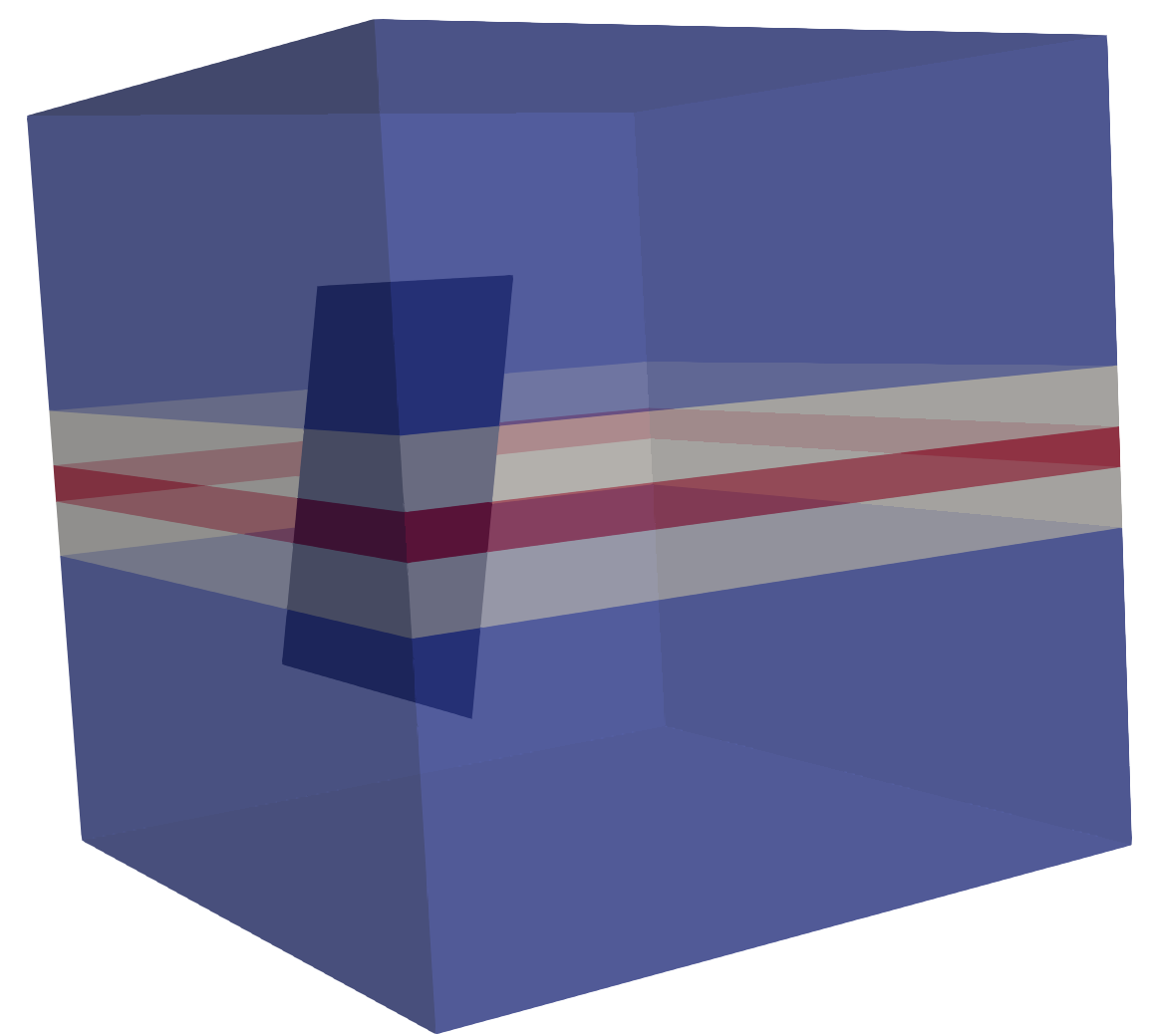}
  \includegraphics[scale=0.25]{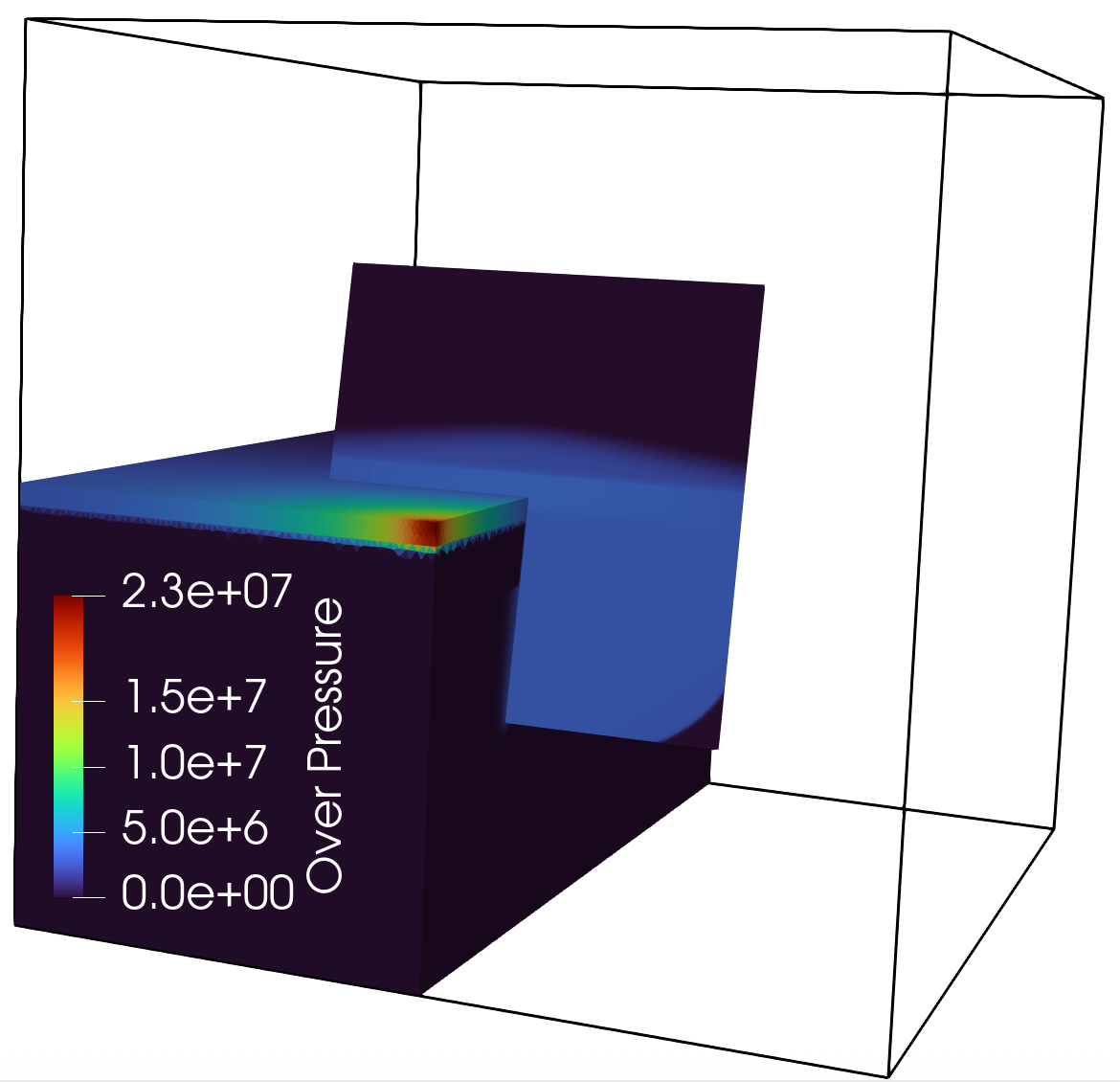}
  \caption{(left) Geometrical configuration of the 3D fluid injection test case with the simulation domain $(-1000\, {\rm m}, 1000\, {\rm m})^3$, the reservoir and barrier layers and the embedded 2D fault of extension $1000$ m in the $y$ direction.  (right) over pressure
    obtained with the $5.21M$ cells mesh and the time stepping defined by $\Delta t_{max} = 10^4$ s, $\Delta t^0 = 10^3$ s and $M_{\Delta t}= 4$.}
       \label{fig:testco23D}
\end{figure}

We extend the previous test case to the 3D geometrical configuration exhibited in Figure  \ref{fig:testco23D}. 
The domain is extruded in the $y$ direction over $2000$ m and the fault is embedded in the domain with an extension of $1000$ m in the $y$ direction. The test case reproduces the previous data set using the rate-strengthening RSF parameters and the shear dilation angle $\psi = \ang{5}$. The injection flow rate has the following Gaussian profile 
    $$
    q(y,z) = 2.68 ~10^{-2}~ e^{-{y^2 + z^2 \over L^2}} \quad \SI{}{kg.m^{-2}.s^{-1}}
    $$
    restricted to the left reservoir boundary and with $L=50$ m. The injection is stopped at time $10$ days and the simulation time is fixed to $T=15$ days.

\begin{figure}[H]
  \centering
  \includegraphics[scale=0.14]{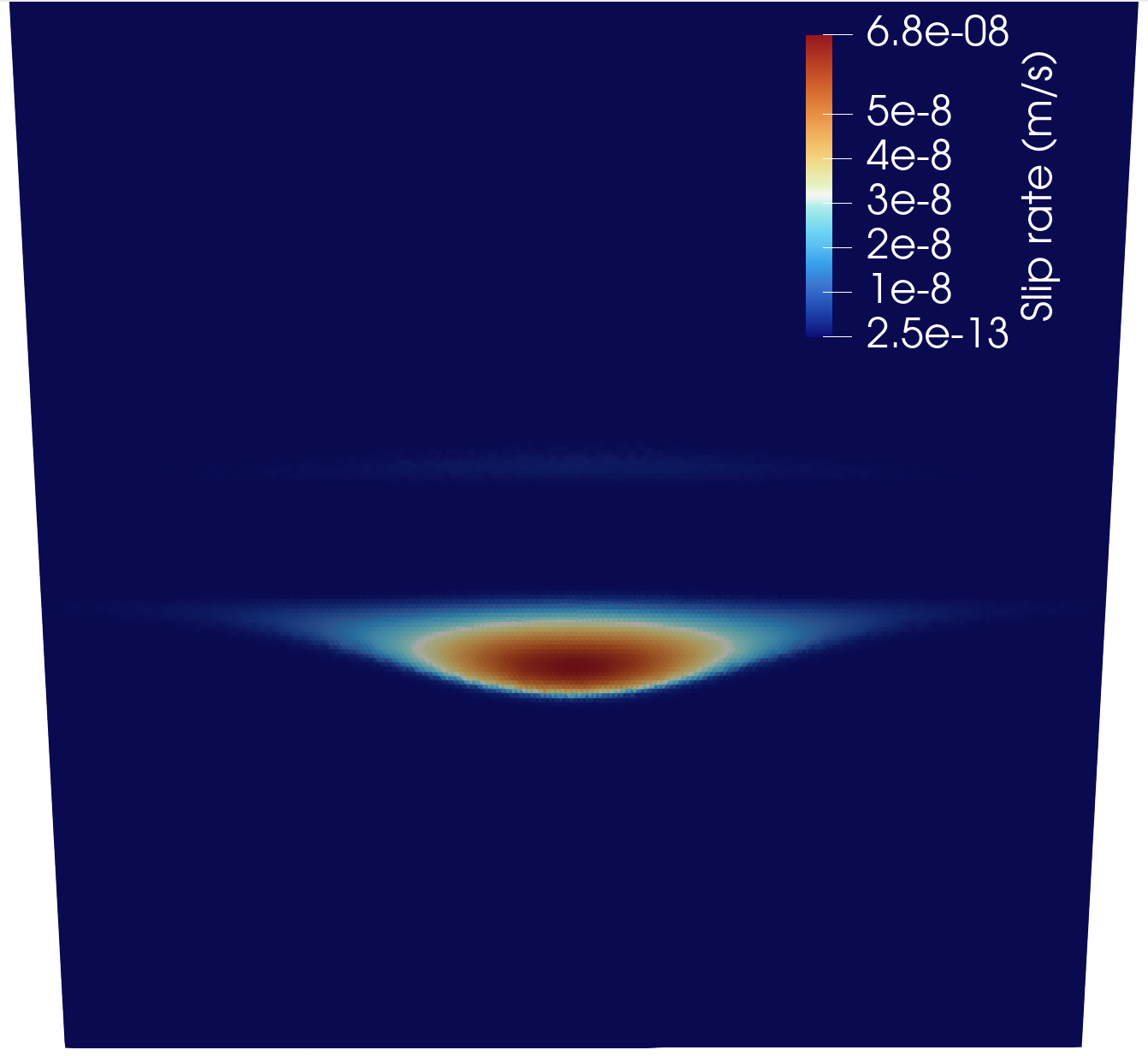}
  \includegraphics[scale=0.14]{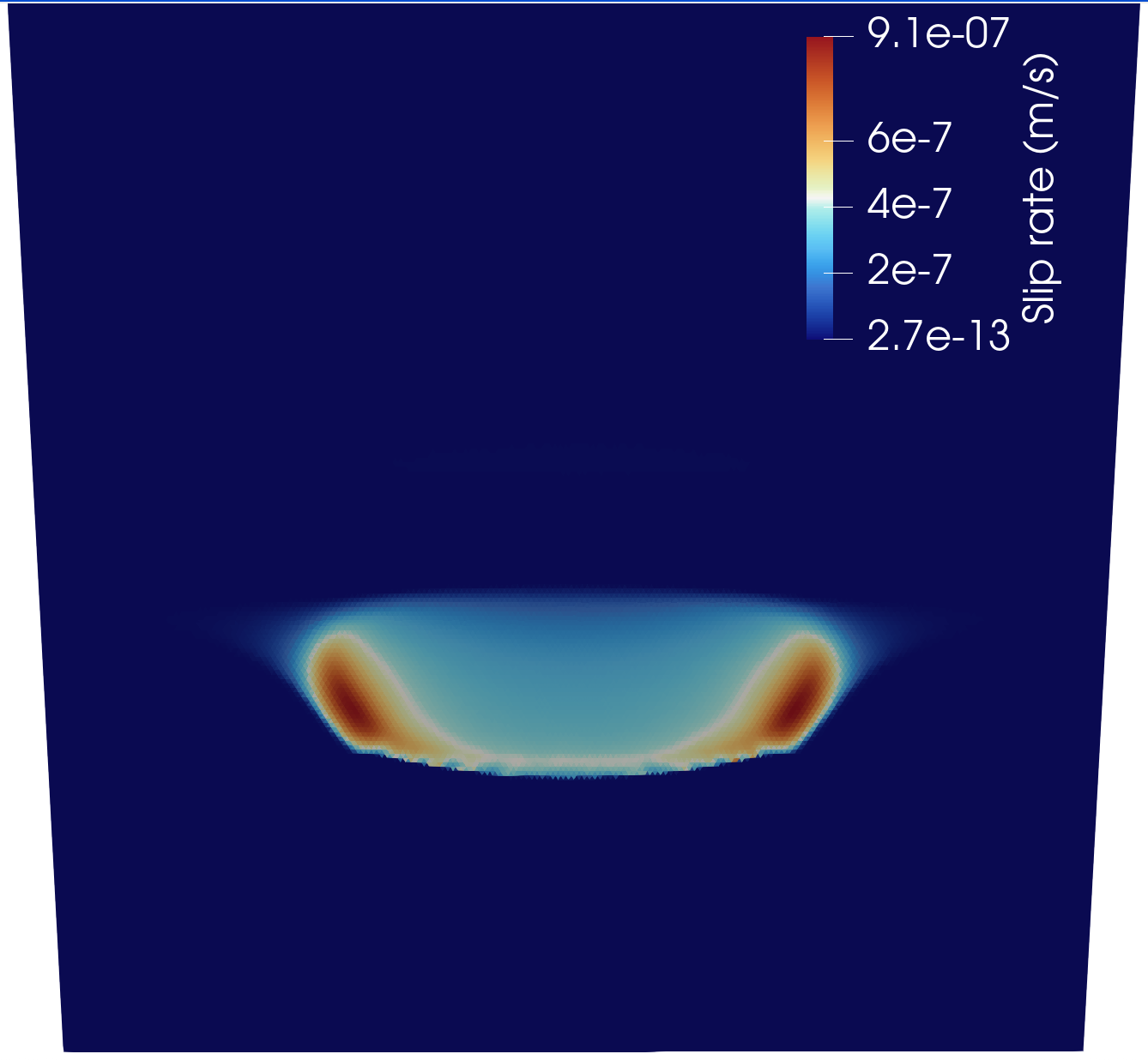}
  \includegraphics[scale=0.14]{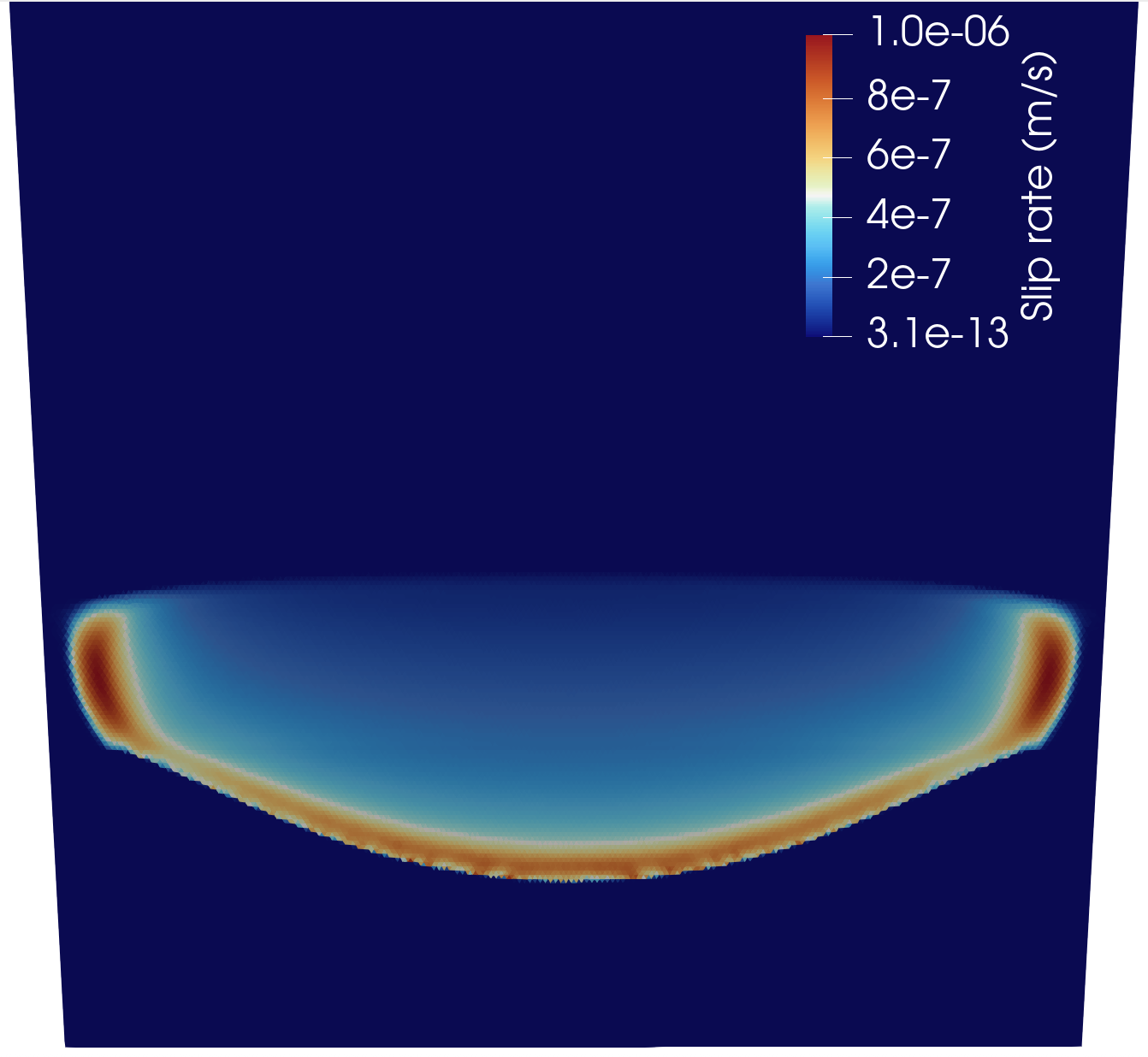}
  \includegraphics[scale=0.14]{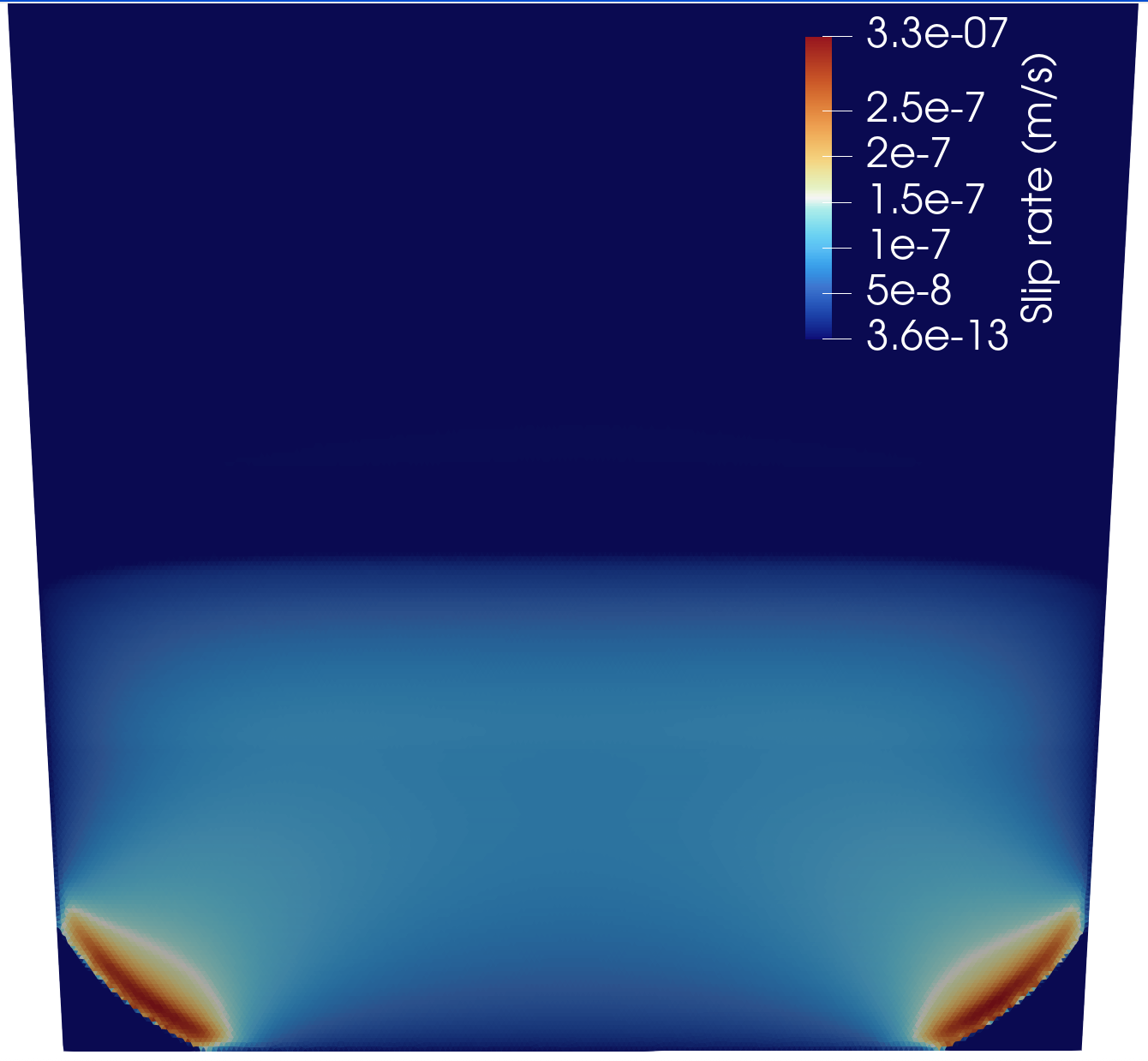}
  \includegraphics[scale=0.14]{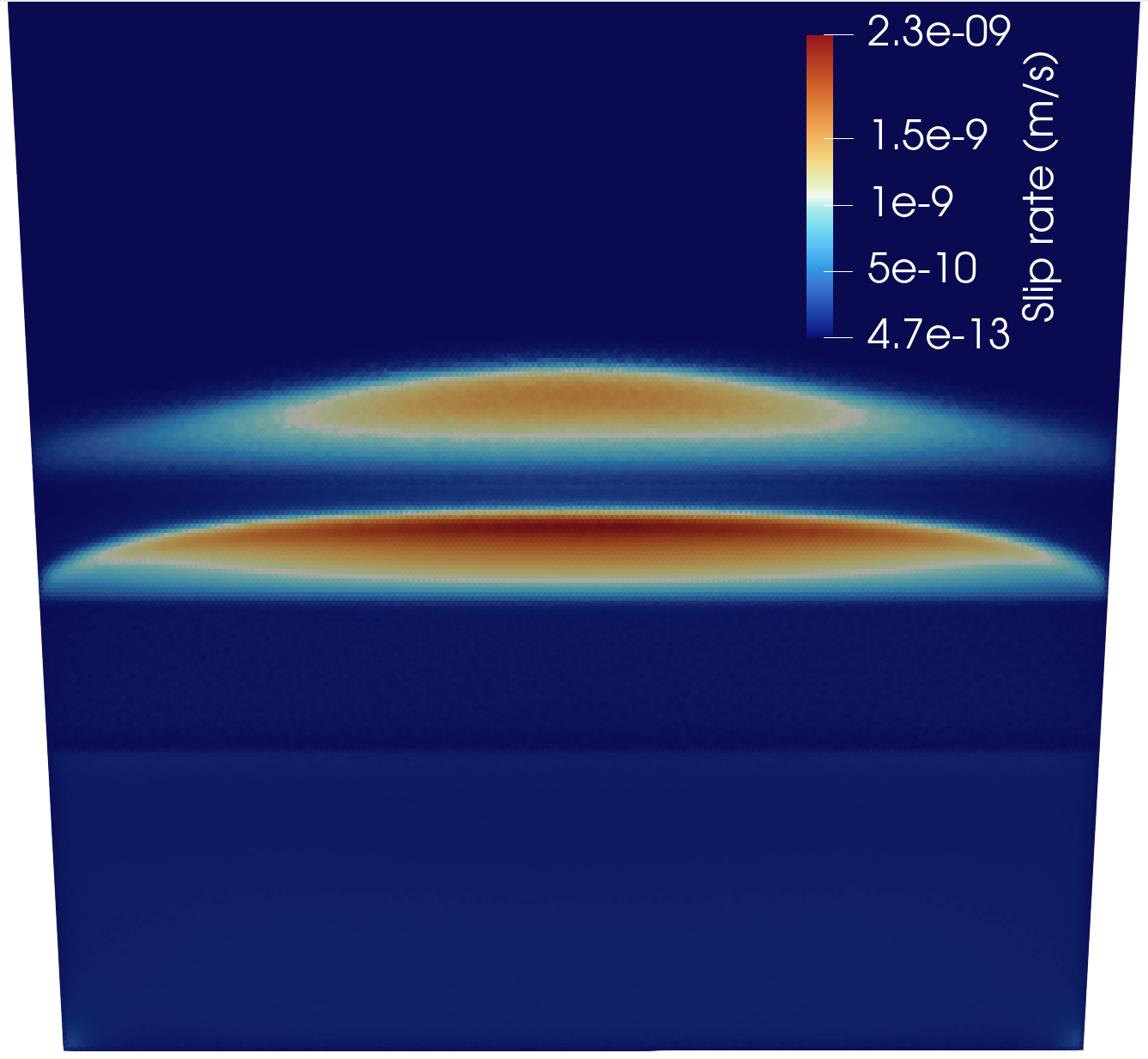}\\
   \includegraphics[scale=0.14]{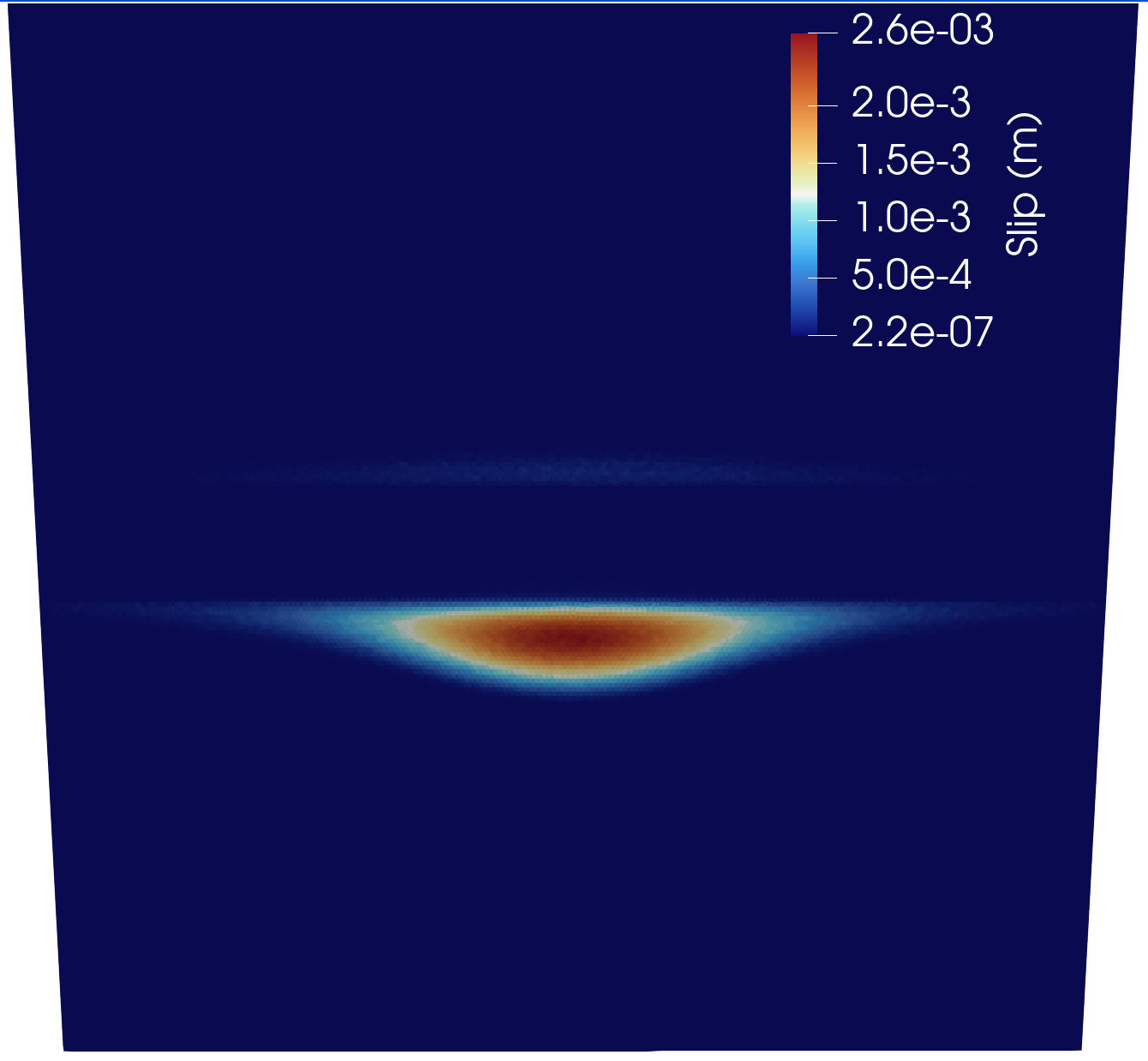}
  \includegraphics[scale=0.14]{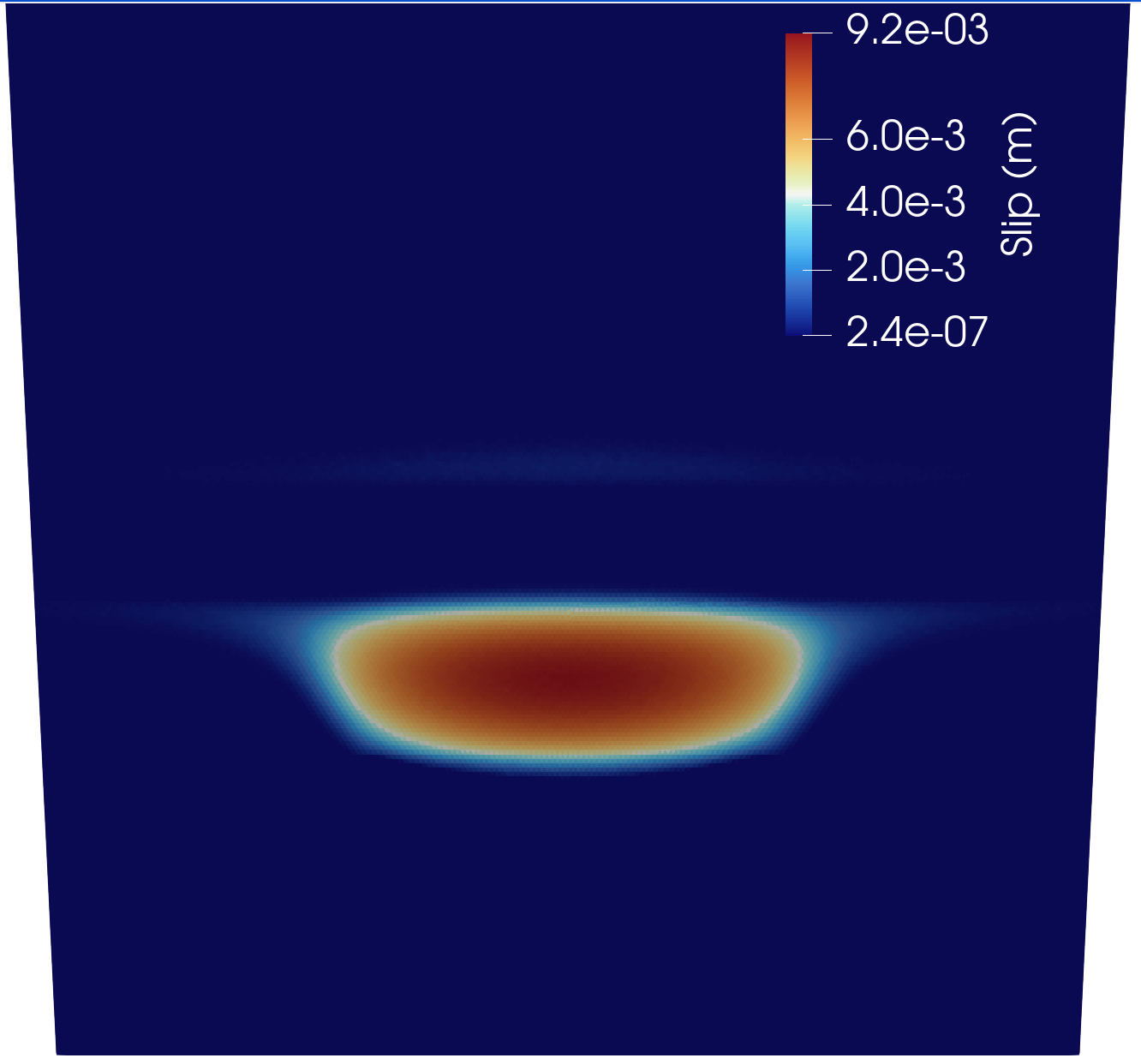}
  \includegraphics[scale=0.14]{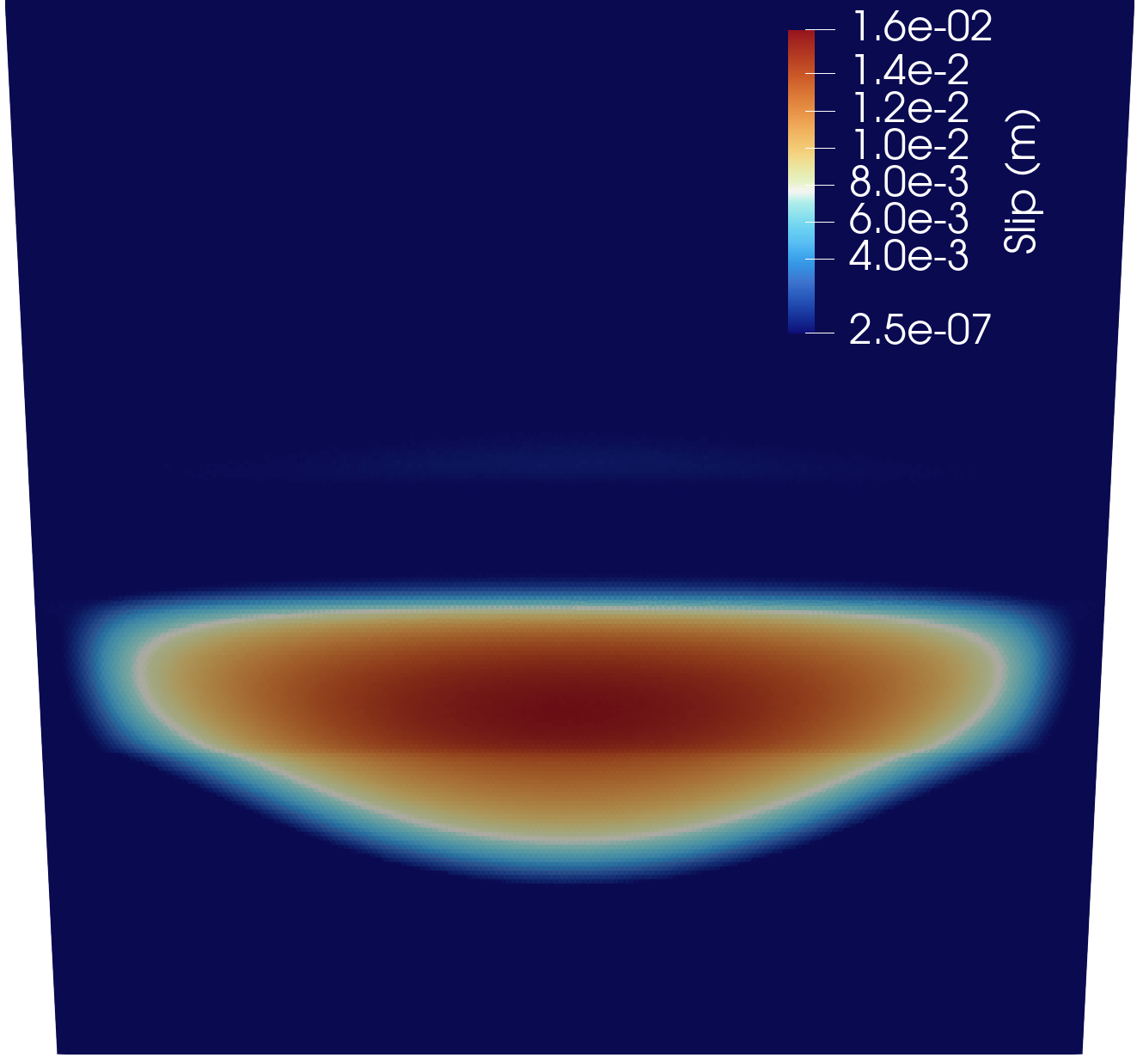}
  \includegraphics[scale=0.14]{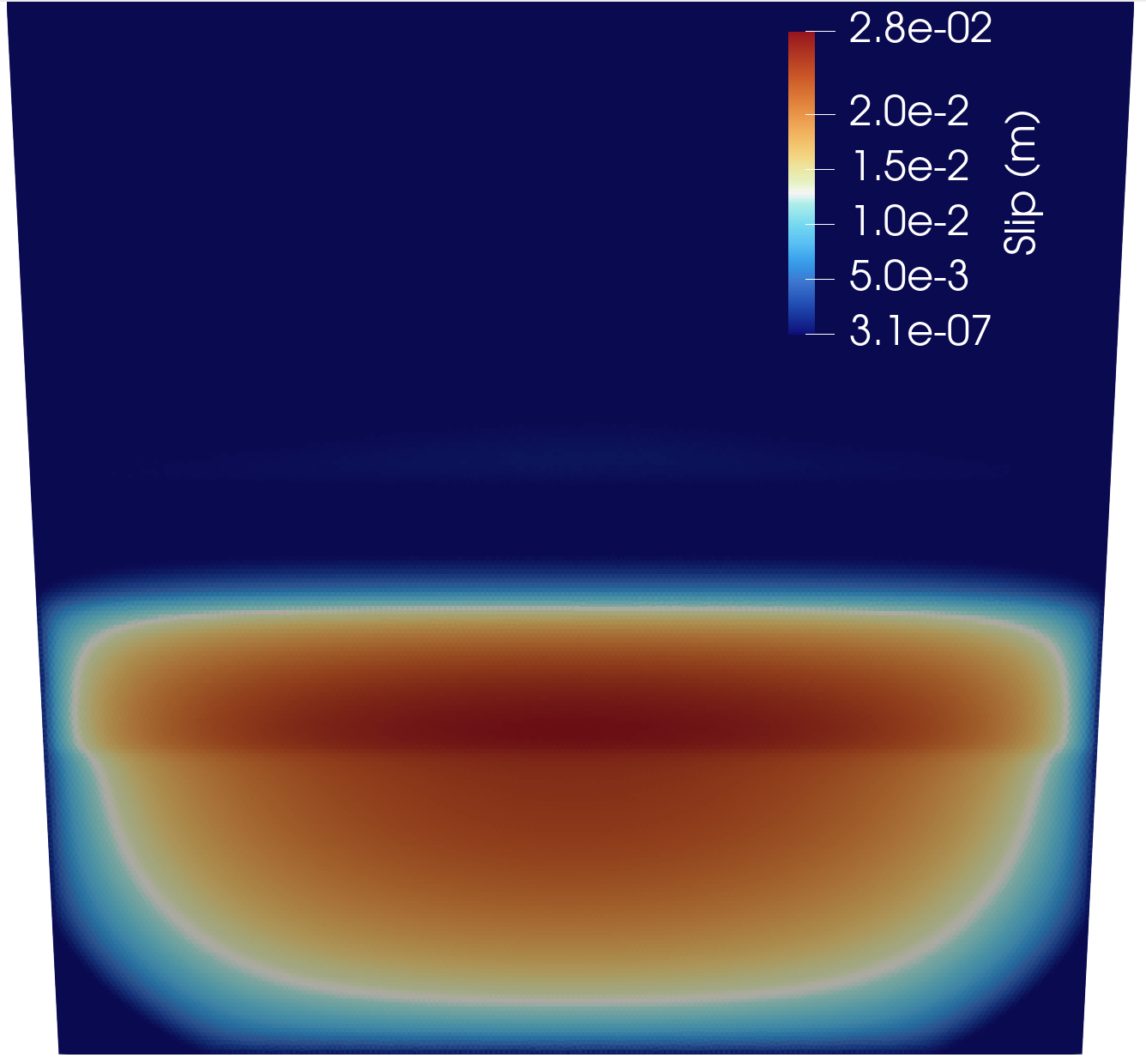}
  \includegraphics[scale=0.14]{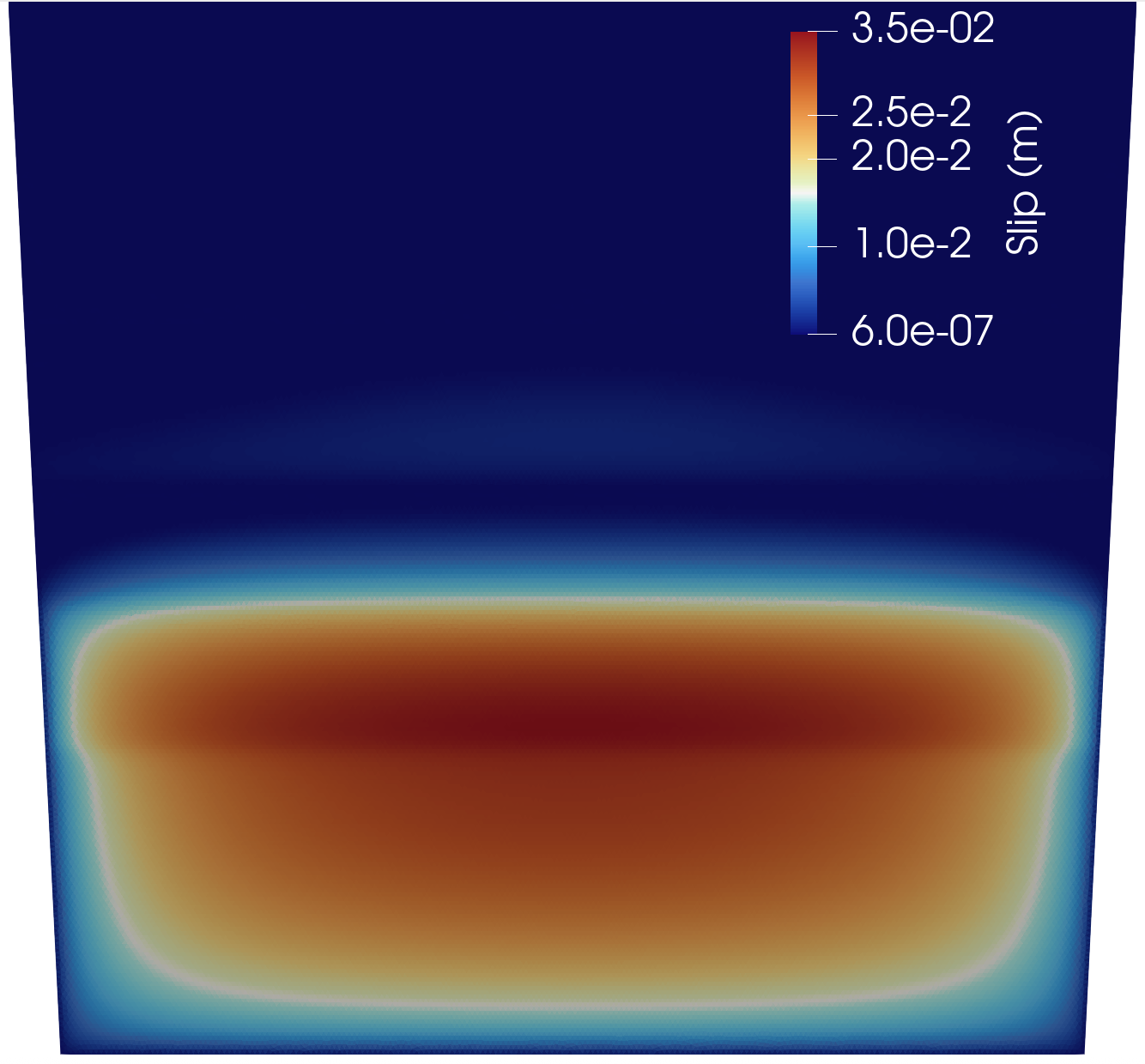}\\
  \includegraphics[scale=0.14]{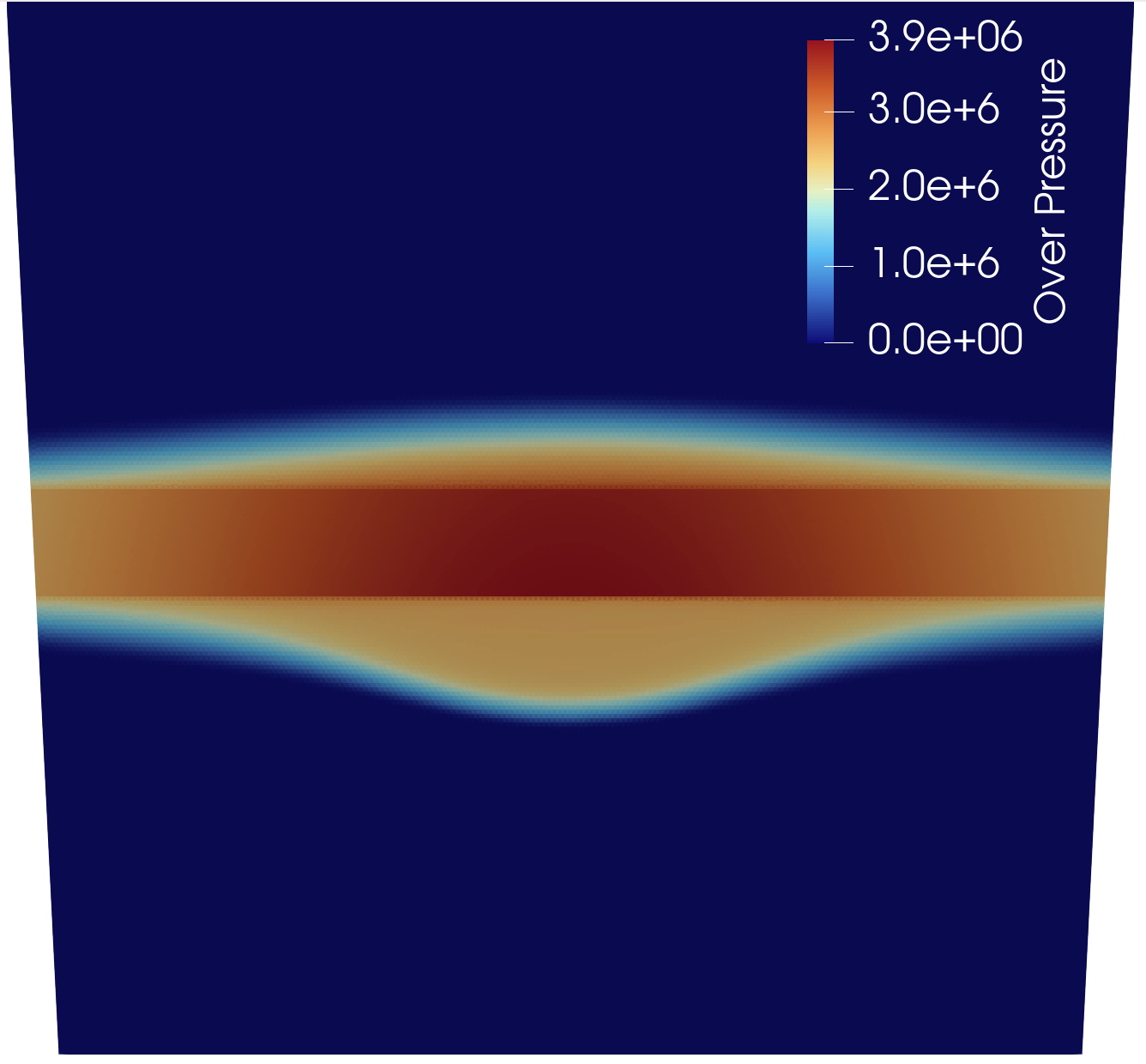}
  \includegraphics[scale=0.14]{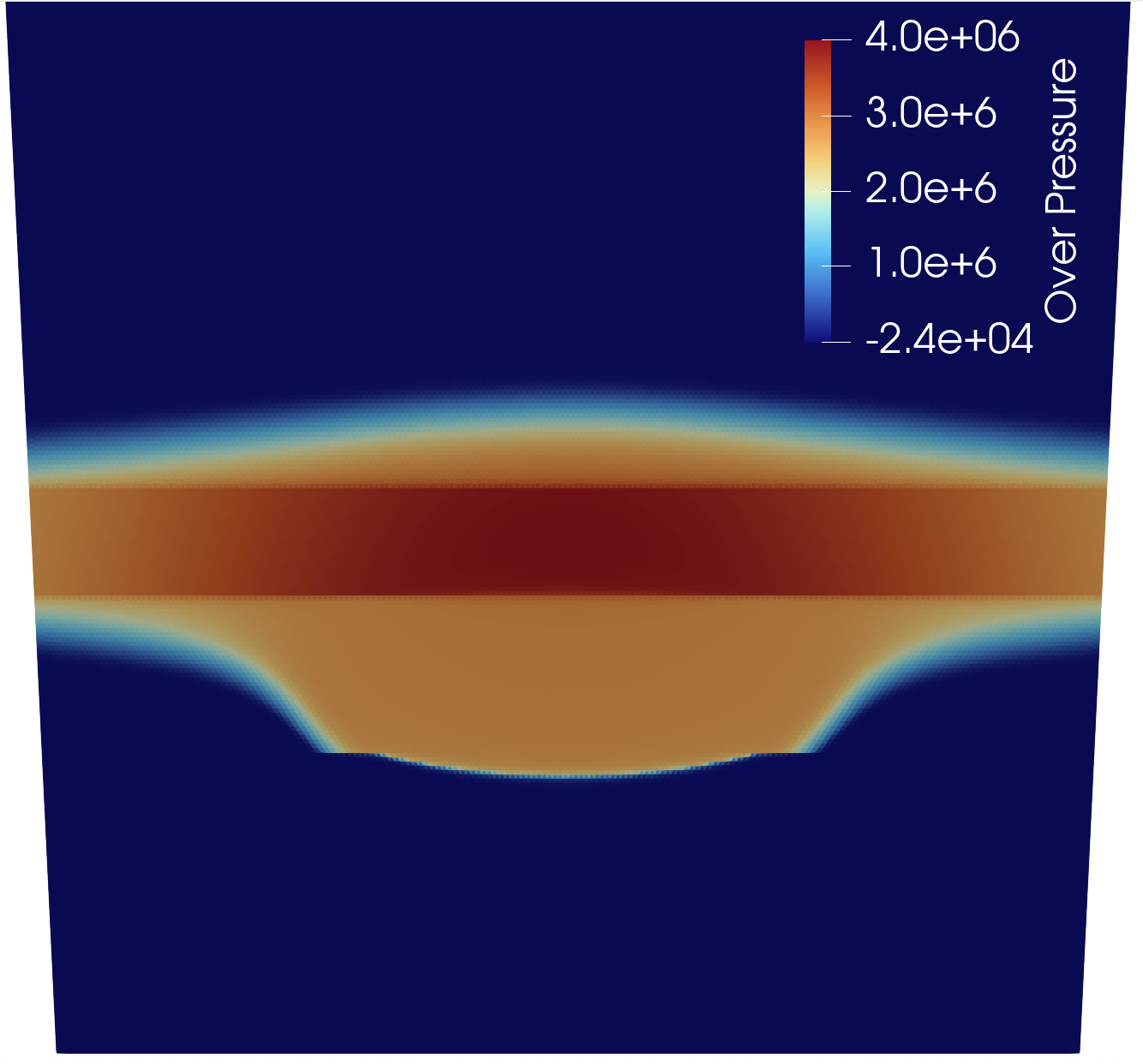}  
  \includegraphics[scale=0.14]{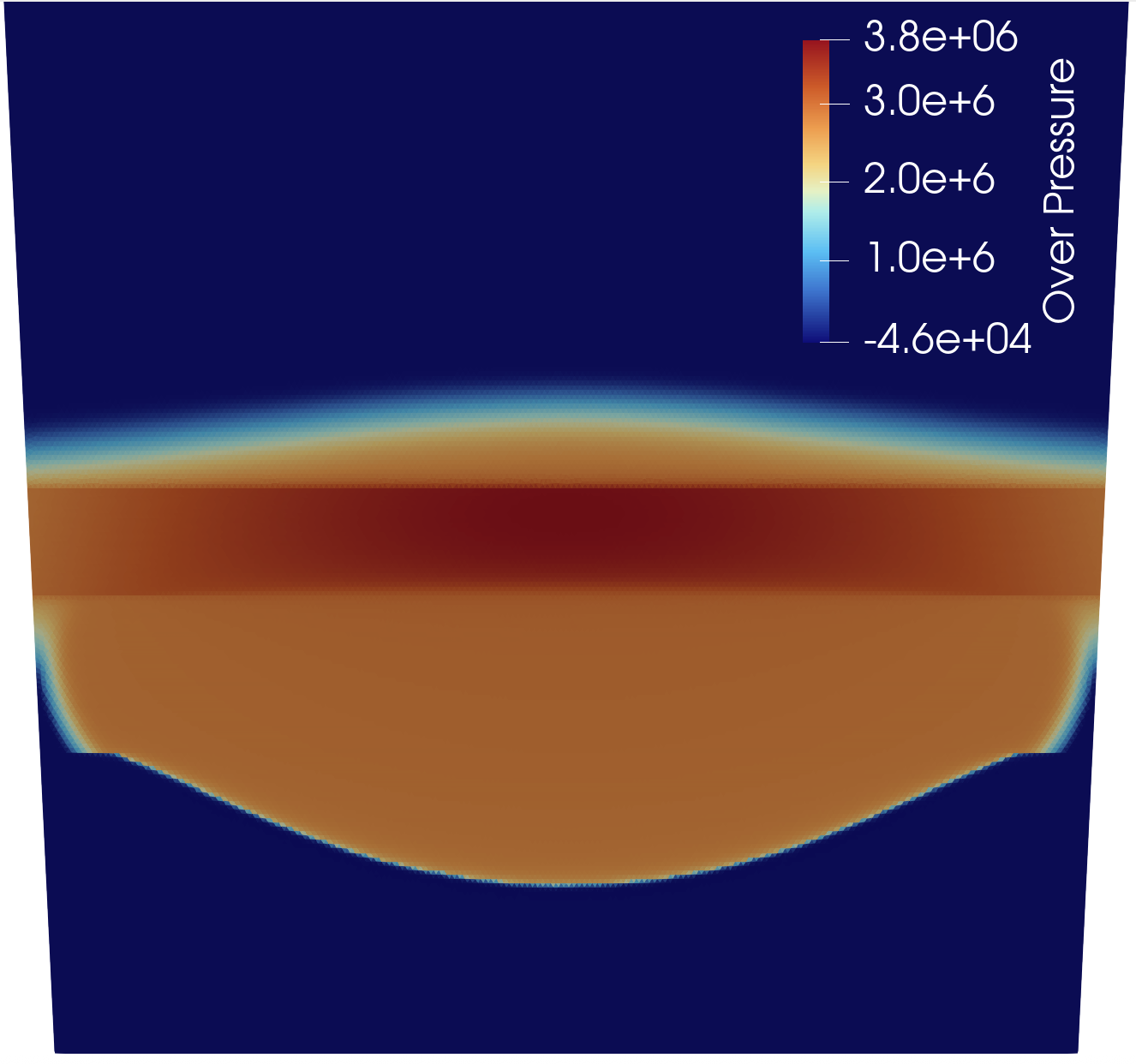}
  \includegraphics[scale=0.14]{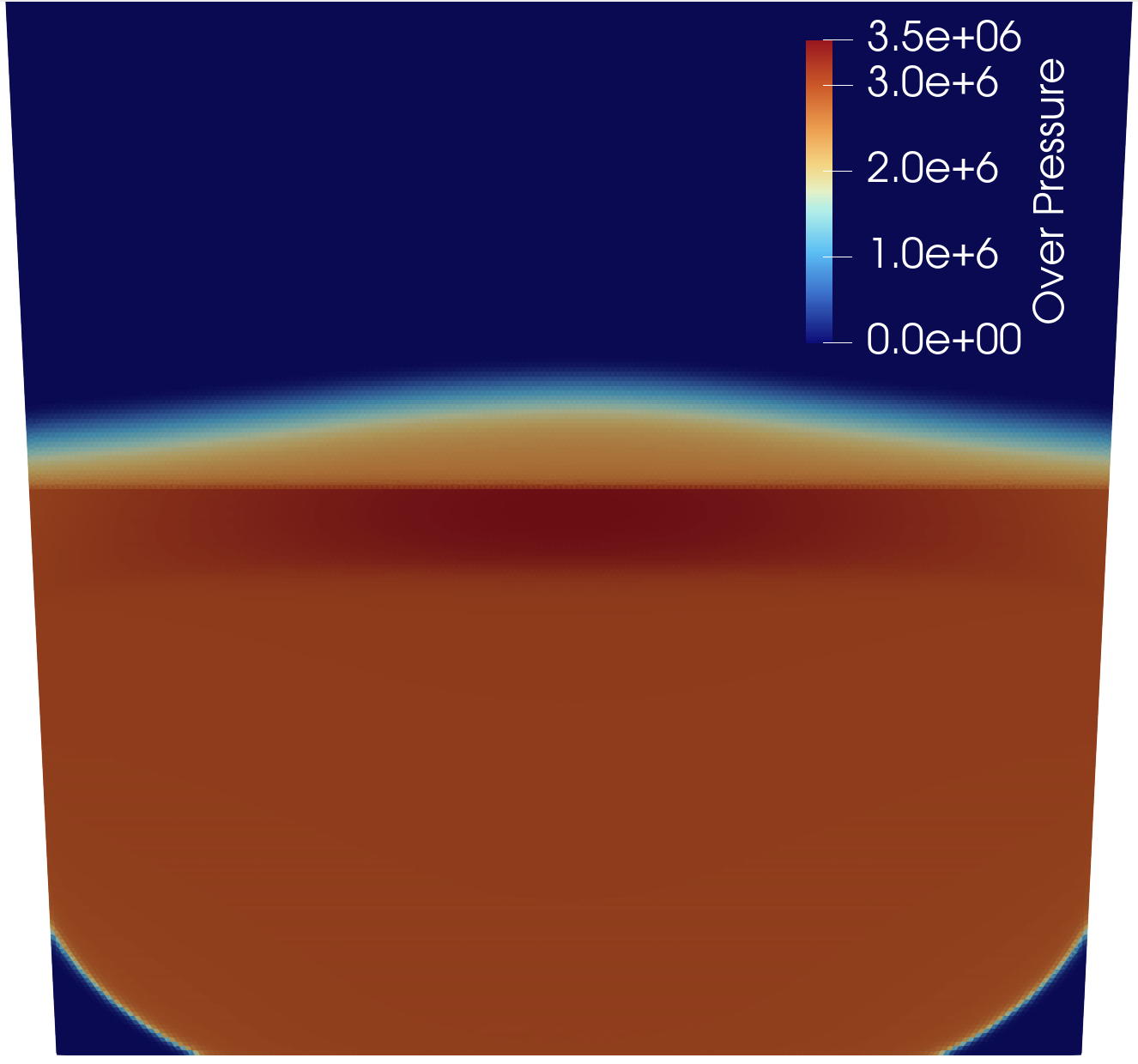}
  \includegraphics[scale=0.14]{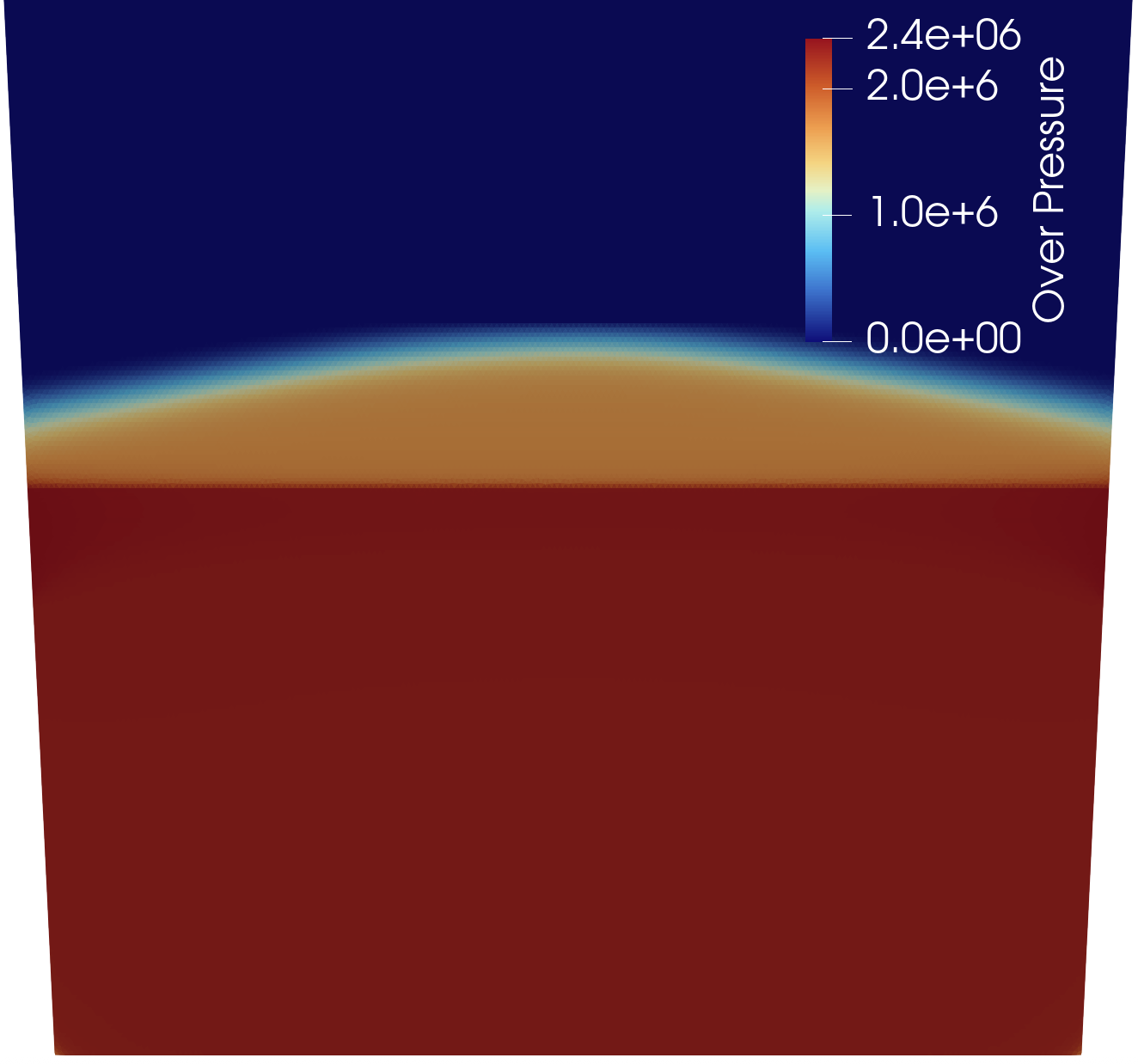}
\caption{3D fluid injection test case simulated on the $5.21M$ cells mesh and with the time stepping defined by $\Delta t_{max} = 10^4$ s, $\Delta t^0 = 10^3$ s and $M_{\Delta t}= 4$.  Slip rate (top line), slip (middle line) and over pressure (bottom line) along the fault at times $t=7.5, 8.2, 8.5, 10$ and $15$ days from left to right. }
       \label{fig:test3D-vspf}
\end{figure}

Figure \ref{fig:test3D-vspf} shows the slip rate, slip, and overpressure along the fault at successive times. The results show that slip initiates at the interface between the reservoir and the bottom barrier, as observed in the 2D case. The slip then propagates laterally along the bottom barrier and downward toward the bottom tip of the fault. As in the 2D case, the strong coupling between slip and shear dilation increases the hydraulic conductivity of the fault and enhances pressure diffusion, leading to a maximum overpressure of approximately 4 MPa and resulting in aseismic slip.

Table \ref{table_strengh_sd_3D} reports the performance of the simulation on a family of three successively refined meshes. We observe very good scalability of both the iterative coupling algorithm and the semi-smooth Newton solver for the contact-mechanics. The contact-mechanics preconditioner also exhibits only a weak sensitivity to mesh refinement. Note that the larger number of time steps on the coarsest mesh is due to the lower spatial resolution of the fault, which triggers oscillations in the maximum slip rate and, consequently, in the time-step size determined by \eqref{eq-timestepping}.    
    

\begin{table}[h!]
		\begin{center}
			\begin{tabular}{l*{7}{c}r}
				\hline
				$\# \cells$  / $\# \nodes$ / $\# \faces_\Gamma$   &&  1.32M/220k/37k       &   2.92M/483k/64k  & 5.21M/855k/94k  \\\hline
				$N_{\Delta t}$        &&  2086      &  1255 &   1122   \\
                                $Mean_{FPt}$         &&  7.0      &  7.0    &  6.9\\
                                 $Max_{FPt}$         &&  13       & 12    &  11   \\
				$N_{Newton}$         &&  12.7        & 12.7   & 12.7      \\
				$N_{GMRES}$          &&  19        &  22   &  24  \\
                                Run times (s)             &&   18k & 24k   &  34k   \\\hline
			\end{tabular}
			\caption{3D fluid injection test cases with rate-strengthening and  shear dilation configuration simulated on the meshes with $\# \cells$  cells, $\# \nodes$ nodes, $\# \faces_\Gamma$ fault faces,  with adaptive time stepping defined by $M_{\Delta t}=4$, $\Delta t_{max} = 10^4$ s, $\Delta t^0 = 10^3$ s. Number of successful time steps $N_{\Delta t}$,  average number of Newton iterations per time step $N_{Newton}$, average number of fixed-point iterations per time step $Mean_{FPt}$, maximum number of fixed-point iterations per time step $Max_{FPt}$, average number of GMRES  iterations per Newton step $N_{GMRES}$, run time for the overall simulation. }
			\label{table_strengh_sd_3D}
		\end{center}
	\end{table}

\section{Conclusion}
In this work, we have presented a fully coupled hydro-mechanical framework for the simulation of fluid-induced fault reactivation, with particular emphasis on the two-way coupling between fault slip and fluid flow. The proposed formulation combines mixed-dimensional Darcy flow in the porous matrix and along the fault with a quasi-dynamic contact-mechanics model governed by a rate-and-state friction law. A key feature of the approach is the explicit incorporation of shear dilation into the evolution of the fault hydraulic aperture, thereby accounting for the feedback between fault slip, hydraulic conductivity and pressure diffusion. This coupling is particularly relevant for subsurface applications involving fluid injection, such as geological carbon storage and underground energy storage.

From a numerical perspective, the framework combines a mixed Virtual Element Method for contact-mechanics with a Hybrid Finite Volume discretisation for the mixed-dimensional flow problem. The use of fully implicit time integration allows the time steps to be controlled primarily by the physical evolution of the system rather than by restrictive numerical stability constraints. The resulting nonlinear problem is solved through an iterative coupling strategy, which preserves the use of dedicated solvers for the mechanical and flow subproblems while providing convergence properties. The analysis establishes fixed-point convergence under a time-step condition independent of the spatial discretisation. In addition, the semi-smooth Newton method developed for the nonlinear contact problem, together with the proposed block preconditioning strategy, provides an efficient and robust treatment of the rate-and-state friction law. 

The numerical experiments demonstrate that the iterative coupling algorithm, the semi-smooth Newton solver, and the contact-mechanics preconditioner exhibit very good scalability with respect to mesh refinement and only moderate sensitivity to the time-step size. This favorable computational behavior constitutes an important step toward simulations of larger and more complex subsurface systems. They also demonstrate the ability of the framework to reproduce the complex interplay between fluid pressure, fault slip, shear dilation, and hydraulic conductivity. In particular, the coupled 2D and 3D fluid injection simulations clearly illustrate the strong influence of shear dilation on fault reactivation. When shear dilation is accounted for, slip enhances fault conductivity and promotes pressure diffusion, resulting in moderate pressure build-up and predominantly aseismic slip. Conversely, when shear dilation is neglected, the reduced hydraulic transmissivity leads to a much stronger pressure increase and to the development of seismic slip. These results highlight the importance of accounting for the two-way hydro-mechanical feedback in simulations of fluid-induced fault reactivation, rather than treating the hydraulic and mechanical responses as one-way coupled processes.

Overall, the results demonstrate that the proposed framework provides a consistent and computationally efficient approach for modeling strongly coupled hydro-mechanical fault processes and offers a promising numerical tool for investigating fault reactivation and assessing the consequences of fluid injection in subsurface energy applications.

\vskip 0.5cm 
\noindent{\bf Acknowledgements}: The authors gratefully acknowledge the financial support provided by the PEPR Maths-Vives program through the Mathsout project. They also warmly thank Fr\'ed\'eric Cappa from Geoazur for his valuable advice and insightful comments throughout the development of this work.

\section*{Appendix}

\begin{table}[!htbp]
    \centering
    \begin{tabular}{lll}
        \hline
        \textbf{Parameter} & \textbf{Value} & \textbf{Description} \\
        \hline
        \multicolumn{3}{l}{\textit{Porous medium and fault properties}} \\
        $k_1$ & $\SI{1e-15}{m^2}$ & Aquifer permeability \\
        $k_2$ & $\SI{1e-13}{m^2}$ & Reservoir permeability \\
        $k_3$ & $\SI{1e-17}{m^2}$ & Caprock permeability \\
        $\phi_m^0$ & $0.1$ & Initial porosity (aquifer) \\
        $\phi_m^0$ & $0.1$ & Initial porosity (reservoir) \\
        $\phi_m^0$ & $0.01$ & Initial porosity (caprock) \\
        $d_f^0$ & $\SI{1e-5}{m}$ & Initial hydraulic aperture \\
        $C_f$ & $d_f^3/12$ & Fault conductivity \\
        \hline
        \multicolumn{3}{l}{\textit{Mechanical and poromechanical properties}} \\
        $\alpha$ & $0.8$ & Biot coefficient \\
        $M$ & $\SI{10}{GPa}$ & Biot modulus \\
        $\lambda$ & $\SI{4}{GPa}$ & First Lam\'e coefficient \\
        $\mu$ & $\SI{4}{GPa}$ & Second Lam\'e coefficient \\
        $K_\n$ & $\SI{60}{GPa.m^{-1}}$ & Fault normal stiffness \\
        \hline
        \multicolumn{3}{l}{\textit{Fluid properties}} \\
        $c_f$ & $\SI{1e-10}{Pa^{-1}}$ & Fluid compressibility \\
        $\rho_{\mathrm{ref}}$ & $\SI{1000}{kg.m^{-3}}$ & Reference fluid density \\
        $\mu_f$ & $\SI{1e-3}{Pa.s}$ & Fluid dynamic viscosity \\
        \hline
        \multicolumn{3}{l}{\textit{Injection and hydraulic boundary conditions}} \\
        $q_{\mathrm{inj}}$ & $\SI{2e-3}{kg.m^{-2}.s^{-1}}$ & Injection flow rate \\
        $t_{\mathrm{inj}}$ & $\SI{10}{days}$ & Injection stopping time \\
        $T$ & $\SI{15}{days}$ & Final simulation time \\
        \hline
        \multicolumn{3}{l}{\textit{Mechanical boundary conditions}} \\
        $\rho_{\mathrm{rock/fluid}}$ & $\SI{2500}{kg.m^{-3}}$ & Average density for gravity loading \\
        $z_{\rm top}$ & $\SI{500}{m}$ & Depth of the top boundary \\
        $\sigma_h/\sigma_v$ & $0.6$ & Horizontal-to-vertical stress ratio \\
        \hline
        \multicolumn{3}{l}{\textit{Rate-and-state friction parameters}} \\
        $d_c$ & $\SI{5e-4}{m}$ & Characteristic slip distance \\
        $v_0$ & $\SI{1e-6}{m.s^{-1}}$ & Reference slip rate \\
        $f_0$ & $\tan(25^\circ)$ & Reference friction coefficient \\
        $c_s$ & $\SI{3464}{m.s^{-1}}$ & Shear-wave velocity \\
        $\eta$ & $\mu/(2c_s)$ & Radiation damping coefficient \\
        \hline
        \multicolumn{3}{l}{\textit{Test configurations}} \\
        $a$ & $0.016$ & Rate-strengthening parameter \\
        $b$ & $0.01$ & Rate-strengthening parameter \\
        $\psi$ & $5^\circ$ & Shear dilation angle \\
        \hline
        $a$ & $0.005$ & Rate-weakening parameter \\
        $b$ & $0.02$ & Rate-weakening parameter \\
        $\psi$ & $0^\circ$ & No shear dilation \\
        \hline
    \end{tabular}
    \caption{Parameters and data set used for the 2D fluid-induced fault reactivation test cases.}
    \label{tab:parameters_2D}
\end{table}
\newpage
%
\bibliographystyle{abbrv}
\bibliography{VB}

\end{document}